\documentclass[10pt,reqno]{amsart}
\usepackage[T1]{fontenc}
\usepackage[utf8]{inputenc}
\usepackage[sc]{mathpazo}
\usepackage{amsmath,amssymb,amsthm,mathrsfs,amsfonts}
\usepackage{esint}
\usepackage{cite}
\usepackage{tikz}
\usetikzlibrary{decorations.pathreplacing,patterns,arrows.meta,calc}
\usepackage{pgfplots}
\pgfplotsset{compat=1.18}
\usepgfplotslibrary{groupplots}
\definecolor{craneNavy}{HTML}{10184A}
\definecolor{craneInk}{HTML}{1B2A70}
\definecolor{craneGold}{HTML}{F5B400}
\definecolor{craneGoldDark}{HTML}{8A5A00}
\definecolor{craneGoldLight}{HTML}{FBEFD2}
\definecolor{craneGray}{HTML}{9CA0AC}
\definecolor{craneGrayLight}{HTML}{F1F2F5}
\usepackage{geometry}
\usepackage{enumitem}
\usepackage{xcolor}
\definecolor{azulA}{RGB}{0,102,204}
\definecolor{azulB}{RGB}{153,204,255}
\newcommand{\la}{\langle}
\newcommand{\ra}{\rangle}

\usepackage[colorlinks=true, linkcolor=blue, citecolor=blue, urlcolor=blue]{hyperref}
\usepackage[capitalize,nameinlink]{cleveref}
\theoremstyle{plain}
\newtheorem{thm}{Theorem}[section]
\newtheorem{prop}[thm]{Proposition}
\newtheorem{lem}[thm]{Lemma}
\newtheorem{cor}[thm]{Corollary}
\theoremstyle{definition}

\newtheorem{rem}[thm]{Remark}
\theoremstyle{remark}

\theoremstyle{plain}
\newtheorem{thmx}{Theorem}

\numberwithin{equation}{section}
\newcommand{\triple}{|\!|\!|}

\title[Decay and Regularity for the DGBO Equation]
{From decay to regularity: the dispersion generalized Benjamin-Ono equation}

\author[G. Fonseca]{Germ\'an Fonseca}
\address{Departamento de Matem\'aticas, Universidad Nacional de Colombia Carrera 45 No. 26-85, Edificio Uriel Guti\'errez Bogot\'a D.C., Colombia}
\curraddr{}
\email{gefonsecab@unal.edu.co}

\author[F.~Linares]{Felipe Linares}
\address{Instituto de Matemática Pura e Aplicada (IMPA), Estrada Dona Castorina, 110, 22460-320 Rio de Janeiro, RJ, Brazil}
\email{linares@impa.br}
\thanks{F. Linares was partially supported by CNPq grant 310329\!/\!2023-0 and FAPERJ grant E\--26\!/\!200236\!/\! 2026}

\author[A.\,J.~Méndez]{Argenis J.\ Méndez}
\address{IMECC-Instituto de Matemática, Estatística e Computação Científica, Universidade Estadual de Campinas (UNICAMP), Rua Sérgio Buarque de Holanda, 651, Cidade Universitária Barão Geraldo, 13083-859 Campinas, São Paulo, Brazil}
\email{ajmendez@unicamp.br}
\thanks{A.\,J.~Méndez gratefully acknowledges financial support from the Fundação de Amparo à Pesquisa do Estado de São Paulo (FAPESP), grant~2024/20513-7, and from the Conselho Nacional de Desenvolvimento Científico e Tecnológico (CNPq), through a \emph{Bolsa de Produtividade em Pesquisa}, grant~304088/2024-2}

\author[G.~Ponce]{Gustavo Ponce}
\address{Department of Mathematics,
	University of California-Santa Barbara,
	Santa Barbara, CA 923106, USA}
\email{ponce@math.ucsb.edu}
\thanks{G. Ponce was support by a grant from the Simons Foundation}

\author[O. Ria\~no]{Oscar Ria\~no}
\address{Departamento de Matem\'aticas, Universidad Nacional de Colombia, Ca\-rre\-ra 30 No. 45-03, 111321, Edificio Yu Takeuchi 404-209, Bogot\'a D.C., Colombia}
\email{ogrianoc@unal.edu.co}
\thanks{ O. Ria\~no gratefully acknowledges financial support from the Fundação de Amparo à Pesquisa do Estado de São Paulo (FAPESP), grant~2024/20513-7.}

\author[G. Rodriguez-Blanco]{Guillermo Rodr\'iguez-Blanco}
\address{Departamento de Matem\'aticas, Universidad Nacional de Colombia Carrera 45 No. 26-85, Edificio Uriel Guti\'errez Bogot\'a D.C., Colombia}
\curraddr{}
\email{grodriguezb@unal.edu.co}


\begin{document}
	
	
	\dedicatory{Dedicated to the memory of Prof. F\'elix Soriano, dear colleague and friend}
	
	\begin{abstract}
		We study the initial value problem associated with the dispersive generalized Benjamin-Ono equation 
		\begin{equation*}
			\partial_{t}u-\partial_{x}|D_{x}|^{a}u+u\,\partial_{x}u=0,\quad 1<a<2,
		\end{equation*}
		in weighted Sobolev spaces. First, we show that if the solution belongs to $L^{2}(|x|^{2\alpha}dx)$ at two different times, then it gains regularity, namely $u\in C([-T,T]:H^{a\alpha}(\mathbb{R}))$, for $\alpha<5/2+a$ and suitable initial regularity. Second, we establish the persistence of the solution flow in $H^{s}(\mathbb{R})\cap L^{2}(|x|^{2\alpha}dx)$ under the sharp regularity assumption $s\ge a\alpha$, for $\alpha\in[1/2,5/2+a)$, where the zero mean condition is required for $\alpha\ge 3/2+a$. Finally, we prove the unique continuation principles associated with the optimal decay rates $|x|^{3/2+a}$ and $|x|^{5/2+a}$ under optimal regularity assumptions. Our approach relies on weighted energy estimates, commutator expansions for the fractional Laplacian and pseudo-differential calculus.
	\end{abstract}
	
	\maketitle
	\section{Introduction}
	In this work we study the relation between decay and regularity of solutions to the initial value problem (IVP) associated with the \emph{dispersive generalized Benjamin-Ono (DGBO) equation}
	\begin{equation}\label{dgbo}
		\begin{cases}
			\partial_{t}u-\partial_{x}|D_{x}|^{a}u+u\,\partial_{x}u=0, & x\in\mathbb{R},\; t\in\mathbb{R},\; 1<a<2,\\
			u(x,0)=u_{0}(x).
		\end{cases}
	\end{equation}
	Here $u=u(x,t)$ is a real-valued function and, for $\beta>0$, $|D_{x}|^{\beta}$ denotes the \emph{fractional Laplacian} of order $\beta$ in the $x$-variable, that is,
	\begin{equation*}
		\widehat{|D_{x}|^{\beta}f}(\xi):=|\xi|^{\beta}\widehat{f}(\xi), \qquad f\in\mathcal{S}(\mathbb{R}).
	\end{equation*}
	We recall that $|D_{x}|=\mathcal{H}\partial_{x}$, where $\mathcal{H}$ denotes the Hilbert transform defined in \eqref{hita} below. The family \eqref{dgbo} interpolates between the \emph{Benjamin-Ono (BO) equation}
	\begin{equation}\label{BO}
		\partial_{t}u-\mathcal{H}\partial_{x}^{2}u+u\,\partial_{x}u=0,
	\end{equation}
	($a=1$), and the \emph{Korteweg-de Vries (KdV) equation}
	\begin{equation*}
		\partial_{t}u+\partial_{x}^{3}u+u\,\partial_{x}u=0,
	\end{equation*}
	($a=2$). These equations model vorticity waves in the coastal zone, see \cite{MST} and references therein.
	
	Real solutions of \eqref{dgbo} satisfy, at least formally, the following \emph{conservation laws}
	\begin{equation}\label{conservation}
		I[u](t)=\int_{\mathbb{R}}u\,dx, \qquad M[u](t)=\int_{\mathbb{R}}u^{2}\,dx, \qquad
		E[u](t)=\frac{1}{2}\int_{\mathbb{R}}\big||D_{x}|^{\frac{a}{2}}u\big|^{2}\,dx-\frac{1}{6}\int_{\mathbb{R}}u^{3}\,dx,
	\end{equation}
	and \eqref{dgbo} has the \emph{Hamiltonian structure} $\partial_{t}u=\partial_{x}\big(\frac{\delta E}{\delta u}\big)$. For $a\in(1,2)$ the equation in \eqref{dgbo} is not known to be completely integrable, and, to our knowledge, these are the only known conserved quantities. Also, if $u$ solves \eqref{dgbo}, then so does $u_{\lambda}(x,t)=\lambda^{a}u(\lambda x,\lambda^{a+1}t)$, $\lambda>0$, and
	\begin{equation*}
		\|u_{\lambda}(\cdot,0)\|_{\dot{H}^{s}}=\lambda^{a+s-\frac{1}{2}}\|u_{0}\|_{\dot{H}^{s}},
	\end{equation*}
	Thus, the \emph{scaling critical space} is $\dot{H}^{\frac{1}{2}-a}(\mathbb{R})$, and \eqref{dgbo} is sub-critical with respect to $M[u]$ and $E[u]$.

\begin{figure}[ht]
		\centering
		\begin{tikzpicture}
			\begin{groupplot}[
				group style={group size=2 by 1, horizontal sep=1.6cm},
				width=0.47\textwidth, height=5.4cm,
				axis lines=middle, axis line style={-{Stealth[length=2mm]}, craneNavy},
				xtick=\empty, ytick=\empty, xmin=-1.55, xmax=1.65,
				xlabel={$\xi$}, xlabel style={anchor=north west},
				every axis plot/.append style={line cap=round, samples=161},
				legend style={draw=none, fill=white, fill opacity=0.9, text opacity=1, font=\scriptsize, cells={anchor=west}}
				]
				\nextgroupplot[ymin=-4.2, ymax=4.2, title={\small $\omega(\xi)=-\xi|\xi|^{a}$}, title style={yshift=-2pt},
				legend style={at={(0.98,0.98)}, anchor=north east}]
				\addplot[craneGray, densely dashed, line width=1.1pt, domain=-1.5:1.5] {-x*abs(x)};
				\addlegendentry{$a=1$ (BO)}
				\addplot[craneGold, line width=1.4pt, domain=-1.5:1.5] {-x*abs(x)^1.5};
				\addlegendentry{$a=3/2$}
				\addplot[craneInk, line width=1.4pt, domain=-1.5:1.5] {-x^3};
				\addlegendentry{$a=2$ (KdV)}
				\nextgroupplot[ymin=-8.8, ymax=1.2, title={\small group velocity $\omega'(\xi)=-(a+1)|\xi|^{a}$}, title style={yshift=-2pt},
				legend style={at={(0.5,0.03)}, anchor=south}]
				\addplot[craneGray, densely dashed, line width=1.1pt, domain=-1.5:1.5] {-2*abs(x)};
				\addplot[craneGold, line width=1.4pt, domain=-1.5:1.5] {-2.5*abs(x)^1.5};
				\addplot[craneInk, line width=1.4pt, domain=-1.5:1.5] {-3*x^2};
			\end{groupplot}
		\end{tikzpicture}
		\caption{Dispersion relation \eqref{dispersionrelation} and group velocity for the BO equation ($a=1$), the DGBO equation with $a=3/2$, and the KdV equation ($a=2$). For $1<a<2$ the dispersion relation is $C^{2}$ but not $C^{3}$ at $\xi=0$, and all non-zero frequencies travel to the left.}
		\label{fig:dispersion}
	\end{figure}
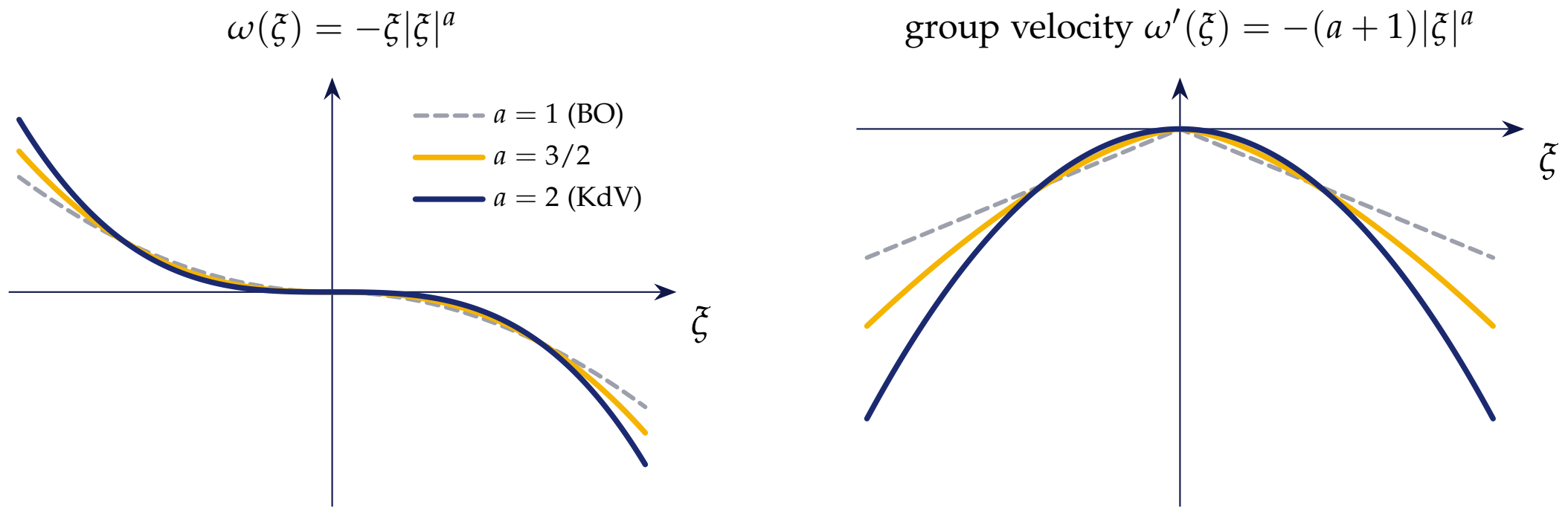

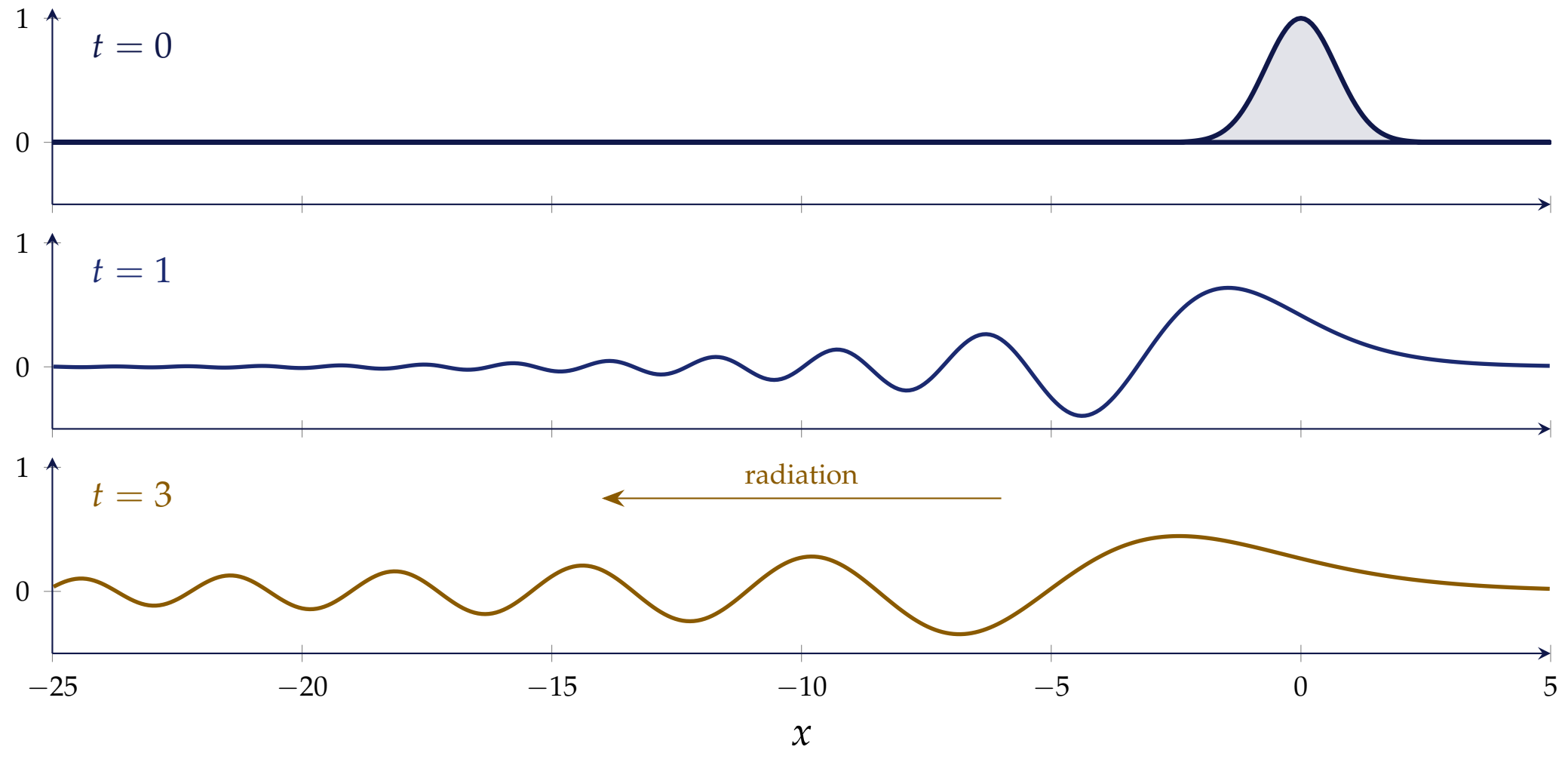
\begin{figure}[ht]
		\centering
		\begin{tikzpicture}
			\begin{groupplot}[
				group style={group size=1 by 3, vertical sep=0.25cm, xticklabels at=edge bottom},
				width=0.92\textwidth, height=3.3cm,
				xmin=-25, xmax=5, ymin=-0.5, ymax=1.08,
				axis lines=left, axis line style={craneNavy},
				ytick={0,1}, yticklabel style={font=\scriptsize},
				xtick={-25,-20,-15,-10,-5,0,5}, xticklabel style={font=\scriptsize},
				every axis plot/.append style={line join=round}
				]
				\nextgroupplot
				\addplot[craneNavy, line width=1.2pt, fill=craneNavy!12] table {
					-24.976 -0.0000
					-24.902 -0.0000
					-24.829 -0.0000
					-24.756 -0.0000
					-24.683 -0.0000
					-24.609 -0.0000
					-24.536 -0.0000
					-24.463 0.0000
					-24.390 -0.0000
					-24.316 -0.0000
					-24.243 -0.0000
					-24.170 0.0000
					-24.097 -0.0000
					-24.023 0.0000
					-23.950 -0.0000
					-23.877 0.0000
					-23.804 -0.0000
					-23.730 -0.0000
					-23.657 0.0000
					-23.584 0.0000
					-23.511 -0.0000
					-23.438 -0.0000
					-23.364 -0.0000
					-23.291 0.0000
					-23.218 0.0000
					-23.145 0.0000
					-23.071 -0.0000
					-22.998 -0.0000
					-22.925 0.0000
					-22.852 -0.0000
					-22.778 -0.0000
					-22.705 0.0000
					-22.632 0.0000
					-22.559 -0.0000
					-22.485 -0.0000
					-22.412 -0.0000
					-22.339 -0.0000
					-22.266 0.0000
					-22.192 -0.0000
					-22.119 0.0000
					-22.046 0.0000
					-21.973 -0.0000
					-21.899 -0.0000
					-21.826 -0.0000
					-21.753 -0.0000
					-21.680 0.0000
					-21.606 -0.0000
					-21.533 0.0000
					-21.460 0.0000
					-21.387 -0.0000
					-21.314 -0.0000
					-21.240 -0.0000
					-21.167 -0.0000
					-21.094 -0.0000
					-21.020 0.0000
					-20.947 0.0000
					-20.874 0.0000
					-20.801 -0.0000
					-20.727 -0.0000
					-20.654 0.0000
					-20.581 -0.0000
					-20.508 0.0000
					-20.435 0.0000
					-20.361 0.0000
					-20.288 -0.0000
					-20.215 -0.0000
					-20.142 0.0000
					-20.068 0.0000
					-19.995 -0.0000
					-19.922 -0.0000
					-19.849 -0.0000
					-19.775 0.0000
					-19.702 0.0000
					-19.629 -0.0000
					-19.556 -0.0000
					-19.482 0.0000
					-19.409 0.0000
					-19.336 0.0000
					-19.263 0.0000
					-19.189 0.0000
					-19.116 0.0000
					-19.043 0.0000
					-18.970 0.0000
					-18.896 -0.0000
					-18.823 0.0000
					-18.750 0.0000
					-18.677 0.0000
					-18.604 0.0000
					-18.530 0.0000
					-18.457 0.0000
					-18.384 0.0000
					-18.311 0.0000
					-18.237 -0.0000
					-18.164 0.0000
					-18.091 0.0000
					-18.018 0.0000
					-17.944 -0.0000
					-17.871 -0.0000
					-17.798 0.0000
					-17.725 0.0000
					-17.651 -0.0000
					-17.578 0.0000
					-17.505 0.0000
					-17.432 -0.0000
					-17.358 0.0000
					-17.285 0.0000
					-17.212 0.0000
					-17.139 0.0000
					-17.065 -0.0000
					-16.992 0.0000
					-16.919 0.0000
					-16.846 -0.0000
					-16.773 0.0000
					-16.699 0.0000
					-16.626 0.0000
					-16.553 -0.0000
					-16.480 -0.0000
					-16.406 -0.0000
					-16.333 0.0000
					-16.260 -0.0000
					-16.186 0.0000
					-16.113 0.0000
					-16.040 0.0000
					-15.967 -0.0000
					-15.894 0.0000
					-15.820 -0.0000
					-15.747 -0.0000
					-15.674 -0.0000
					-15.601 0.0000
					-15.527 0.0000
					-15.454 -0.0000
					-15.381 0.0000
					-15.308 -0.0000
					-15.234 -0.0000
					-15.161 -0.0000
					-15.088 -0.0000
					-15.015 0.0000
					-14.941 0.0000
					-14.868 -0.0000
					-14.795 -0.0000
					-14.722 0.0000
					-14.648 0.0000
					-14.575 -0.0000
					-14.502 0.0000
					-14.429 -0.0000
					-14.355 0.0000
					-14.282 0.0000
					-14.209 -0.0000
					-14.136 -0.0000
					-14.062 0.0000
					-13.989 -0.0000
					-13.916 -0.0000
					-13.843 0.0000
					-13.770 0.0000
					-13.696 0.0000
					-13.623 -0.0000
					-13.550 -0.0000
					-13.477 0.0000
					-13.403 -0.0000
					-13.330 0.0000
					-13.257 0.0000
					-13.184 0.0000
					-13.110 0.0000
					-13.037 0.0000
					-12.964 0.0000
					-12.891 -0.0000
					-12.817 0.0000
					-12.744 -0.0000
					-12.671 -0.0000
					-12.598 0.0000
					-12.524 0.0000
					-12.451 0.0000
					-12.378 0.0000
					-12.305 0.0000
					-12.231 0.0000
					-12.158 -0.0000
					-12.085 0.0000
					-12.012 -0.0000
					-11.938 -0.0000
					-11.865 -0.0000
					-11.792 -0.0000
					-11.719 -0.0000
					-11.646 -0.0000
					-11.572 -0.0000
					-11.499 0.0000
					-11.426 -0.0000
					-11.352 0.0000
					-11.279 0.0000
					-11.206 -0.0000
					-11.133 -0.0000
					-11.060 -0.0000
					-10.986 -0.0000
					-10.913 0.0000
					-10.840 -0.0000
					-10.767 -0.0000
					-10.693 -0.0000
					-10.620 -0.0000
					-10.547 -0.0000
					-10.474 -0.0000
					-10.400 -0.0000
					-10.327 0.0000
					-10.254 -0.0000
					-10.181 -0.0000
					-10.107 -0.0000
					-10.034 0.0000
					-9.961 0.0000
					-9.888 0.0000
					-9.815 0.0000
					-9.741 -0.0000
					-9.668 -0.0000
					-9.595 -0.0000
					-9.521 0.0000
					-9.448 0.0000
					-9.375 0.0000
					-9.302 0.0000
					-9.229 0.0000
					-9.155 0.0000
					-9.082 0.0000
					-9.009 0.0000
					-8.935 -0.0000
					-8.862 0.0000
					-8.789 0.0000
					-8.716 0.0000
					-8.643 -0.0000
					-8.569 -0.0000
					-8.496 -0.0000
					-8.423 -0.0000
					-8.350 0.0000
					-8.276 -0.0000
					-8.203 -0.0000
					-8.130 -0.0000
					-8.057 -0.0000
					-7.983 0.0000
					-7.910 -0.0000
					-7.837 -0.0000
					-7.764 0.0000
					-7.690 -0.0000
					-7.617 0.0000
					-7.544 0.0000
					-7.471 -0.0000
					-7.397 -0.0000
					-7.324 -0.0000
					-7.251 0.0000
					-7.178 0.0000
					-7.104 0.0000
					-7.031 0.0000
					-6.958 -0.0000
					-6.885 0.0000
					-6.811 -0.0000
					-6.738 -0.0000
					-6.665 0.0000
					-6.592 -0.0000
					-6.519 -0.0000
					-6.445 -0.0000
					-6.372 0.0000
					-6.299 0.0000
					-6.226 0.0000
					-6.152 0.0000
					-6.079 0.0000
					-6.006 0.0000
					-5.933 0.0000
					-5.859 0.0000
					-5.786 0.0000
					-5.713 0.0000
					-5.640 0.0000
					-5.566 0.0000
					-5.493 0.0000
					-5.420 0.0000
					-5.347 0.0000
					-5.273 0.0000
					-5.200 0.0000
					-5.127 0.0000
					-5.054 0.0000
					-4.981 0.0000
					-4.907 0.0000
					-4.834 0.0000
					-4.761 0.0000
					-4.688 0.0000
					-4.614 0.0000
					-4.541 0.0000
					-4.468 0.0000
					-4.394 0.0000
					-4.321 0.0000
					-4.248 0.0000
					-4.175 0.0000
					-4.102 0.0000
					-4.028 0.0000
					-3.955 0.0000
					-3.882 0.0000
					-3.809 0.0000
					-3.735 0.0000
					-3.662 0.0000
					-3.589 0.0000
					-3.516 0.0000
					-3.442 0.0000
					-3.369 0.0000
					-3.296 0.0000
					-3.223 0.0000
					-3.149 0.0001
					-3.076 0.0001
					-3.003 0.0001
					-2.930 0.0002
					-2.856 0.0003
					-2.783 0.0004
					-2.710 0.0006
					-2.637 0.0010
					-2.563 0.0014
					-2.490 0.0020
					-2.417 0.0029
					-2.344 0.0041
					-2.271 0.0058
					-2.197 0.0080
					-2.124 0.0110
					-2.051 0.0149
					-1.978 0.0200
					-1.904 0.0266
					-1.831 0.0350
					-1.758 0.0455
					-1.685 0.0586
					-1.611 0.0745
					-1.538 0.0939
					-1.465 0.1170
					-1.392 0.1442
					-1.318 0.1759
					-1.245 0.2122
					-1.172 0.2533
					-1.099 0.2991
					-1.025 0.3494
					-0.952 0.4039
					-0.879 0.4619
					-0.806 0.5225
					-0.732 0.5848
					-0.659 0.6476
					-0.586 0.7094
					-0.513 0.7689
					-0.440 0.8244
					-0.366 0.8745
					-0.293 0.9177
					-0.220 0.9529
					-0.146 0.9788
					-0.073 0.9947
					0.000 1.0000
					0.073 0.9947
					0.146 0.9788
					0.220 0.9529
					0.293 0.9177
					0.366 0.8745
					0.440 0.8244
					0.513 0.7689
					0.586 0.7094
					0.659 0.6476
					0.732 0.5848
					0.806 0.5225
					0.879 0.4619
					0.952 0.4039
					1.025 0.3494
					1.099 0.2991
					1.172 0.2533
					1.245 0.2122
					1.318 0.1759
					1.392 0.1442
					1.465 0.1170
					1.538 0.0939
					1.611 0.0745
					1.685 0.0586
					1.758 0.0455
					1.831 0.0350
					1.904 0.0266
					1.978 0.0200
					2.051 0.0149
					2.124 0.0110
					2.197 0.0080
					2.271 0.0058
					2.344 0.0041
					2.417 0.0029
					2.490 0.0020
					2.563 0.0014
					2.637 0.0010
					2.710 0.0006
					2.783 0.0004
					2.856 0.0003
					2.930 0.0002
					3.003 0.0001
					3.076 0.0001
					3.149 0.0001
					3.223 0.0000
					3.296 0.0000
					3.369 0.0000
					3.442 0.0000
					3.516 0.0000
					3.589 0.0000
					3.662 0.0000
					3.735 0.0000
					3.809 0.0000
					3.882 0.0000
					3.955 0.0000
					4.028 0.0000
					4.102 0.0000
					4.175 0.0000
					4.248 0.0000
					4.321 0.0000
					4.394 0.0000
					4.468 0.0000
					4.541 0.0000
					4.614 0.0000
					4.688 0.0000
					4.761 0.0000
					4.834 0.0000
					4.907 0.0000
					4.981 0.0000
				} \closedcycle;
				\node[anchor=north west, font=\small, craneNavy] at (axis cs:-24.5,1.0) {$t=0$};
				\nextgroupplot
				\addplot[craneInk, line width=1.0pt] table {
					-24.976 0.0018
					-24.902 0.0012
					-24.829 0.0004
					-24.756 -0.0004
					-24.683 -0.0012
					-24.609 -0.0019
					-24.536 -0.0024
					-24.463 -0.0027
					-24.390 -0.0027
					-24.316 -0.0024
					-24.243 -0.0018
					-24.170 -0.0009
					-24.097 0.0000
					-24.023 0.0010
					-23.950 0.0019
					-23.877 0.0027
					-23.804 0.0032
					-23.730 0.0033
					-23.657 0.0031
					-23.584 0.0026
					-23.511 0.0017
					-23.438 0.0006
					-23.364 -0.0006
					-23.291 -0.0017
					-23.218 -0.0028
					-23.145 -0.0036
					-23.071 -0.0040
					-22.998 -0.0041
					-22.925 -0.0037
					-22.852 -0.0030
					-22.778 -0.0019
					-22.705 -0.0005
					-22.632 0.0009
					-22.559 0.0023
					-22.485 0.0036
					-22.412 0.0045
					-22.339 0.0050
					-22.266 0.0050
					-22.192 0.0045
					-22.119 0.0035
					-22.046 0.0022
					-21.973 0.0006
					-21.899 -0.0012
					-21.826 -0.0029
					-21.753 -0.0043
					-21.680 -0.0055
					-21.606 -0.0061
					-21.533 -0.0062
					-21.460 -0.0057
					-21.387 -0.0046
					-21.314 -0.0031
					-21.240 -0.0012
					-21.167 0.0009
					-21.094 0.0029
					-21.020 0.0048
					-20.947 0.0063
					-20.874 0.0073
					-20.801 0.0076
					-20.727 0.0073
					-20.654 0.0063
					-20.581 0.0046
					-20.508 0.0025
					-20.435 0.0001
					-20.361 -0.0024
					-20.288 -0.0049
					-20.215 -0.0070
					-20.142 -0.0085
					-20.068 -0.0094
					-19.995 -0.0095
					-19.922 -0.0088
					-19.849 -0.0073
					-19.775 -0.0051
					-19.702 -0.0024
					-19.629 0.0006
					-19.556 0.0036
					-19.482 0.0065
					-19.409 0.0089
					-19.336 0.0107
					-19.263 0.0117
					-19.189 0.0117
					-19.116 0.0107
					-19.043 0.0089
					-18.970 0.0062
					-18.896 0.0030
					-18.823 -0.0006
					-18.750 -0.0043
					-18.677 -0.0079
					-18.604 -0.0109
					-18.530 -0.0132
					-18.457 -0.0145
					-18.384 -0.0148
					-18.311 -0.0139
					-18.237 -0.0120
					-18.164 -0.0090
					-18.091 -0.0053
					-18.018 -0.0010
					-17.944 0.0035
					-17.871 0.0079
					-17.798 0.0119
					-17.725 0.0151
					-17.651 0.0174
					-17.578 0.0184
					-17.505 0.0181
					-17.432 0.0165
					-17.358 0.0137
					-17.285 0.0097
					-17.212 0.0049
					-17.139 -0.0004
					-17.065 -0.0059
					-16.992 -0.0111
					-16.919 -0.0158
					-16.846 -0.0196
					-16.773 -0.0222
					-16.699 -0.0233
					-16.626 -0.0229
					-16.553 -0.0209
					-16.480 -0.0175
					-16.406 -0.0128
					-16.333 -0.0070
					-16.260 -0.0006
					-16.186 0.0060
					-16.113 0.0126
					-16.040 0.0185
					-15.967 0.0234
					-15.894 0.0270
					-15.820 0.0290
					-15.747 0.0292
					-15.674 0.0276
					-15.601 0.0242
					-15.527 0.0192
					-15.454 0.0128
					-15.381 0.0053
					-15.308 -0.0026
					-15.234 -0.0108
					-15.161 -0.0185
					-15.088 -0.0255
					-15.015 -0.0311
					-14.941 -0.0351
					-14.868 -0.0372
					-14.795 -0.0372
					-14.722 -0.0351
					-14.648 -0.0308
					-14.575 -0.0247
					-14.502 -0.0170
					-14.429 -0.0080
					-14.355 0.0017
					-14.282 0.0116
					-14.209 0.0212
					-14.136 0.0300
					-14.062 0.0374
					-13.989 0.0430
					-13.916 0.0466
					-13.843 0.0477
					-13.770 0.0463
					-13.696 0.0425
					-13.623 0.0362
					-13.550 0.0279
					-13.477 0.0179
					-13.403 0.0066
					-13.330 -0.0054
					-13.257 -0.0175
					-13.184 -0.0291
					-13.110 -0.0397
					-13.037 -0.0487
					-12.964 -0.0556
					-12.891 -0.0601
					-12.817 -0.0617
					-12.744 -0.0605
					-12.671 -0.0564
					-12.598 -0.0495
					-12.524 -0.0401
					-12.451 -0.0286
					-12.378 -0.0153
					-12.305 -0.0010
					-12.231 0.0139
					-12.158 0.0286
					-12.085 0.0424
					-12.012 0.0548
					-11.938 0.0652
					-11.865 0.0731
					-11.792 0.0780
					-11.719 0.0798
					-11.646 0.0782
					-11.572 0.0733
					-11.499 0.0651
					-11.426 0.0540
					-11.352 0.0403
					-11.279 0.0246
					-11.206 0.0073
					-11.133 -0.0107
					-11.060 -0.0289
					-10.986 -0.0465
					-10.913 -0.0629
					-10.840 -0.0773
					-10.767 -0.0892
					-10.693 -0.0981
					-10.620 -0.1036
					-10.547 -0.1055
					-10.474 -0.1035
					-10.400 -0.0976
					-10.327 -0.0881
					-10.254 -0.0752
					-10.181 -0.0592
					-10.107 -0.0407
					-10.034 -0.0203
					-9.961 0.0014
					-9.888 0.0237
					-9.815 0.0458
					-9.741 0.0670
					-9.668 0.0866
					-9.595 0.1038
					-9.521 0.1182
					-9.448 0.1292
					-9.375 0.1363
					-9.302 0.1394
					-9.229 0.1381
					-9.155 0.1326
					-9.082 0.1229
					-9.009 0.1091
					-8.935 0.0917
					-8.862 0.0712
					-8.789 0.0479
					-8.716 0.0226
					-8.643 -0.0040
					-8.569 -0.0312
					-8.496 -0.0584
					-8.423 -0.0847
					-8.350 -0.1093
					-8.276 -0.1317
					-8.203 -0.1513
					-8.130 -0.1673
					-8.057 -0.1795
					-7.983 -0.1874
					-7.910 -0.1908
					-7.837 -0.1894
					-7.764 -0.1834
					-7.690 -0.1728
					-7.617 -0.1578
					-7.544 -0.1386
					-7.471 -0.1158
					-7.397 -0.0897
					-7.324 -0.0609
					-7.251 -0.0301
					-7.178 0.0022
					-7.104 0.0351
					-7.031 0.0681
					-6.958 0.1004
					-6.885 0.1313
					-6.811 0.1603
					-6.738 0.1865
					-6.665 0.2097
					-6.592 0.2291
					-6.519 0.2445
					-6.445 0.2555
					-6.372 0.2618
					-6.299 0.2634
					-6.226 0.2601
					-6.152 0.2519
					-6.079 0.2391
					-6.006 0.2217
					-5.933 0.2000
					-5.859 0.1744
					-5.786 0.1453
					-5.713 0.1132
					-5.640 0.0785
					-5.566 0.0417
					-5.493 0.0035
					-5.420 -0.0357
					-5.347 -0.0751
					-5.273 -0.1143
					-5.200 -0.1527
					-5.127 -0.1898
					-5.054 -0.2250
					-4.981 -0.2579
					-4.907 -0.2880
					-4.834 -0.3149
					-4.761 -0.3384
					-4.688 -0.3581
					-4.614 -0.3738
					-4.541 -0.3853
					-4.468 -0.3925
					-4.394 -0.3954
					-4.321 -0.3939
					-4.248 -0.3881
					-4.175 -0.3780
					-4.102 -0.3639
					-4.028 -0.3458
					-3.955 -0.3240
					-3.882 -0.2988
					-3.809 -0.2703
					-3.735 -0.2390
					-3.662 -0.2051
					-3.589 -0.1690
					-3.516 -0.1309
					-3.442 -0.0913
					-3.369 -0.0505
					-3.296 -0.0088
					-3.223 0.0334
					-3.149 0.0758
					-3.076 0.1181
					-3.003 0.1600
					-2.930 0.2011
					-2.856 0.2414
					-2.783 0.2804
					-2.710 0.3180
					-2.637 0.3540
					-2.563 0.3883
					-2.490 0.4205
					-2.417 0.4508
					-2.344 0.4788
					-2.271 0.5047
					-2.197 0.5282
					-2.124 0.5493
					-2.051 0.5682
					-1.978 0.5846
					-1.904 0.5987
					-1.831 0.6105
					-1.758 0.6200
					-1.685 0.6274
					-1.611 0.6326
					-1.538 0.6357
					-1.465 0.6369
					-1.392 0.6363
					-1.318 0.6339
					-1.245 0.6298
					-1.172 0.6242
					-1.099 0.6172
					-1.025 0.6088
					-0.952 0.5993
					-0.879 0.5887
					-0.806 0.5771
					-0.732 0.5647
					-0.659 0.5515
					-0.586 0.5376
					-0.513 0.5232
					-0.440 0.5083
					-0.366 0.4931
					-0.293 0.4776
					-0.220 0.4619
					-0.146 0.4461
					-0.073 0.4302
					0.000 0.4143
					0.073 0.3985
					0.146 0.3828
					0.220 0.3673
					0.293 0.3520
					0.366 0.3370
					0.440 0.3222
					0.513 0.3078
					0.586 0.2937
					0.659 0.2800
					0.732 0.2666
					0.806 0.2537
					0.879 0.2412
					0.952 0.2290
					1.025 0.2174
					1.099 0.2061
					1.172 0.1953
					1.245 0.1849
					1.318 0.1749
					1.392 0.1654
					1.465 0.1563
					1.538 0.1476
					1.611 0.1393
					1.685 0.1313
					1.758 0.1238
					1.831 0.1167
					1.904 0.1099
					1.978 0.1034
					2.051 0.0973
					2.124 0.0915
					2.197 0.0860
					2.271 0.0808
					2.344 0.0760
					2.417 0.0713
					2.490 0.0670
					2.563 0.0629
					2.637 0.0590
					2.710 0.0554
					2.783 0.0519
					2.856 0.0487
					2.930 0.0457
					3.003 0.0428
					3.076 0.0402
					3.149 0.0377
					3.223 0.0353
					3.296 0.0331
					3.369 0.0311
					3.442 0.0291
					3.516 0.0273
					3.589 0.0256
					3.662 0.0240
					3.735 0.0226
					3.809 0.0212
					3.882 0.0199
					3.955 0.0186
					4.028 0.0175
					4.102 0.0164
					4.175 0.0155
					4.248 0.0145
					4.321 0.0137
					4.394 0.0128
					4.468 0.0121
					4.541 0.0114
					4.614 0.0107
					4.688 0.0101
					4.761 0.0095
					4.834 0.0089
					4.907 0.0084
					4.981 0.0080
				};
				\node[anchor=north west, font=\small, craneInk] at (axis cs:-24.5,1.0) {$t=1$};
				\nextgroupplot[xlabel={$x$}]
				\addplot[craneGoldDark, line width=1.0pt] table {
					-24.976 0.0374
					-24.902 0.0522
					-24.829 0.0657
					-24.756 0.0776
					-24.683 0.0875
					-24.609 0.0952
					-24.536 0.1005
					-24.463 0.1032
					-24.390 0.1033
					-24.316 0.1007
					-24.243 0.0954
					-24.170 0.0877
					-24.097 0.0776
					-24.023 0.0655
					-23.950 0.0516
					-23.877 0.0363
					-23.804 0.0200
					-23.730 0.0030
					-23.657 -0.0143
					-23.584 -0.0313
					-23.511 -0.0476
					-23.438 -0.0630
					-23.364 -0.0770
					-23.291 -0.0891
					-23.218 -0.0992
					-23.145 -0.1070
					-23.071 -0.1123
					-22.998 -0.1149
					-22.925 -0.1148
					-22.852 -0.1120
					-22.778 -0.1065
					-22.705 -0.0984
					-22.632 -0.0880
					-22.559 -0.0753
					-22.485 -0.0609
					-22.412 -0.0449
					-22.339 -0.0277
					-22.266 -0.0097
					-22.192 0.0086
					-22.119 0.0268
					-22.046 0.0447
					-21.973 0.0616
					-21.899 0.0773
					-21.826 0.0914
					-21.753 0.1036
					-21.680 0.1135
					-21.606 0.1211
					-21.533 0.1260
					-21.460 0.1282
					-21.387 0.1276
					-21.314 0.1243
					-21.240 0.1182
					-21.167 0.1095
					-21.094 0.0984
					-21.020 0.0851
					-20.947 0.0699
					-20.874 0.0531
					-20.801 0.0351
					-20.727 0.0161
					-20.654 -0.0033
					-20.581 -0.0228
					-20.508 -0.0420
					-20.435 -0.0605
					-20.361 -0.0779
					-20.288 -0.0939
					-20.215 -0.1080
					-20.142 -0.1202
					-20.068 -0.1300
					-19.995 -0.1373
					-19.922 -0.1419
					-19.849 -0.1438
					-19.775 -0.1428
					-19.702 -0.1391
					-19.629 -0.1326
					-19.556 -0.1235
					-19.482 -0.1119
					-19.409 -0.0982
					-19.336 -0.0824
					-19.263 -0.0649
					-19.189 -0.0461
					-19.116 -0.0263
					-19.043 -0.0059
					-18.970 0.0149
					-18.896 0.0355
					-18.823 0.0556
					-18.750 0.0748
					-18.677 0.0928
					-18.604 0.1093
					-18.530 0.1239
					-18.457 0.1364
					-18.384 0.1466
					-18.311 0.1542
					-18.237 0.1591
					-18.164 0.1613
					-18.091 0.1607
					-18.018 0.1573
					-17.944 0.1512
					-17.871 0.1424
					-17.798 0.1311
					-17.725 0.1174
					-17.651 0.1017
					-17.578 0.0842
					-17.505 0.0651
					-17.432 0.0448
					-17.358 0.0236
					-17.285 0.0019
					-17.212 -0.0200
					-17.139 -0.0418
					-17.065 -0.0630
					-16.992 -0.0834
					-16.919 -0.1026
					-16.846 -0.1203
					-16.773 -0.1363
					-16.699 -0.1502
					-16.626 -0.1619
					-16.553 -0.1711
					-16.480 -0.1778
					-16.406 -0.1818
					-16.333 -0.1830
					-16.260 -0.1815
					-16.186 -0.1772
					-16.113 -0.1703
					-16.040 -0.1608
					-15.967 -0.1488
					-15.894 -0.1346
					-15.820 -0.1183
					-15.747 -0.1002
					-15.674 -0.0805
					-15.601 -0.0595
					-15.527 -0.0376
					-15.454 -0.0150
					-15.381 0.0080
					-15.308 0.0311
					-15.234 0.0538
					-15.161 0.0760
					-15.088 0.0972
					-15.015 0.1173
					-14.941 0.1359
					-14.868 0.1528
					-14.795 0.1678
					-14.722 0.1806
					-14.648 0.1911
					-14.575 0.1992
					-14.502 0.2047
					-14.429 0.2076
					-14.355 0.2079
					-14.282 0.2055
					-14.209 0.2005
					-14.136 0.1929
					-14.062 0.1829
					-13.989 0.1705
					-13.916 0.1559
					-13.843 0.1393
					-13.770 0.1210
					-13.696 0.1010
					-13.623 0.0797
					-13.550 0.0573
					-13.477 0.0341
					-13.403 0.0104
					-13.330 -0.0136
					-13.257 -0.0376
					-13.184 -0.0613
					-13.110 -0.0845
					-13.037 -0.1069
					-12.964 -0.1282
					-12.891 -0.1482
					-12.817 -0.1667
					-12.744 -0.1835
					-12.671 -0.1983
					-12.598 -0.2112
					-12.524 -0.2218
					-12.451 -0.2302
					-12.378 -0.2362
					-12.305 -0.2397
					-12.231 -0.2408
					-12.158 -0.2394
					-12.085 -0.2356
					-12.012 -0.2294
					-11.938 -0.2208
					-11.865 -0.2100
					-11.792 -0.1971
					-11.719 -0.1821
					-11.646 -0.1653
					-11.572 -0.1468
					-11.499 -0.1269
					-11.426 -0.1056
					-11.352 -0.0832
					-11.279 -0.0600
					-11.206 -0.0360
					-11.133 -0.0116
					-11.060 0.0130
					-10.986 0.0376
					-10.913 0.0620
					-10.840 0.0860
					-10.767 0.1093
					-10.693 0.1318
					-10.620 0.1533
					-10.547 0.1735
					-10.474 0.1924
					-10.400 0.2098
					-10.327 0.2255
					-10.254 0.2394
					-10.181 0.2514
					-10.107 0.2615
					-10.034 0.2695
					-9.961 0.2753
					-9.888 0.2791
					-9.815 0.2807
					-9.741 0.2802
					-9.668 0.2775
					-9.595 0.2727
					-9.521 0.2658
					-9.448 0.2569
					-9.375 0.2461
					-9.302 0.2335
					-9.229 0.2192
					-9.155 0.2032
					-9.082 0.1857
					-9.009 0.1668
					-8.935 0.1468
					-8.862 0.1256
					-8.789 0.1035
					-8.716 0.0806
					-8.643 0.0570
					-8.569 0.0330
					-8.496 0.0087
					-8.423 -0.0158
					-8.350 -0.0403
					-8.276 -0.0646
					-8.203 -0.0887
					-8.130 -0.1123
					-8.057 -0.1354
					-7.983 -0.1577
					-7.910 -0.1792
					-7.837 -0.1998
					-7.764 -0.2194
					-7.690 -0.2378
					-7.617 -0.2549
					-7.544 -0.2708
					-7.471 -0.2852
					-7.397 -0.2982
					-7.324 -0.3097
					-7.251 -0.3196
					-7.178 -0.3279
					-7.104 -0.3346
					-7.031 -0.3397
					-6.958 -0.3432
					-6.885 -0.3450
					-6.811 -0.3452
					-6.738 -0.3438
					-6.665 -0.3408
					-6.592 -0.3362
					-6.519 -0.3302
					-6.445 -0.3227
					-6.372 -0.3138
					-6.299 -0.3035
					-6.226 -0.2920
					-6.152 -0.2792
					-6.079 -0.2653
					-6.006 -0.2503
					-5.933 -0.2343
					-5.859 -0.2174
					-5.786 -0.1996
					-5.713 -0.1811
					-5.640 -0.1619
					-5.566 -0.1421
					-5.493 -0.1218
					-5.420 -0.1011
					-5.347 -0.0800
					-5.273 -0.0586
					-5.200 -0.0370
					-5.127 -0.0154
					-5.054 0.0063
					-4.981 0.0280
					-4.907 0.0496
					-4.834 0.0710
					-4.761 0.0922
					-4.688 0.1131
					-4.614 0.1336
					-4.541 0.1537
					-4.468 0.1734
					-4.394 0.1926
					-4.321 0.2113
					-4.248 0.2293
					-4.175 0.2468
					-4.102 0.2636
					-4.028 0.2797
					-3.955 0.2952
					-3.882 0.3099
					-3.809 0.3239
					-3.735 0.3371
					-3.662 0.3495
					-3.589 0.3612
					-3.516 0.3721
					-3.442 0.3822
					-3.369 0.3916
					-3.296 0.4001
					-3.223 0.4079
					-3.149 0.4149
					-3.076 0.4212
					-3.003 0.4267
					-2.930 0.4314
					-2.856 0.4355
					-2.783 0.4388
					-2.710 0.4415
					-2.637 0.4435
					-2.563 0.4449
					-2.490 0.4456
					-2.417 0.4457
					-2.344 0.4453
					-2.271 0.4443
					-2.197 0.4427
					-2.124 0.4407
					-2.051 0.4382
					-1.978 0.4352
					-1.904 0.4318
					-1.831 0.4280
					-1.758 0.4238
					-1.685 0.4192
					-1.611 0.4143
					-1.538 0.4091
					-1.465 0.4036
					-1.392 0.3978
					-1.318 0.3918
					-1.245 0.3856
					-1.172 0.3792
					-1.099 0.3726
					-1.025 0.3659
					-0.952 0.3590
					-0.879 0.3520
					-0.806 0.3449
					-0.732 0.3377
					-0.659 0.3304
					-0.586 0.3231
					-0.513 0.3158
					-0.440 0.3084
					-0.366 0.3010
					-0.293 0.2937
					-0.220 0.2863
					-0.146 0.2790
					-0.073 0.2717
					0.000 0.2645
					0.073 0.2573
					0.146 0.2502
					0.220 0.2432
					0.293 0.2362
					0.366 0.2294
					0.440 0.2226
					0.513 0.2160
					0.586 0.2094
					0.659 0.2030
					0.732 0.1967
					0.806 0.1905
					0.879 0.1844
					0.952 0.1784
					1.025 0.1726
					1.099 0.1669
					1.172 0.1613
					1.245 0.1559
					1.318 0.1506
					1.392 0.1454
					1.465 0.1404
					1.538 0.1355
					1.611 0.1307
					1.685 0.1260
					1.758 0.1215
					1.831 0.1171
					1.904 0.1129
					1.978 0.1087
					2.051 0.1047
					2.124 0.1008
					2.197 0.0970
					2.271 0.0934
					2.344 0.0899
					2.417 0.0864
					2.490 0.0831
					2.563 0.0799
					2.637 0.0769
					2.710 0.0739
					2.783 0.0710
					2.856 0.0682
					2.930 0.0655
					3.003 0.0629
					3.076 0.0604
					3.149 0.0580
					3.223 0.0557
					3.296 0.0535
					3.369 0.0513
					3.442 0.0493
					3.516 0.0473
					3.589 0.0454
					3.662 0.0435
					3.735 0.0418
					3.809 0.0401
					3.882 0.0384
					3.955 0.0369
					4.028 0.0354
					4.102 0.0339
					4.175 0.0325
					4.248 0.0312
					4.321 0.0299
					4.394 0.0287
					4.468 0.0275
					4.541 0.0264
					4.614 0.0253
					4.688 0.0243
					4.761 0.0233
					4.834 0.0223
					4.907 0.0214
					4.981 0.0205
				};
				\node[anchor=north west, font=\small, craneGoldDark] at (axis cs:-24.5,1.0) {$t=3$};
				\draw[-{Stealth[length=2mm]}, craneGoldDark] (axis cs:-6,0.75) -- node[above, font=\scriptsize] {radiation} (axis cs:-14,0.75);
			\end{groupplot}
		\end{tikzpicture}
		\caption{Solution of the linear problem $\partial_{t}v-\partial_{x}|D_{x}|^{a}v=0$ with $a=3/2$ and $v(x,0)=e^{-x^{2}}$, computed with the explicit Fourier representation $\widehat{v}(\xi,t)=e^{it\xi|\xi|^{a}}\widehat{v}(\xi,0)$. The oscillatory tail moves to the left, and its wavelength decreases with $|x|$, since the waves located at $x\sim-(a+1)|\xi|^{a}t$ have frequency $|\xi|$.}
		\label{fig:linear}
	\end{figure}

	Before stating our results, let us describe the linear part of \eqref{dgbo}, which explains most of the phenomena considered here. A plane wave $e^{i(x\xi-\omega t)}$ solves $\partial_{t}v-\partial_{x}|D_{x}|^{a}v=0$ if and only if
	\begin{equation}\label{dispersionrelation}
		\omega(\xi)=-\xi|\xi|^{a}, \qquad\text{so that}\qquad \omega'(\xi)=-(a+1)|\xi|^{a}.
	\end{equation}
	Thus, the \emph{phase velocity} $\omega(\xi)/\xi=-|\xi|^{a}$ and the \emph{group velocity} $\omega'(\xi)=-(a+1)|\xi|^{a}$ are negative, and the frequencies $|\xi|\sim N$ travel to the left with speed $(a+1)N^{a}$, see Figure~\ref{fig:dispersion}. Since this speed increases with $|\xi|$, a \emph{localized perturbation spreads out immediately}: the high frequencies move towards $x=-\infty$ and leave behind an oscillatory tail, see Figure~\ref{fig:linear}. The exponent $a$ measures this effect. In the normalization of \eqref{dgbo}, the speed $(a+1)|\xi|^{a}$ equals $2|\xi|$ for $a=1$ and $3\xi^{2}$ for $a=2$. In this sense, the dispersion in \eqref{dgbo} is stronger than that of the BO equation and weaker than that of the KdV equation.
	
	The equation in \eqref{dgbo} also has \emph{travelling waves} $u(x,t)=\phi_{c}(x-ct)$, $c>0$, which move to the right, in the opposite direction of the linear radiation. The profile $\phi_{c}$ solves
	\begin{equation}\label{soliton}
		|D_{x}|^{a}\phi_{c}+c\,\phi_{c}=\frac{1}{2}\phi_{c}^{2}, \qquad \phi_{c}(x)=c\,\phi_{1}\big(c^{1/a}x\big).
	\end{equation}
	Frank and Lenzmann \cite{FL} proved that the \emph{ground states} of this equation are unique (up to translations) and that $|\phi_{1}(x)|\le c_{0}\,(1+|x|^{1+a})^{-1}$. Unlike the KdV case, where solitary waves decay exponentially, this decay is \emph{algebraic}. This comes from the lack of smoothness at the origin of the symbol $i\xi|\xi|^{a}$ of $\partial_{x}|D_{x}|^{a}$. In fact, \eqref{soliton} gives
	\begin{equation*}
		\widehat{\phi_{c}}(\xi)=\frac{\widehat{\phi_{c}^{2}}(\xi)}{2(c+|\xi|^{a})},
	\end{equation*}
	and the expansion of this expression near $\xi=0$ contains the term $-\frac{1}{2c^{2}}\|\phi_{c}\|_{2}^{2}|\xi|^{a}$. We will see below that \emph{this singularity determines the maximal decay} that a solution of \eqref{dgbo} can preserve in time.
	
	The first moment of the solution also has a simple dynamical meaning. Integrating by parts in \eqref{dgbo} and using that $\int_{\mathbb{R}}|D_{x}|^{a}u\,dx=0$, we get
	\begin{equation}\label{virial}
		\frac{d}{dt}\int_{\mathbb{R}}x\,u(x,t)\,dx=\frac{1}{2}\int_{\mathbb{R}}u^{2}(x,t)\,dx=\frac{1}{2}\|u_{0}\|_{2}^{2}.
	\end{equation}
	Hence, $\int_{\mathbb{R}} x\,u\,dx$ is an affine function of time with slope $\frac{1}{2}\|u_{0}\|_{2}^{2}$, which is positive unless $u\equiv 0$. If $\int_{\mathbb{R}} u_{0}\,dx\neq 0$, this means that the center of mass moves with constant velocity,
	\begin{equation*}
		\frac{d}{dt}\bigg(\frac{\int_{\mathbb{R}}x\,u(x,t)\,dx}{\int_{\mathbb{R}}u_{0}(x)\,dx}\bigg)=\frac{\|u_{0}\|_{2}^{2}}{2\int_{\mathbb{R}}u_{0}(x)\,dx}.
	\end{equation*} The identity \eqref{virial} plays a role in the unique continuation principles below.

	The IVP \eqref{dgbo} in the Sobolev spaces $H^{s}(\mathbb{R})$ has been extensively studied. Unlike the KdV equation, the dispersion in \eqref{dgbo} is not strong enough to treat the nonlinear term by a fixed point argument. Molinet, Saut and Tzvetkov \cite{MST} showed that for $1\le a<2$ the flow map data-solution of \eqref{dgbo} is not of class $C^{2}$ at the origin from $H^{s}(\mathbb{R})$ to $C([0,T]:H^{s}(\mathbb{R}))$, for any $s\in\mathbb{R}$ and $T>0$. For the BO equation, Koch and Tzvetkov \cite{Koch-Tzvetkov} showed that the flow map is not even uniformly continuous on bounded sets of $H^{s}(\mathbb{R})$, $s>0$. Thus, the local theory of \eqref{dgbo} requires more refined arguments. Herr \cite{Herr} proved well-posedness for data in $H^{s}(\mathbb{R})$ intersected with a weighted low-frequency space, using spaces of Bourgain type, and Guo \cite{ZG} adapted to \eqref{dgbo} the less perturbative approach of Ionescu, Kenig and Tataru. Herr, Ionescu, Kenig and Koch \cite{HIKK} established global well-posedness in the space of real-valued $L^{2}(\mathbb{R})$ functions by means of a para-differential re-normalization that controls the low-high frequency interactions. Molinet and Vento \cite{MolinetVento2015}, without a gauge transformation, obtained unconditional well-posedness in the energy space $H^{a/2}(\mathbb{R})$ by an improvement of the energy method for strongly non-resonant dispersive equations. Here we use the local theory in $H^{s}(\mathbb{R})$ of \cite{HIKK,MolinetVento2015}.

	Let us introduce the \emph{weighted Sobolev spaces}
	\begin{equation}\label{Zspaces}
		Z_{s,r}=H^{s}(\mathbb{R})\cap L^{2}(|x|^{2r}dx), \qquad
		\dot{Z}_{s,r}=\Big\{f\in Z_{s,r}\,:\,\widehat{f}(0)=\int_{\mathbb{R}}f(x)\,dx=0\Big\}, \qquad s,r\in\mathbb{R}.
	\end{equation}
	In general, for dispersive equations the decay of the data is preserved by the solution only if the data are regular enough. For \eqref{dgbo} this follows from the commutator identity
	\begin{equation}\label{commGamma}
		x\,e^{t\partial_{x}|D_{x}|^{a}}f=e^{t\partial_{x}|D_{x}|^{a}}(xf)-(1+a)\,t\,|D_{x}|^{a}e^{t\partial_{x}|D_{x}|^{a}}f,
	\end{equation}
	where $e^{t\partial_{x}|D_{x}|^{a}}$ denotes the solution operator of the linear problem $\partial_{t}v-\partial_{x}|D_{x}|^{a}v=0$. Equivalently, $\Gamma=x+(1+a)t|D_{x}|^{a}$ commutes with $\partial_{t}-\partial_{x}|D_{x}|^{a}$. Hence, along the linear flow and for $t\neq 0$, \emph{one unit of decay is equivalent to $a$ derivatives in $L^{2}(\mathbb{R})$}, and persistence of the solution of \eqref{dgbo} in $Z_{s,r}$ can only hold if $s\ge ar$. We emphasize that the factor $a$ is exactly the exponent in the group velocity $(a+1)|\xi|^{a}$.
	
	For the KdV equation, the persistence of polynomial decay goes back to Kato \cite{KATO1}. For the BO equation, Iorio \cite{IORIO} showed that polynomial decay is not preserved by the flow in general. The results in \cite{IORIO} were extended to fractional weights and sharpened as follows.
	
	\begin{thmx}[\cite{FP,FLP2}]\label{TheoremA}
		Let $u$ be the solution of the IVP associated with the BO equation \eqref{BO}.
		\begin{itemize}
			\item[(i)] If $s\ge 1$, $r\in[0,s]$, $r<5/2$ and $u_{0}\in Z_{s,r}$, then $u\in C([0,\infty):Z_{s,r})$.
			\item[(ii)] If $r\in[5/2,7/2)$, $r\le s$ and $u_{0}\in\dot{Z}_{s,r}$, then $u\in C([0,\infty):\dot{Z}_{s,r})$.
			\item[(iii)] Let $u\in C([0,T]:Z_{2,2})$. If there exist two different times $t_{1},t_{2}\in[0,T]$ such that $u(\cdot,t_{j})\in Z_{5/2,5/2}$, $j=1,2$, then $\widehat{u}_{0}(0)=0$.
			\item[(iv)] Let $u\in C([0,T]:\dot{Z}_{3,3})$. If there exist three different times $t_{1},t_{2},t_{3}\in[0,T]$ such that $u(\cdot,t_{j})\in\dot{Z}_{7/2,7/2}$, $j=1,2,3$, then $u\equiv 0$.
			\item[(v)] If $u_{0}\in\dot{Z}_{5,4}$ and $\int_{\mathbb{R}} x\,u_{0}(x)\,dx\neq 0$, then $u(\cdot,t^{*})\in\dot{Z}_{5,4}$, where
			\begin{equation*}
				t^{*}=-\frac{4}{\|u_{0}\|_{2}^{2}}\int_{\mathbb{R}}x\,u_{0}(x)\,dx.
			\end{equation*}
		\end{itemize}
	\end{thmx}
	Part (v) of Theorem~\ref{TheoremA} shows that part (iv) fails if one assumes the condition at only two times. This is in contrast with the generalized KdV and the semi-linear Schr\"odinger equations, where \emph{unique continuation principles} of this kind hold with two times (see the references in \cite{FONLP}). For the IVP \eqref{dgbo}, the analogous results were obtained by Fonseca, Linares and Ponce \cite{FONLP}. In \cite{FONLP} the equation is written as $$\partial_{t}w+|D_{x}|^{1+b}\partial_{x}w+w\,\partial_{x}w=0,\quad 0<b<1.$$
	
	After the change of variables $w(x,t)=-u(x,-t)$ and $a=1+b$, their results read as follows.
	\begin{thmx}[\cite{FONLP}]\label{TheoremB}
		Let $a\in(1,2)$ and let $u$ be the solution of the IVP \eqref{dgbo}.
		\begin{itemize}
			\item[(i)] If $u_{0}\in Z_{s,r}$, then $u\in C([-T,T]:Z_{s,r})$ if either $s\ge a$ and $r\in(0,1]$, or $s\ge 2a$ and $r\in(1,2]$, or $s\ge[(r+1)^{-}]\,a$ and $2<r<3/2+a$, where $[\cdot]$ denotes the integer part function.
			\item[(ii)] If $u_{0}\in\dot{Z}_{s,r}$, $s\ge[(r+1)^{-}]\,a$ and $3/2+a\le r<5/2+a$, then $u\in C([-T,T]:\dot{Z}_{s,r})$.
			\item[(iii)] Let $u\in C([-T,T]:Z_{s,(3/2+a)^{-}})$ with $s\ge a(3/2+a)+(2-a)/2$. If there exist two different times $t_{1},t_{2}\in[-T,T]$ such that $u(\cdot,t_{j})\in Z_{s,3/2+a}$, $j=1,2$, then $$\int_{\mathbb{R}} u_{0}(x)\,dx=0.$$
			\item[(iv)] Let $u\in C([-T,T]:Z_{s,(5/2+a)^{-}})$ with $s\ge a(5/2+a)+(2-a)/2$. If there exist three different times $t_{1},t_{2},t_{3}\in[-T,T]$ such that $u(\cdot,t_{j})\in\dot{Z}_{s,5/2+a}$, $j=1,2,3$, then $u\equiv 0$. The same conclusion holds with two different times $t_{1},t_{2}$, provided that 
			\begin{equation*}
				\int_{\mathbb{R}} x\,u(x,t_{1})\,dx=0\quad \mbox{or}\quad  \int_{\mathbb{R}} x\,u(x,t_{2})\,dx=0.
			\end{equation*}
			
		\end{itemize}
	\end{thmx}
	In \cite{FONLP} it was also shown that two times are not enough in part (iv) of Theorem~\ref{TheoremB}. More precisely, setting $\tilde{a}=[2a-1]/2$, if $u_{0}\in\dot{Z}_{s,7/2+\tilde{a}}$ with $s\ge a(7/2+\tilde{a})+(2-a)/2$ and ${\displaystyle \int_{\mathbb{R}} x\,u_{0}(x)\,dx\neq 0}$, then $u(\cdot,t^{*})\in\dot{Z}_{s,7/2+\tilde{a}}$ with $t^{*}$ as in Theorem~\ref{TheoremA} (v). For the semi-linear Schr\"odinger equation and the fractional KdV equation see \cite{NP} and \cite{Riano2021}, respectively.
	
	Notice that the regularity in Theorem~\ref{TheoremB} is in general larger than the one given by \eqref{commGamma}. For instance, for $r\in(2,3]$, Theorem~\ref{TheoremB} (i) requires $s\ge 3a$, while \eqref{commGamma} suggests $s\ge ar$. The \emph{converse question}, that is, \emph{whether the decay of the solution implies some regularity}, was recently studied by Linares and Ponce \cite{LPBO} for the BO equation.
	\begin{thmx}[\cite{LPBO}]\label{TheoremC}
		Let $u$ be a solution of the IVP associated with the BO equation \eqref{BO}.
		\begin{itemize}
			\item[(i)] Let $u\in C([0,T]:L^{2}(\mathbb{R}))$ be the solution obtained in \cite{MOLINET3}. If there exist $t_{1},t_{2}\in[0,T]$, $t_{1}<t_{2}$, and $\alpha\in(0,1/2]$ such that $|x|^{\alpha}u(x,t_{j})\in L^{2}(\mathbb{R})$, $j=1,2$, then $u\in C([0,T]:H^{\alpha}(\mathbb{R}))$.
			\item[(ii)] Let $u\in C([0,T]:H^{s}(\mathbb{R}))$, $s>3/2$. If there exist $t_{1},t_{2}\in[0,T]$, $t_{1}<t_{2}$, and $\alpha\in(s,7/2)$ such that $|x|^{\alpha}u(x,t_{j})\in L^{2}(\mathbb{R})$, $j=1,2$, then $u\in C([0,T]:H^{\alpha}(\mathbb{R}))$.
		\end{itemize}
	\end{thmx}
	Roughly, Theorem~\ref{TheoremC} says that, for $\beta<7/2$, the $L^{2}(|x|^{2\beta}dx)$-norm of the solution of the BO equation remains bounded in the existence interval if and only if the data belong to this space and $u_{0}\in H^{\beta}(\mathbb{R})$. Results of this type were obtained for the generalized KdV equation by Isaza, Linares and Ponce \cite{ILP3}, and for the intermediate long wave equation in \cite{LPILW}. As pointed out in \cite{LPBO}, for the BO equation the gain of regularity comes from an operator with non-smooth symbol and the range of decay is a priori limited, which makes the argument more involved. It is natural to ask whether this characterization holds for a model with dispersion between those of the BO and the KdV equations, that is, for the IVP \eqref{dgbo}, and whether the regularity in Theorem~\ref{TheoremB} can be lowered to the sharp level $s\ge ar$. Our main results give a positive answer to both questions, except for the cases described in Remark~\ref{remarkAA2} and at the end of this section, see also Figure~\ref{fig:diagram}.

	We present our first main result, which shows that \emph{decay at two different times implies a gain of regularity} of exactly the order predicted by \eqref{commGamma}.
	\begin{thm}\label{main}
		We consider $a\in (1,2)$.
		\begin{itemize}
			\item[(i)] Let $u\in C([-T,T]:H^s(\mathbb{R}))$, with $s\geq 1/6$, be the solution of \eqref{dgbo} provided by \cite{HIKK,MolinetVento2015}. If there exist $t_1,t_2\in [-T,T]$ with $t_1<t_2$, and $\alpha\in (s/a,1/2]$ such that
			\begin{equation*}
				|x|^{\alpha}u(x,t_j)\in L^2(\mathbb{R}), \quad j=1,2,  
			\end{equation*}
			then
			\begin{equation*}
				u\in C([-T,T]:H^{a \alpha }(\mathbb{R})).  
			\end{equation*}
			\item[(ii)] Let $u\in C([-T,T]:H^s(\mathbb{R}))$, with $s>3/2$, be the solution of \eqref{dgbo}. If there exist $t_1,t_2\in [-T,T]$ with $t_1<t_2$, and $\alpha\in (s/a,5/2+a)$ such that
			\begin{equation*}
				|x|^{\alpha}u(x,t_j)\in L^2(\mathbb{R}), \quad j=1,2,  
			\end{equation*}
			then
			\begin{equation*}
				u\in C([-T,T]:H^{a \alpha}(\mathbb{R})).  
			\end{equation*}
		\end{itemize}
	\end{thm}
	Our second main result concerns persistence in the spaces \eqref{Zspaces} with the sharp regularity, and the corresponding unique continuation principles.
	\begin{thm}\label{AA2}
		Let $a\in(1,2)$.
		\begin{enumerate}
			\item[(i)] Let $\alpha\in(0,3/2+a)$. Let $u\in C([-T,T] : H^s(\mathbb{R}))$, with $s\geq \max\{a/2,\, a \alpha\}$, be a solution of the IVP \eqref{dgbo}. If $|x|^{\alpha}u_0\in L^2(\mathbb{R})$, then  $|x|^{\alpha}u\in  C([-T,T] : L^2(\mathbb{R})).$
			\vskip.1in
			\item[(ii)]
			Let  $u\in C([-T,T] : H^{a(3/2+a)^{-}}(\mathbb{R}))$ be a solution of the IVP \eqref{dgbo}. If there exist two different times $t_1, t_2\in [-T,T]$ such that $|x|^{3/2+a}u(x,t_j)\in L^2(\mathbb{R}),\,j=1,2$, then $\int_{\mathbb{R}} u_0(x)dx=0$.
			\vskip.1in
			\item[(iii)]
			Let $\alpha\in[3/2+a,5/2+a)$. Let $u\in C([-T,T] : H^{s}(\mathbb{R}))$, with $s\geq a\alpha $, be a solution of the IVP \eqref{dgbo} with $\int_{\mathbb{R}} u_0(x)dx=0$.  If $|x|^{\alpha}u_0\in L^2(\mathbb{R})$, then  $|x|^{\alpha}u\in  C([-T,T] : L^2(\mathbb{R})).$
			
			\vskip.1in
			\item[(iv)]
			Let  $u\in C([-T,T] : H^{a(5/2+a)^{-}}(\mathbb{R}))$ be a solution of the IVP \eqref{dgbo}. If there exist three different times $t_1, t_2,t_3\in [-T,T]$ such that $|x|^{5/2+a}u(x,t_j)\in L^2(\mathbb{R}),\,j=1,2,3$, then $u\equiv 0$.
		\end{enumerate}
	\end{thm}
	\begin{rem}\label{remarkmain}
		\begin{itemize}
			\item[(a)] Formally, taking $a=1$ in Theorem~\ref{main} (ii) and using the change of variables relating \eqref{dgbo} and \eqref{BO}, we recover the range $\alpha\in(s,7/2)$ and the regularity $H^{\alpha}(\mathbb{R})$ in Theorem~\ref{TheoremC} (ii), since in this case $5/2+a=7/2$ and $a\alpha=\alpha$.
			\item[(b)] In part (i), the nonlinear term in the weighted energy estimate is controlled by $\|u\|_{3}^{3}$. This explains the condition $s\ge 1/6$, by the embedding $H^{1/6}(\mathbb{R})\hookrightarrow L^{3}(\mathbb{R})$.
		\end{itemize}
	\end{rem}
	\begin{rem}\label{remarkAA2}
		\begin{itemize}
			\item[(a)] In view of the local theory in $H^{s}(\mathbb{R})$ (see \cite{HIKK,MolinetVento2015}), when $\frac{1}{2}\le\alpha<\frac{5}{2}+a$, parts (i) and (iii) of Theorem~\ref{AA2} improve the persistence results in Theorem~\ref{TheoremB} to the sharp regularity $s\ge a\alpha$. \emph{To our knowledge, it is an open problem whether $s\ge a\alpha$ suffices for $\alpha\in(0,\frac{1}{2})$.}
			\item[(b)] The unique continuation principles in parts (ii) and (iv) of Theorem~\ref{AA2} are optimal with respect to the regularity, since we only assume $s\ge a(3/2+a)^{-}$ and $s\ge a(5/2+a)^{-}$, respectively. Theorem~\ref{TheoremB} requires $(2-a)/2$ additional derivatives.
			\item[(c)] By part (ii) of Theorem~\ref{AA2}, persistence in $Z_{s,\alpha}$ with $\alpha=(3/2+a)^{-}$ is the best possible for data with $\int_{\mathbb{R}} u_{0}\,dx\neq 0$. By parts (iii) and (iv), $\alpha=(5/2+a)^{-}$ is the largest possible decay for non-trivial solutions in $C([-T,T]:H^{a(5/2+a)^{-}}(\mathbb{R}))$.
		\end{itemize}
	\end{rem}
	Combining Theorems~\ref{main} and \ref{AA2}, we obtain the following \emph{characterization}. Let $u\in C([-T,T]:H^{s}(\mathbb{R}))$, $s>3/2$, be a solution of \eqref{dgbo} and let $\alpha\in(s/a,3/2+a)$. Then $|x|^{\alpha}u\in L^{\infty}([-T,T]:L^{2}(\mathbb{R}))$ if and only if $|x|^{\alpha}u_{0}\in L^{2}(\mathbb{R})$ and $u_{0}\in H^{a\alpha}(\mathbb{R})$. For $\alpha\in[3/2+a,5/2+a)$, $|x|^{\alpha}u\in L^{\infty}([-T,T]:L^{2}(\mathbb{R}))$ if and only if $|x|^{\alpha}u_{0}\in L^{2}(\mathbb{R})$, $u_{0}\in H^{a\alpha}(\mathbb{R})$ and $\int_{\mathbb{R}} u_{0}\,dx=0$, where the zero mean condition is necessary by Theorem~\ref{AA2} (ii). This extends Theorem~\ref{TheoremC} to the DGBO equation. Figure~\ref{fig:diagram} shows the regions of the $(\alpha,s)$ plane covered by our results. Moving to the right means more decay, and moving upwards means more regularity. Theorem~\ref{AA2} describes the points preserved by the flow, that is, those on or above the line $s=a\alpha$ (regions A and B). Theorem~\ref{main} describes a solution that starts below this line and is localized at two different times: its regularity is raised up to the line $s=a\alpha$ (regions C and D). The vertical lines $\alpha=\frac{3}{2}+a$ and $\alpha=\frac{5}{2}+a$ do not depend on the regularity. They come from the low frequencies, namely, the zero mean condition and the maximal decay.
	\begin{figure}[ht]
		\centering
		\begin{tikzpicture}[
			lab/.style={circle, draw=#1, fill=white, text=#1, line width=0.7pt, inner sep=1.2pt, font=\footnotesize\bfseries},
			lab/.default=craneNavy]
			\begin{axis}[
				width=0.92\textwidth, height=9cm,
				xmin=0, xmax=4.6, ymin=0, ymax=6.9,
				axis lines=left, axis line style={-{Stealth[length=2.2mm]}, craneNavy},
				xlabel={decay exponent $\alpha$ (weight $|x|^{\alpha}$)}, ylabel={regularity $s$ (space $H^{s}(\mathbb{R})$)},
				xtick={0.5,1,3,4}, xticklabels={$\frac{1}{2}$,$\frac{3}{2a}$,$\frac{3}{2}+a$,$\frac{5}{2}+a$},
				ytick={0.1667,0.75,1.5}, yticklabels={$\frac{1}{6}$,$\frac{a}{2}$,$\frac{3}{2}$},
				tick style={craneNavy}, clip=false
				]
				\fill[craneInk!18] (axis cs:0,0.75) -- (axis cs:0.5,0.75) -- (axis cs:3,4.5) -- (axis cs:3,6.75) -- (axis cs:0,6.75) -- cycle;
				\fill[craneGoldLight] (axis cs:3,4.5) -- (axis cs:4,6) -- (axis cs:4,6.75) -- (axis cs:3,6.75) -- cycle;
				\fill[craneGold!35] (axis cs:1,1.5) -- (axis cs:4,1.5) -- (axis cs:4,6) -- cycle;
				\fill[craneGold!35] (axis cs:0.1111,0.1667) -- (axis cs:0.5,0.1667) -- (axis cs:0.5,0.75) -- cycle;
				\fill[pattern=north east lines, pattern color=craneGray] (axis cs:0,0) -- (axis cs:0.5,0.75) -- (axis cs:0,0.75) -- cycle;
				\fill[pattern=dots, pattern color=craneGray] (axis cs:0,0) -- (axis cs:4,0) -- (axis cs:4,1.5) -- (axis cs:1,1.5) -- (axis cs:0.5,0.75) -- (axis cs:0.5,0.1667) -- (axis cs:0.1111,0.1667) -- cycle;
				\fill[pattern=crosshatch, pattern color=craneNavy!55] (axis cs:4,0) rectangle (axis cs:4.45,6.75);
				\addplot[craneNavy, line width=1.4pt, domain=0:4.45, samples=2] {1.5*x};
				\node[craneNavy, rotate=34, font=\small, anchor=south, fill=white, fill opacity=0.8, text opacity=1, inner sep=1pt] at (axis cs:1.45,2.18) {sharp line $s=a\alpha$};
				\draw[craneNavy, densely dashed, line width=0.9pt] (axis cs:3,0) -- (axis cs:3,6.75);
				\draw[craneNavy, line width=1.1pt] (axis cs:4,0) -- (axis cs:4,6.75);
				\draw[craneGray!80!black, densely dotted] (axis cs:0,1.5) -- (axis cs:4,1.5);
				\draw[craneGray!80!black, densely dotted] (axis cs:0,0.75) -- (axis cs:0.5,0.75) -- (axis cs:0.5,0);
				\draw[craneGray!80!black, densely dotted] (axis cs:1,0) -- (axis cs:1,1.5);
				\fill[craneGoldDark] (axis cs:2.4,2.0) circle (1.8pt);
				\fill[craneGoldDark] (axis cs:2.4,3.6) circle (1.8pt);
				\draw[-{Stealth[length=2.6mm]}, line width=1.1pt, craneGoldDark] (axis cs:2.4,2.06) -- (axis cs:2.4,3.52);
				\node[craneGoldDark, font=\scriptsize, anchor=north] at (axis cs:2.4,1.97) {$(\alpha,s)$};
				\node[craneGoldDark, font=\scriptsize, anchor=east] at (axis cs:2.36,3.68) {$(\alpha,a\alpha)$};
				\node[craneGoldDark, font=\scriptsize, align=left, anchor=west, fill=white, fill opacity=0.85, text opacity=1, inner sep=1.5pt] at (axis cs:2.48,2.75) {$|x|^{\alpha}u(t_{1}),|x|^{\alpha}u(t_{2})\in L^{2}$\\ $\Longrightarrow\ u\in C([-T,T]:H^{a\alpha})$};
				\node[craneNavy, font=\scriptsize, align=center, anchor=south, fill=white, inner sep=1pt] at (axis cs:3,6.78) {Thm.~\ref{AA2} (ii):\\ $\int_{\mathbb{R}} u_{0}\,dx=0$};
				\node[craneNavy, font=\scriptsize, align=center, anchor=south, fill=white, inner sep=1pt] at (axis cs:4.2,6.78) {Thm.~\ref{AA2} (iv):\\ $u\equiv 0$};
				\node[lab=craneInk] at (axis cs:1.0,4.6) {A};
				\node[lab=craneGoldDark] at (axis cs:3.55,6.25) {B};
				\node[lab=craneGoldDark] at (axis cs:3.55,2.1) {C};
				\node[lab=craneGray!70!black] at (axis cs:2.6,0.8) {F};
				\node[lab=craneNavy] at (axis cs:4.225,3.3) {G};
				\draw[craneGoldDark, -{Stealth[length=1.6mm]}] (axis cs:0.85,0.35) node[lab=craneGoldDark, anchor=west] {D} -- (axis cs:0.42,0.38);
				\draw[craneGray!70!black, -{Stealth[length=1.6mm]}] (axis cs:-0.45,1.2) node[lab=craneGray!70!black, anchor=east] {E} -- (axis cs:0.12,0.6);
			\end{axis}
		\end{tikzpicture}
		
		\vspace{2mm}
		{\footnotesize
			\begin{tabular}{c l l}
				\hline
				Region & Parameters & Result \\
				\hline
				A & $0<\alpha<\frac{3}{2}+a$, $s\geq\max\{\frac{a}{2},a\alpha\}$ & persistence: $|x|^{\alpha}u_{0}\in L^{2}\Rightarrow|x|^{\alpha}u\in C([-T,T]:L^{2})$, Thm.~\ref{AA2} (i)\\
				B & $\frac{3}{2}+a\leq\alpha<\frac{5}{2}+a$, $s\geq a\alpha$ & persistence for data with $\int_{\mathbb{R}} u_{0}\,dx=0$, Thm.~\ref{AA2} (iii)\\
				C & $s>\frac{3}{2}$, $\frac{s}{a}<\alpha<\frac{5}{2}+a$ & decay at two times $\Rightarrow u\in C([-T,T]:H^{a\alpha})$, Thm.~\ref{main} (ii)\\
				D & $s\geq\frac{1}{6}$, $\frac{s}{a}<\alpha\leq\frac{1}{2}$ & decay at two times $\Rightarrow u\in C([-T,T]:H^{a\alpha})$, Thm.~\ref{main} (i)\\
				E & $0<\alpha<\frac{1}{2}$, $a\alpha\leq s<\frac{a}{2}$ & persistence with $s\geq a\alpha$ open\\
				F & $\alpha>\frac{s}{a}$, $s\leq\frac{3}{2}$, outside D & not covered by Thm.~\ref{main}\\
				G & $\alpha\geq\frac{5}{2}+a$ & no non-trivial persistence, Thm.~\ref{AA2} (iv)\\
				\hline
		\end{tabular}}
		\caption{Summary of the results in the plane $(\alpha,s)$, where $\alpha$ is the decay exponent of the weight $|x|^{\alpha}$ and $s$ is the Sobolev regularity of the solution (drawn for $a=3/2$, so that $3/(2a)=1$). The line $s=a\alpha$ is the sharp relation between decay and regularity dictated by \eqref{commGamma}. Above this line (regions A and B) the decay of the data is preserved by the flow, and in region B the zero mean condition is required. Below this line (regions C and D) the regularity is too low for persistence, and Theorem~\ref{main} shows that decay at two different times moves the solution vertically up to the line $s=a\alpha$, as indicated by the arrow. The dashed line $\alpha=\frac{3}{2}+a$ is the threshold of Theorem~\ref{AA2} (ii): decay of order $\frac{3}{2}+a$ at two times forces $\int_{\mathbb{R}} u_{0}\,dx=0$. The solid line $\alpha=\frac{5}{2}+a$ is the maximal decay rate: by Theorem~\ref{AA2} (iv), a solution with decay of order $\frac{5}{2}+a$ at three times vanishes identically (region G). Regions E and F correspond to the cases left open.}
		\label{fig:diagram}
	\end{figure}
	From the physical point of view, our results can be explained by the dispersive mechanism in \eqref{dispersionrelation}. Theorem~\ref{main} says that a wave which is localized, in the sense that $|x|^{\alpha}u(\cdot,t)\in L^{2}(\mathbb{R})$, at two different times cannot carry much energy at high frequencies. In fact, the frequencies $|\xi|\sim N$ travel to the left with speed $(a+1)N^{a}$, so in the interval $[t_{1},t_{2}]$ they move a distance of order $N^{a}(t_{2}-t_{1})$, and they would destroy the localization at time $t_{2}$ unless their amplitude is small. In particular, \emph{rough data cannot remain localized}. Under the assumptions of Theorem~\ref{main}, if $|x|^{\alpha}u_{0}\in L^{2}(\mathbb{R})$ but $u_{0}\notin H^{a\alpha}(\mathbb{R})$, then $|x|^{\alpha}u(\cdot,t)\notin L^{2}(\mathbb{R})$ for every $t\neq 0$, that is, the high-frequency part of the solution is immediately radiated to the left. On the other hand, by Theorem~\ref{AA2} (i) and (iii), \emph{regular data remain localized}. In this sense, the relation $s=a\alpha$ is the balance between the regularity of the wave and the size of the region it occupies.
	
	The thresholds $3/2+a$ and $5/2+a$ have a different origin. They come from the \emph{long waves}, that is, from the frequencies near $\xi=0$, where the symbol $i\xi|\xi|^{a}$ is not smooth. Formally, a solution with $\int_{\mathbb{R}}u_{0}\,dx\neq 0$ develops immediately a tail of the form $c\,t\,\big(\int_{\mathbb{R}}u_{0}\,dx\big)|x|^{-(a+2)}$, since the inverse Fourier transform of $\xi|\xi|^{a}$ near the origin behaves like $|x|^{-(a+2)}$ (compare with \eqref{ee2}). This tail belongs to $L^{2}(|x|^{2\alpha}dx)$ if and only if $\alpha<a+3/2$, which explains Theorem~\ref{AA2} (ii). When the mean is zero, the next tail decays like $|x|^{-(a+3)}$ with coefficient $\int_{0}^{t}\int_{\mathbb{R}}x\,u(x,\tau)\,dx\,d\tau$, which by \eqref{virial} is determined by the motion of the first moment. Since the nonlinear term makes the first moment grow at the rate $\frac{1}{2}\|u_{0}\|_{2}^{2}$, this coefficient is a quadratic polynomial in time, and \emph{for a non-trivial solution it vanishes at most at two times}. This is the content of Theorem~\ref{AA2} (iv). For the KdV equation, $a=2$, the symbol $i\xi^{3}$ is a polynomial, no tails appear, and fast decay is preserved by the flow. For $1\leq a<2$ the decay of solutions is limited, as it is for the solitary waves $\phi_{c}$ in \eqref{soliton}, which decay like $|x|^{-(1+a)}$.
	
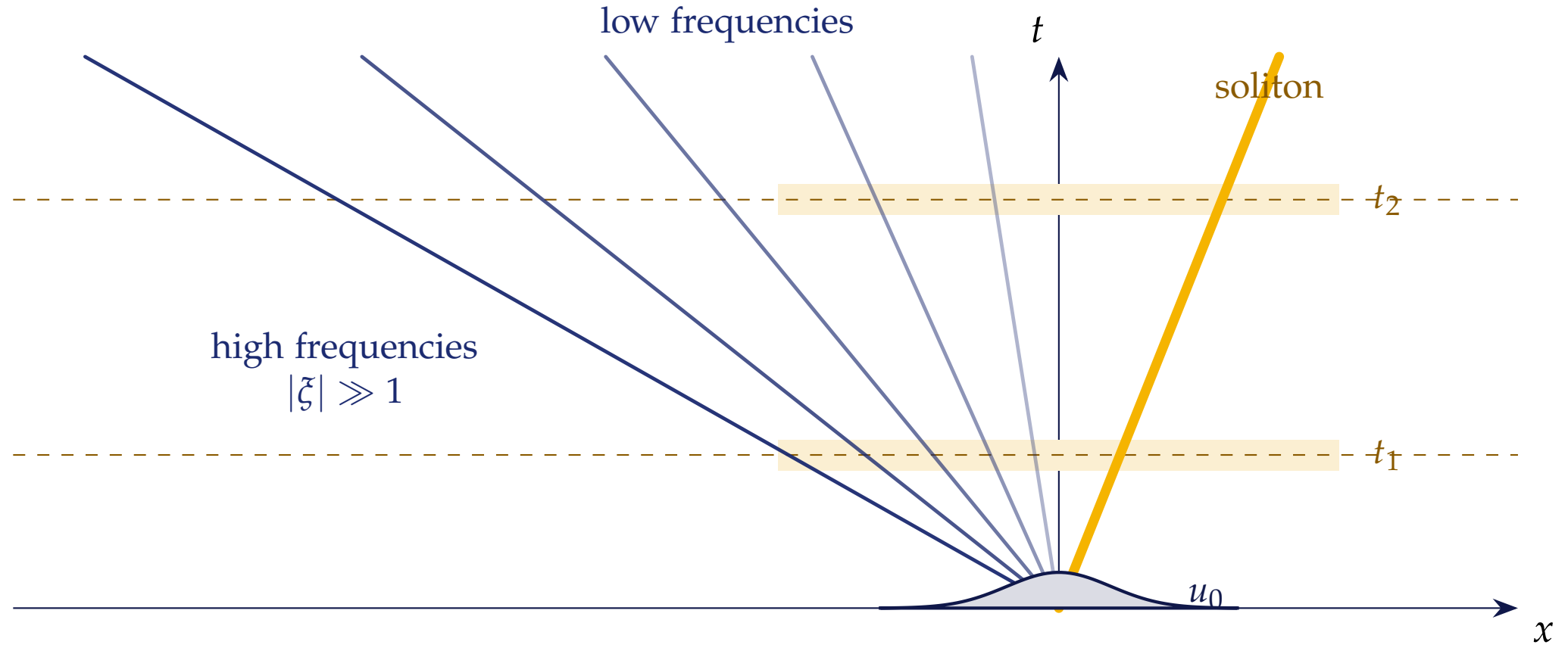
\begin{figure}[ht]
		\centering
		\begin{tikzpicture}
			\begin{axis}[
				width=0.92\textwidth, height=6.4cm,
				xmin=-8.2, xmax=3.6, ymin=0, ymax=1.08,
				axis lines=middle, axis line style={-{Stealth[length=2.2mm]}, craneNavy},
				xlabel={$x$}, ylabel={$t$},
				xlabel style={at={(ticklabel* cs:1.0)}, anchor=north west},
				ylabel style={at={(ticklabel* cs:1.0)}, anchor=south east},
				xtick=\empty, ytick=\empty,
				clip=false, every axis plot/.append style={line cap=round}
				]
				\fill[craneGoldLight] (axis cs:-2.2,0.27) rectangle (axis cs:2.2,0.33);
				\fill[craneGoldLight] (axis cs:-2.2,0.77) rectangle (axis cs:2.2,0.83);
				\draw[craneGoldDark, dashed, thin] (axis cs:-8.2,0.30) -- (axis cs:3.6,0.30);
				\draw[craneGoldDark, dashed, thin] (axis cs:-8.2,0.80) -- (axis cs:3.6,0.80);
				\foreach \v/\op in {0.63/0.35, 1.79/0.5, 3.29/0.65, 5.06/0.8, 7.07/0.95}{
					\edef\temp{\noexpand\addplot[craneInk, line width=0.9pt, opacity=\op, domain=0:1.08, samples=2] ({-\v*x},{x});}
					\temp
				}
				\addplot[craneGold, line width=2.4pt, domain=0:1.08, samples=2] ({1.6*x},{x});
				\addplot[craneNavy, fill=craneNavy!15, line width=0.8pt, domain=-1.4:1.4, samples=80]
				({x},{0.07*exp(-3*x*x)}) \closedcycle;
				\node[craneGoldDark, anchor=west, font=\small] at (axis cs:2.35,0.30) {$t_{1}$};
				\node[craneGoldDark, anchor=west, font=\small] at (axis cs:2.35,0.80) {$t_{2}$};
				\node[craneGoldDark, font=\small, anchor=south] at (axis cs:1.65,0.97) {soliton};
				\node[craneInk, font=\small, anchor=north, align=center, fill=white, inner sep=1pt] at (axis cs:-5.6,0.55)
				{high frequencies\\ $|\xi|\gg 1$};
				\node[craneInk, font=\small, anchor=south, align=center] at (axis cs:-2.6,1.08)
				{low frequencies};
				\node[craneNavy, font=\small, anchor=north west] at (axis cs:0.9,0.075) {$u_{0}$};
			\end{axis}
		\end{tikzpicture}
		\caption{Space-time picture of the linear dynamics of \eqref{dgbo} (drawn for $a=3/2$). A wave packet concentrated at frequency $\xi$ moves along the ray $x=-(a+1)|\xi|^{a}t$, while solitary waves travel to the right. If $|x|^{\alpha}u(\cdot,t_{j})\in L^{2}(\mathbb{R})$ at two different times (shaded windows), the high-frequency part of the solution, which crosses arbitrarily long distances in the time interval $[t_{1},t_{2}]$, must be small. This is the heuristic behind Theorem~\ref{main}.}
		\label{fig:rays}
	\end{figure}

	Next, we describe the main ideas in the proof of Theorem~\ref{main}. At the linear level, the identity \eqref{commGamma} already gives the result for $\alpha=1$. If $v(t)=e^{t\partial_{x}|D_{x}|^{a}}v_{0}$ and $t\neq 0$, then
	\begin{equation*}
		|D_{x}|^{a}v(t)=\frac{1}{(1+a)t}\Big(e^{t\partial_{x}|D_{x}|^{a}}(xv_{0})-x\,v(t)\Big),
	\end{equation*}
	and therefore $xv_{0},\,xv(t)\in L^{2}(\mathbb{R})$ imply $|D_{x}|^{a}v(t)\in L^{2}(\mathbb{R})$. In view of \eqref{dispersionrelation}, this means that the frequencies of size $N$ move a distance $(a+1)N^{a}(t_{2}-t_{1})$ in $[t_{1},t_{2}]$, so a solution localized at the times $t_{1}$ and $t_{2}$ cannot carry much energy at high frequencies, see Figure~\ref{fig:rays}.
	For the nonlinear problem and fractional weights, we follow the ideas in \cite{LPBO}, that is, an iterative argument combined with weighted energy estimates in the spirit of the Kato smoothing effect \cite{KATO1}. More precisely, we multiply the equation in \eqref{dgbo} by $x\langle x\rangle_{N}^{2\alpha-1}u$, where $\langle x\rangle_{N}$ is a bounded approximation of $\langle x\rangle=(1+x^{2})^{1/2}$, and integrate in space. The dispersive term gives
	\begin{equation*}
		\begin{split}
			\int_{\mathbb{R}}|D_{x}|^{a}\partial_{x}u\;u\,x\langle x\rangle_{N}^{2\alpha-1}\,dx
			&=-\Big(\frac{a+1}{2}\Big)\int_{\mathbb{R}}\big||D_{x}|^{a/2}u\big|^{2}\big(x\langle x\rangle_{N}^{2\alpha-1}\big)'\,dx\\
			&\quad+\big(\text{lower order terms}\big).
		\end{split}
	\end{equation*}
	Since the weight $x\langle x\rangle_{N}^{2\alpha-1}$ is increasing, its derivative measures the flux of energy carried by the left-going waves, and integrating between $t_{1}$ and $t_{2}$ we gain $a/2$ derivatives at the price of $1/2$ of the weight, that is,
	\begin{equation*}
		\langle x\rangle^{\alpha-1/2}|D_{x}|^{a/2}u\in L^{2}([t_{1},t_{2}]:L^{2}(\mathbb{R})).
	\end{equation*}
	This estimate holds at almost every time in $[t_{1},t_{2}]$. Thus, we can choose two times $t_{1}\le t_{3}<t_{4}\le t_{2}$ where it holds, and repeat the argument with $J^{a/2}u$ in place of $u$ and with the weight reduced by $1/2$. Iterating this procedure and using interpolation between weights and derivatives, we get $J^{a\alpha}u(\cdot,t^{*})\in L^{2}(\mathbb{R})$ for some $t^{*}\in(t_{1},t_{2})$, and the local theory completes the proof. The main difficulty is the \emph{non-local character} of $|D_{x}|^{a}$. Unlike the KdV equation, the smoothing effect involves $|D_{x}|^{a/2}$, whose symbol is not smooth at the origin, and the lower order terms above are not local. To deal with this, we combine the commutator expansions of Ginibre and Velo \cite{GV1,GV2}, the fractional Leibniz rules and commutator estimates of D.~Li \cite{Dongli}, and weighted estimates for pseudo-differential operators. The remainders are then controlled, uniformly in $N$, by quantities of the form $\big\|\big(|D_{x}|^{a+1+2\sigma}(x\langle x\rangle_{N}^{2\alpha-1})\big)^{\wedge}\big\|_{L^{1}_{\xi}}$. This restricts the size of the weight at each step of the iteration and explains the different cases in the proof.
	
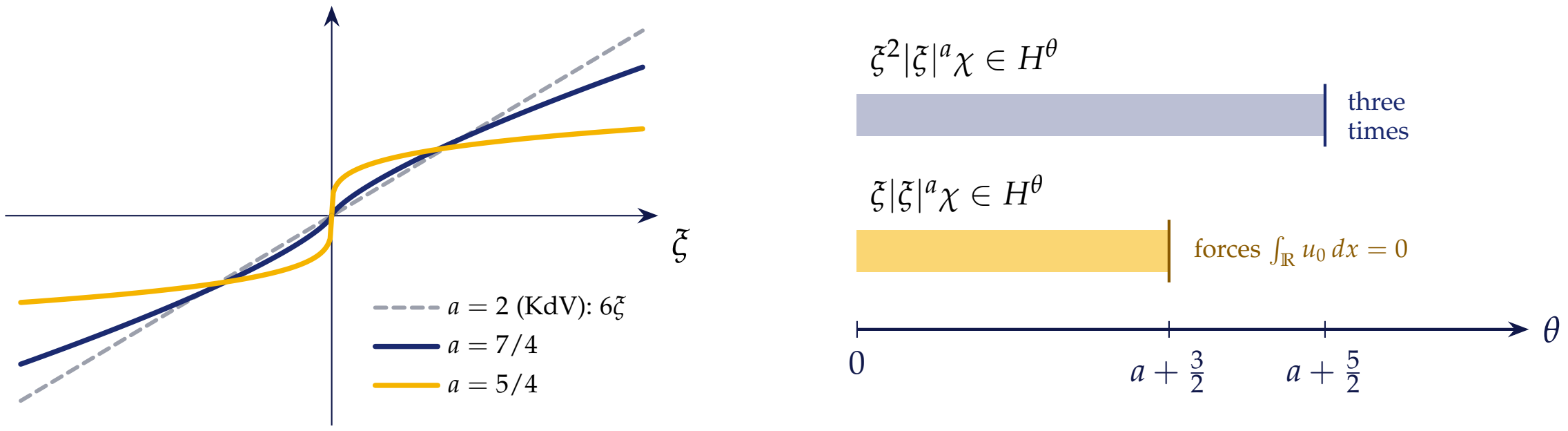
\begin{figure}[ht]
		\centering
		\begin{minipage}[c]{0.49\textwidth}
			\centering
			\begin{tikzpicture}
				\begin{axis}[
					width=\textwidth, height=5.6cm,
					xmin=-1.05, xmax=1.05, ymin=-6.8, ymax=6.8,
					axis lines=middle, axis line style={-{Stealth[length=2mm]}, craneNavy},
					xlabel={$\xi$}, xtick=\empty, ytick=\empty,
					xlabel style={anchor=north west},
					legend style={draw=none, fill=none, font=\scriptsize, at={(0.98,0.03)}, anchor=south east, fill=white, fill opacity=0.9, text opacity=1, cells={anchor=west}},
					every axis plot/.append style={line cap=round}
					]
					\addplot[craneGray, line width=1.1pt, densely dashed, domain=-1:1, samples=2] {6*x};
					\addlegendentry{$a=2$ (KdV): $6\xi$}
					\addplot[craneInk, line width=1.4pt, domain=0:1, samples=201] {1.75*2.75*x^(0.75)};
					\addlegendentry{$a=7/4$}
					\addplot[craneGold, line width=1.4pt, domain=0:1, samples=201] {1.25*2.25*x^(0.25)};
					\addlegendentry{$a=5/4$}
					\addplot[craneInk, line width=1.4pt, domain=0:1, samples=201, forget plot] ({-x},{-1.75*2.75*x^(0.75)});
					\addplot[craneGold, line width=1.4pt, domain=0:1, samples=201, forget plot] ({-x},{-1.25*2.25*x^(0.25)});
				\end{axis}
			\end{tikzpicture}
		\end{minipage}\hfill
		\begin{minipage}[c]{0.49\textwidth}
			\centering
			\begin{tikzpicture}[x=1.15cm, y=1cm]
				\draw[-{Stealth[length=2mm]}, craneNavy, line width=0.8pt] (0,0) -- (5.6,0) node[anchor=west] {$\theta$};
				\foreach \p/\l in {0/0, 2.6/$a+\frac{3}{2}$, 3.9/$a+\frac{5}{2}$}{
					\draw[craneNavy] (\p,0.08) -- (\p,-0.08) node[anchor=north, font=\small] {\l};
				}
				\fill[craneGold!55] (0,0.55) rectangle (2.6,0.95);
				\draw[craneGoldDark, line width=0.8pt] (2.6,0.45) -- (2.6,1.05);
				\node[anchor=west, font=\small] at (0,1.3) {$\xi|\xi|^{a}\chi\in H^{\theta}$};
				\fill[craneInk!30] (0,1.85) rectangle (3.9,2.25);
				\draw[craneInk, line width=0.8pt] (3.9,1.75) -- (3.9,2.35);
				\node[anchor=west, font=\small] at (0,2.6) {$\xi^{2}|\xi|^{a}\chi\in H^{\theta}$};
				\node[anchor=west, font=\scriptsize, craneGoldDark, align=left] at (2.7,0.75) {forces $\int_{\mathbb{R}} u_{0}\,dx=0$};
				\node[anchor=west, font=\scriptsize, craneInk, align=left] at (3.98,2.05) {three\\ times};
			\end{tikzpicture}
		\end{minipage}
		\caption{Left: the function $\partial_{\xi}^{2}(\xi|\xi|^{a})=a(a+1)|\xi|^{a-1}\operatorname{sgn}(\xi)$ for two values of $a\in(1,2)$, compared with the smooth KdV case $a=2$. The vertical tangent at $\xi=0$ reflects the limited smoothness of the symbol. Right: the Sobolev regularity of the first two singular terms in the expansion of $\widehat{u}(\xi,t)$ near the origin, which determines the decay thresholds $3/2+a$ and $5/2+a$ in Theorem~\ref{AA2}.}
		\label{fig:symbol}
	\end{figure}

	The decay thresholds in Theorem~\ref{AA2} are easier to understand in the frequency variable, where the weight $|x|^{\alpha}$ corresponds to $\alpha$ derivatives of $\widehat{u}$. Formally, from the integral equation associated with \eqref{dgbo} we have
	\begin{equation*}
		\widehat{u}(\xi,t)=e^{it\xi|\xi|^{a}}\widehat{u}_{0}(\xi)-\frac{i}{2}\int_{0}^{t}e^{i(t-\tau)\xi|\xi|^{a}}\,\xi\,\widehat{u^{2}}(\xi,\tau)\,d\tau.
	\end{equation*}
	Expanding the exponential near $\xi=0$, the first non-smooth term is $it\,\widehat{u}_{0}(0)\,\xi|\xi|^{a}$. Let $\chi\in C^{\infty}_{0}(\mathbb{R})$ with $\chi=1$ near the origin. Since $\xi|\xi|^{a}\chi(\xi)\in H^{\theta}(\mathbb{R})$ if and only if $\theta<a+3/2$ (see \cite{FONLP}), decay of order $3/2+a$ at two different times implies $\widehat{u}_{0}(0)=0$, which is Theorem~\ref{AA2} (ii). Under the zero mean condition, the next singular term is $\xi^{2}|\xi|^{a}$, which belongs to $H^{\theta}_{\rm loc}$ only for $\theta<a+5/2$. Taking $t_{1}=0$ for simplicity, by \eqref{virial} its coefficient is
	\begin{equation*}
		\int_{0}^{t}\int_{\mathbb{R}}x\,u(x,\tau)\,dx\,d\tau
		=t\int_{\mathbb{R}}x\,u_{0}(x)\,dx+\frac{t^{2}}{4}\|u_{0}\|_{2}^{2}.
	\end{equation*}
	This polynomial in $t$ vanishes only at $t=0$ and at $t^{*}=-4\|u_{0}\|_{2}^{-2}\int_{\mathbb{R}} x\,u_{0}(x)\,dx$. Thus, decay of order $5/2+a$ at three different times implies $\|u_{0}\|_{2}=0$. This heuristic agrees with Theorem~\ref{AA2} (iv) and with the value of $t^{*}$ in Theorem~\ref{TheoremA} (v) and in \cite{FONLP}, which shows that in general two times are not enough. We also notice that the ground states $\phi_{c}$ in \eqref{soliton} are positive, so they have non-zero mean. By the bound above and the expansion of $\widehat{\phi_{c}}$, $\phi_{c}\in L^{2}(|x|^{2\alpha}dx)$ if and only if $\alpha<a+1/2$, which is below the threshold $3/2+a$ in Theorem~\ref{AA2} (ii). Figure~\ref{fig:symbol} illustrates the singularity of the symbol at the origin.
	The proof of parts (i) and (iii) of Theorem~\ref{AA2} is based on an iteration in the size of the weight. We write $\alpha=m+\theta$ with $m\in\mathbb{N}_{0}$ and $0<\theta\le 1$, multiply the equation in \eqref{dgbo} by $x^{2m}\langle x\rangle_{N}^{2\theta}u$, and use the identity
	\begin{equation*}
		x^{m}\partial_{x}|D_{x}|^{a}f-\partial_{x}|D_{x}|^{a}(x^{m}f)=\sum_{k=1}^{m}\binom{m}{k}\widetilde{|D_{x}|^{a+1,k}}\big(x^{m-k}f\big),
	\end{equation*}
	where $\widetilde{|D_{x}|^{a+1,k}}$ denotes the Fourier multiplier with symbol $i^{k}\partial_{\xi}^{k}(i\xi|\xi|^{a})$. The main step is to control the dispersive term with a constant independent of $N$, keeping the regularity at the level $s\ge\max\{a/2,a\alpha\}$. For this we decompose carefully the commutators between $\langle x\rangle_{N}^{\theta}$, $\mathcal{H}$ and $|D_{x}|^{a}$. When $\alpha\ge 3/2+a$, the zero mean condition enters through the cancellation $\widehat{u}(0,t)=0$. Finally, to obtain the unique continuation principles in Theorem~\ref{AA2} (ii) and (iv), we apply negative powers $J^{-m}$, $m>0$, to the integral formulation of \eqref{dgbo}. This allows us to isolate the singular term described above without extra regularity on the solution.
	
	\emph{We remark that, for weights $\alpha>1/2$, the case $s\le 3/2$ is not covered by Theorem~\ref{main}. In our argument the condition $s>3/2$ is used to control $\|\partial_{x}u\|_{L^{\infty}_{x}}$ in the nonlinear part of the weighted energy estimates. This gap may seem to be a consequence of the method. Nevertheless, as observed in \cite{LPBO} for the BO equation, a similar gap appears for the Burgers equation, where the a priori estimate $\sup_{[0,T]}\|u(t)\|_{H^{s}}\le c_{0}\|u_{0}\|_{H^{s}}$ holds only for $s=0$ and $s>3/2$.}
	
	The paper is organized as follows. In Section 2 we establish the preliminary estimates: commutator and weighted estimates, results from pseudo-differential calculus, the commutator expansions of Ginibre and Velo, and some nonlinear weighted estimates. Section 3 is devoted to the proof of Theorem~\ref{main}. In Section 4 we prove Theorem~\ref{AA2}.
	\section{Notation}
	For two operators $A$ and $B$, we define their \emph{commutator} by $[A;B] := AB - BA$. Throughout this document, we  will make extensive use of the \emph{Fourier transform}, which is defined as 
	\begin{equation*}
		\widehat{f}(\xi)
		=
		\int_{\mathbb{R}} e^{-ix\xi} \, f(x)\,dx,
	\end{equation*}
	with the inverse transform given by 
	\begin{equation*}
		g^{\vee}(x)=\frac{1}{2\pi}\int_{\mathbb{R}} e^{ix\xi}\, g(\xi)\,d\xi.
	\end{equation*}
	We denote by $\mathcal{H}$  the \emph{Hilbert transform}, defined by
	\begin{equation}\label{hita}
		\mathcal {H} f(x)=\frac{1}{\pi} {\rm p.v.}\big(\frac{1}{x}\ast f\big)(x)
		=\frac{1}{\pi}\lim_{\epsilon\downarrow 0}\int\limits_{|y|\ge \epsilon} \frac{f(x-y)}{y}\,dy=(-i\,\operatorname{sgn}(\xi) \widehat{f}(\xi))^{\vee}(x).
	\end{equation}
	Given $s\in \mathbb{R}$, we define the \emph{Bessel potential operator} $J^s$ through its Fourier multiplier $\widehat{J^s f}(\xi)=(1+|\xi|^2)^{s/2}\widehat{f}(\xi)$, whenever this expression is well-defined. 
	
	We denote by $C_0^{\infty}(\mathbb{R})$ the space of smooth, compactly supported functions. $\mathcal{S}(\mathbb{R})$ denotes the \emph{Schwartz space} of rapidly decreasing smooth functions and $\mathcal{S}'(\mathbb{R})$ denotes its dual space of tempered distributions. For $1 \leq p \leq \infty$, $L^p(\mathbb{R})$ denotes the usual
	\emph{Lebesgue space} with norm $\|\cdot\|_{p}$. For $s \in \mathbb{R}$, we denote by $H^s(\mathbb{R})$ the standard $L^2$-based \emph{Sobolev space} of order $s$, equipped with its usual norm $\|f\|_{H^s}=\|J^s f\|_2$.
	
	\maketitle

	\maketitle
	\section{Preliminary Estimates}	

	In this section, we establish the preliminary estimates needed to deduce our main results. We begin with commutator and weighted estimates, followed by several results from pseudo-differential calculus. We then develop further commutator expansions and conclude by establishing general nonlinear estimates.

	We start by collecting several commutator estimates involving fractional derivatives, the Hilbert transform and multiplication by smooth functions. These estimates will be used repeatedly to localize the non-local operator $|D_x|^{a}$ in the weighted energy estimates of Sections 3 and 4.
	
	\begin{lem}\label{lowcommu}
		Let $0\le \alpha<1$, $0<\beta\le 1-\alpha$, and $1<p<\infty$. Then
		\begin{equation*}
			\bigl\|
			|D_{x}|^{\alpha+\beta}\bigl(a |D_{x}|^{1-(\alpha+\beta)}f\bigr)
			-
			|D_{x}|^{\alpha}\bigl(a |D_{x}|^{1-\alpha}f\bigr)
			\bigr\|_{p}
			\leq c_{\alpha,\beta,p}
			\|\partial_x a\|_{\infty}
			\|f\|_{p}.
		\end{equation*}
	\end{lem}
	\begin{proof}
		The estimate was established in \cite{Dongli} as an improvement of one found in  \cite{DMP}. 
	\end{proof}
	The next lemma provides the first terms in the expansion of $|D_x|^{s}(fg)$ for $1<s<3$, with remainders controlled by $\||D_x|^{s}f\|_{\infty}$. This will be particularly useful when $f$ is a weight.
	\begin{lem}\label{extendenlemma}
		Let $1<p<\infty$. Then
		\begin{equation}\label{eqextenden1}
			\||D_x|^{s}(fg)-f|D_x|^{s}g-s\partial_x f \mathcal{H}|D_x|^{s-1} g\|_{p}\leq c_{s,p} \||D_x|^{s}f\|_{\infty}\|g\|_{p},
		\end{equation}
		if $1<s<2$. Moreover, it holds
		\begin{equation}\label{eqextenden2}
			\begin{aligned}
				\||D_x|^{s}(fg)-f|D_x|^{s}g-s\partial_x f \mathcal{H}|D_x|^{s-1}g+&\frac{s(s-1)}{2} \partial_x^2 f|D_x|^{s-2} g\|_{p}\leq c_{s,p} \||D_x|^{s}f\|_{\infty}\|g\|_{p},
			\end{aligned}
		\end{equation}
		whenever $2\leq s <3$.
	\end{lem}
	\begin{proof}
		The estimates are the one-dimensional version of inequality (1.6) in \cite[Theorem 1.2]{Dongli}  with $s_1=s$ and $s_2=0$.
	\end{proof}
	Next, we recall a commutator estimate established in \cite{KATOP2}, and an extension given in  \cite{GO}.
	\begin{thm} Let $s>0$, $1<p<\infty$, and  $1/p=1/p_1+1/p_2=1/p_3+1/p_4,$ with $p_2, p_3<\infty$. Then there exists $c=c_{s,p,p_1,p_2,p_3,p_4}>0$ such that for all $f, g\in \mathcal{S}(\mathbb{R})$, it holds that
		\begin{equation}\label{kp-comm}
			\|[J^s; f]g\|_p\le c \big(\|\partial_x f\|_{p_1}\|J^{s-1}g\|_{p_2}+\|J^sf\|_{p_3}\|g\|_{p_4}\big).
		\end{equation} 
		Moreover, one has 
		\begin{equation}\label{kp-comm1}
			\|[J^s;f]g\|_p\le c (\|\partial_x f\|_{\infty}\|J^{s-1}g\|_{p}+\|J^{s-1}\partial_xf\|_{p}\|g\|_{\infty}).
		\end{equation}
		Additionally, $s\in(0,1]$, then
		\begin{equation}\label{kp-comm2}
			\|[J^s; f]g\|_p\le c (\|J^{s-1}\partial_xf\|_{p}\|g\|_{\infty}).
		\end{equation}
	\end{thm}
	\medskip
	In the case $p=p_2=p_3$, the proof of \eqref{kp-comm} was given in \cite{KATOP2}. The general case in \eqref{kp-comm} follows the same argument using the version of the multiplier result of Coifman-Meyer \cite{CoMe} stated in  \cite{GO} (Theorem A). The estimate \eqref{kp-comm1} was proved in \cite{Dongli} answering a question raised in \cite{KATOP2}.
	We shall also need the fractional Leibniz rule in the scale of the Bessel potentials.
	\begin{lem}\label{lem1}
		Let $p\in[1,\infty]$ and $s\geq 0$, with $s$ not an odd integer in the case $p=1$. Then 
		\begin{equation}\label{comm-go-1}
			\|J^{s} (fg)\|_p \leq c (\|f\|_{p_1}\|J^{s}g\|_{p_2} +\|g\|_{q_1}\|J^{s}f\|_{q_2}),
		\end{equation}
		with $1/p_1+1/p_2=1/q_1+1/q_2=1/p$.
	\end{lem}
	For a proof of Lemma \ref{lem1} with $p\in(1,\infty)$ we refer to  \cite{KATOP2}, \cite{KPV1}, and \cite{GO}. 
	The case $r=p_1=p_2=q_1=q_2=\infty$ was established in \cite{BLi}, see also \cite{GMN}. 
	The case $p=p_2=q_2=1$ was settled in \cite{OW}.

	Additionally we have some refinement of the commutator estimates found in \cite{KATOP2} and \cite{KPV1}.
	\begin{thm}\label{kprefi}
		Let $s > 0$ and $1 < p < \infty$. Then for any $f, g \in \mathcal{S}(\mathbb{R})$,
		\begin{equation*}
			\| |D_{x}|^{s}(fg) - f |D_{x}|^{s} g \|_{p}
			\leq c_{s,p_1,p_2,p_3,p_4}\big(
			\| |D_{x}|^{s} f \|_{p_1} \| g \|_{p_{2}}
			+
			\| \partial_{x} f \|_{p_3} \| |D_{x}|^{s-1} g \|_{p_{4}}\big),
		\end{equation*}
		where $1 < p_{1}, p_{2}, p_{3}, p_{4} \le \infty$, and
		\begin{equation*}
			\frac{1}{p_{1}} + \frac{1}{p_{2}}
			=
			\frac{1}{p_{3}} + \frac{1}{p_{4}}
			=
			\frac{1}{p}.
		\end{equation*}
		For $0 < s < 1$ and $1 < p < \infty$, the following inequality holds:
		\begin{equation*}
			\| |D_{x}|^{s}(fg) - f |D_{x}|^{s} g \|_{p}
			\leq c_{s,p_1,p_2}
			\| |D_{x}|^{s} f \|_{p_1} \| g \|_{p_2},
		\end{equation*}
		where
		\begin{equation*}
			\frac{1}{p_{1}} + \frac{1}{p_{2}} = \frac{1}{p},
			\qquad
			1 < p_{1}, p_{2} \le \infty.
		\end{equation*}
		For $s = 1$ and $1 < p < \infty$,
		\begin{equation*}
			\| |D_{x}|(fg) - f |D_{x}| g \|_{p}
			\leq c_{s,p_1,p_2}
			\| \partial_{x} f \|_{p_1} \| g \|_{p_2},
		\end{equation*}
		where
		\begin{equation*}
			\frac{1}{p_{1}} + \frac{1}{p_{2}} = \frac{1}{p},
			\qquad
			1 < p_{1}, p_{2} \le \infty.
		\end{equation*}
		For $p_{1} = \infty$, $p_{2} = p$ and $1 < p < \infty$, we have the following BMO-refinements:
		\begin{itemize}
			\item If $0 < s < 1$, then
			\begin{equation*}
				\| |D_{x}|^{s}(fg) - f |D_{x}|^{s} g \|_{p}
				\leq c_s
				\| |D_{x}|^{s} f \|_{\mathrm{BMO}} \| g \|_{p.}
			\end{equation*}
			\item If $s \ge 1$, then
			\begin{equation*}
				\| |D_{x}|^{s}(fg) - f |D_{x}|^{s} g \|_{p}
				\leq c_{s,p,p_3,p_4}\big(
				\| \partial_{x} f \|_{p_3} \| |D_{x}|^{s-1} g \|_{p_{4}}
				+
				\| |D_{x}|^{s} f \|_{\mathrm{BMO}} \| g \|_{p}\big),
			\end{equation*}
			where
			\begin{equation*}
				\frac{1}{p_{3}} + \frac{1}{p_{4}} = \frac{1}{p},
				\qquad
				1 < p_{3} \le \infty,
				\qquad
				1 < p_{4} < \infty.
			\end{equation*}
			\item Also for $s \ge 1$,
			\begin{equation*}
				\| |D_{x}|^{s}(fg) - f |D_{x}|^{s} g \|_{p}
				\leq c_{s,p}\big(
				\| \partial_{x} f \|_{p} \| |D_{x}|^{s-1} g \|_{\mathrm{BMO}}
				+
				\| |D_{x}|^{s} f \|_{\mathrm{BMO}} \| g \|_{p}\big).
			\end{equation*}
		\end{itemize}
	\end{thm}
	Since the dispersive operator can be written as $\partial_x|D_x|^{a}=-\mathcal{H}|D_x|^{a+1}$, commutators between the Hilbert transform and smooth functions appear naturally in our arguments.
	Next, we recall a commutator estimate for the Hilbert transform, see \eqref{hita}.
	\begin{thm}\label{calderon1}
		For any $j, k \in \mathbb N\cup\{0\},\,j+k\geq 1$ and any $p\in(1,\infty)$
		\begin{equation} \label{calderon}
			\| \partial_x^j \big[\mathcal H;\eta\big]\partial_x^kf\|_p\leq c_{p,j,k}\|\eta^{(j+k)}\|_{\infty} \|f\|_p.
		\end{equation}
	\end{thm}
	\begin{proof}
		The case $j+k=1$ corresponds to Calder\'on first commutator estimate \cite{calderon}. The general case of \eqref{I1} was established in \cite{BC}.  For a different proof see \cite{DMP}.
	\end{proof}
	We also require fractional versions of the previous estimate, in which the derivatives are replaced by homogeneous fractional derivatives.
	\begin{lem}\label{calderon7}
		For any $a,b >0$
		\begin{equation} \label{calderon2}
			\| |D_x|^a \big[\mathcal H;\varphi\big] |D_x|^b f\|_2\leq c_{a,b}\| \widehat{|D_x|^{a+b}\varphi}\|_{L^1_{\xi}} \|f\|_2.
		\end{equation}
		Let $1<p<\infty$. If $0\leq a,b \leq 1$, $b>0$ with $a+b= 1$, then
		\begin{equation}\label{fracComwell}
			\||D_x|^{a}[\mathcal{H};g]|D_x|^{b}f\|_{p}\leq c_{a,b} \|\partial_x g\|_{\infty}\|f\|_{p}.
		\end{equation}
		Additionally, for any $b>0$,
		\begin{equation}\label{fracComwell2}
			\|[\mathcal{H};g]|D_x|^{b}f\|_{p} \leq c_{b} \||D_x|^{b} g\|_{\infty}\|f\|_{p}.
		\end{equation}
	\end{lem}
	\begin{proof}
		The proof of \eqref{fracComwell} can be found in \cite{RianoCap2021}, while \eqref{fracComwell2} is obtained in \cite[Proposition 2.4]{Riano2021}. To deduce \eqref{calderon2}, we have that
		\begin{equation*}
			\widehat {(|D_x|^a \big[\mathcal H;\varphi\big] |D_x|^b f})(\xi)=i \int_{\mathbb{R}} |\xi|^a|\eta|^b(\operatorname{sgn}(\xi)-\operatorname{sgn}(\eta))\widehat{\varphi}(\xi-\eta)\widehat{f}(\eta)\,d\eta.
		\end{equation*}
		Since 
		\begin{equation*}
			\operatorname{sgn}(\xi)-\operatorname{sgn}(\eta)\neq 0 \iff \eta\,\xi<0\implies |\xi-\eta|=|\xi|+|\eta|
		\end{equation*}
		one has  in this set that 
		\begin{equation*}
			|\xi|^a|\eta|^b\leq (|\xi|+|\eta|)^{a+b}=|\xi-\eta|^{a+b}.
		\end{equation*}
		Thus,
		\begin{equation*}
			|\widehat {(|D_x|^a \big[\mathcal H;\eta\big] |D_x|^b f})(\xi)|\leq 2\int_{\mathbb{R}} |\xi-\eta|^{a+b}|\widehat{\varphi}(\xi-\eta)||\widehat{f}(\eta)|\,d\eta,
		\end{equation*}
		therefore Young's inequality yields the desired result.
	\end{proof}
	Next, we turn to the relation between weights and derivatives. The interpolation inequalities below allow us to exchange decay and regularity, and they are the basic tool to deduce regularity from the weighted estimates obtained in Section 3.
	We recall the following interpolation results, which will be used throughout the paper.
	\begin{lem}\label{lem2-interpolation} Let $p\in (1,\infty)$. Let $a, b\in\mathbb{R}$. Then for any $\theta\in (0,1)$,
		\begin{equation*}\label{interpol-1}
			\| \langle x\rangle^{\theta b}J^{(1-\theta)a} f\|_p \le c_p\|\langle x\rangle^{b}f\|_p^{\theta}\|J^a f\|_p^{1-\theta}.
		\end{equation*}
		and
		\begin{equation}\label{inter1x}
			\|J^{a \theta }(\langle x \rangle^{(1-\theta)b}f)\|_{p}\leq c_2 \|\langle  x \rangle^{b}f\|_{p}^{1-\beta}\|J^{a}f\|_{p}^{\beta}.
		\end{equation}
	\end{lem}
	\begin{rem} The case $p=2$ and $a, b>0$ was proven in \cite{NP}.
	\end{rem}
	We shall also use the following generalization:
	\begin{lem}\label{lem2b-interpolation} Let $p\in (1,\infty)$. Let $a,b,c,d\in\mathbb{R}$. Then for $\theta\in[0,1]$ it holds that
		\begin{equation}\label{interpol-complete}
			\|\langle x\rangle^{\theta a+(1-\theta)c}J^{\theta b+(1-\theta)d} f\|_p \le c_p\|\langle x\rangle^{a}J^{b} f\|_p^{\theta}\|\langle x\rangle^{c}J^{d} f\|_p^{(1-\theta)}.
		\end{equation}
	\end{lem}
	See \cite{LPILW} for a proof.
	When the weight satisfies the Muckenhoupt $A_2$ condition, the interpolation can be performed with homogeneous derivatives.
	\begin{prop}\label{Dint}  Let $\alpha\in(-1/2,1/2)$ and $b>0$, then
		\begin{equation}\label{Dint1}
			\|\langle x \rangle^{\theta \alpha}|D_x|^{(1-\theta)b}f\|_2\leq c\|\langle x\rangle^{\alpha}f\|_2^{\theta} \,\||D_x|^bf\|_2^{1-\theta}.
		\end{equation}
	\end{prop}
	The proof of \eqref{Dint1} follows the Three Line Theorem argument and the fact that  $\langle x\rangle^\beta$ and $|x|^\beta$ are $A_2$ weights for $\beta\in(-1,1)$, that is, they satisfy the Muckenhoupt $A_2$ condition. See, for instance, \cite{HMW,Petermichl}.
	The following elementary lemma will be used to exploit the zero mean condition $\widehat{u}(0)=0$ in the frequency variable.
	\begin{lem}\label{Hardy-type} Let $f\in H^s(\mathbb{R})$, $s>1/2$, with compact support and $f(0)=0$. Then there exists $g\in H^{s-1}(\mathbb{R})$ such that
		\begin{equation*}
			f(x)=xg(x)\hskip10pt\text{and}\hskip10pt \|g\|_{H^{s-1}}\leq c\|f\|_{H^s}.
		\end{equation*}
	\end{lem}
	\begin{proof}
		Let
		\begin{equation*}
			f(x)=\int_0^xf'(t)\,dt=x\,\int_0^1f'(\theta x)\,d\theta\equiv x\,g(x).
		\end{equation*}
		Thus, 
		\begin{equation*}
			\|g\|_{H^{s-1}}\leq c\Big\|\frac{f}{x}\Big\|_{H^{s-1}}\leq c\|f'\|_{H^{s-1}}\leq c\|f\|_{H^s}.
		\end{equation*}
	\end{proof}
	We conclude this section by presenting the following identity for the commutator of $x^m$ and $|D_x|^{a}$.
	\begin{prop}\label{GeneralConmutator}
		Let $a>0$, $m\in \mathbb{N}$, then for an appropriated function $f$ the following identity holds 
		\begin{equation}\label{generalconmmutator}
			x^m\partial_x |D_x|^{a}f-\partial_x |D_x|^{a}(x^mf)=\sum_{k=1}^m \binom{m}{k} \widetilde{|D_x|^{a+1,k}}\big(x^{m-k}f\big),    
		\end{equation}
		where each $\widetilde{|D_x|^{a+1,k}}$ is defined via the Fourier multiplier by the function $i^k\, \partial_{\xi}^{k}(i\xi|\xi|^{a})$. In particular,
		\begin{equation}\label{conmderivIdent}
			\begin{aligned}
				&\widetilde{|D_x|^{a+1,1}}=-(a+1)|D_x|^{a},\\
				&\widetilde{|D_x|^{a+1,2}}=-a(a+1)\partial_x |D_x|^{a-2}=a(a+1)\mathcal{H}|D_x|^{a-1}, \\
				&\widetilde{|D_x|^{a+1,3}}=a(a+1)(a-1)|D_x|^{a-2}, \\
				&\widetilde{|D_x|^{a+1,4}}=a(a+1)(a-1)(a-2)\partial_x |D_x|^{a-4}=-a(a+1)(a-1)(a-2)\mathcal{H} |D_x|^{a-3}.
			\end{aligned}    
		\end{equation}
		Moreover, for each $k\geq 1$, it follows
		\begin{equation}\label{partexampleconmm}
			\begin{aligned}
				\, [x,\widetilde{|D_x|^{a+1,k}}]=\widetilde{|D_x|^{a+1,k+1}}.  
			\end{aligned}    
		\end{equation}
	\end{prop}
	\begin{proof}
		The proof follows by observing that
		\begin{equation*}
			\mathcal{F}\bigl(x^m\partial_x|D_x|^{a} f\bigr)(\xi)
			=i^m\partial_\xi^m
			\bigl(i\xi|\xi|^{a}\widehat f(\xi)\bigr),
		\end{equation*}
		and then applying the Leibniz rule to the right-hand side.    
	\end{proof}	
	\subsection{Pseudo-differential Calculus and Applications}\hskip15pt

	In this subsection we recall some basic facts about pseudo-differential operators whose symbols involve weights in both variables, and we deduce some consequences relating homogeneous and non-homogeneous derivatives in weighted spaces.

	For $m,q\in\mathbb{R}$, the class of symbols $\mathcal{S}^{m,q}=\mathcal{S}^{m,q}(\mathbb{R}^2)$ consists of all $a\in C^\infty(\mathbb{R}_x\times\mathbb{R}_\xi)$ such that, for every $k,\ell\in\mathbb{N}_0$, there is a constant $c_{k,\ell}>0$ for which
	\begin{equation*}
		\bigl|\partial_x^k\partial_\xi^\ell a(x,\xi)\bigr|
		\leq c_{k,\ell}\langle x\rangle^{q-k}
		\langle\xi\rangle^{m-\ell},
		\qquad (x,\xi)\in\mathbb{R}^2.
	\end{equation*}
	If $a\in\mathcal{S}^{m,q}$ and $f\in\mathcal{S}(\mathbb{R})$, we set
	\begin{equation*}
		T_af(x)
		=\frac{1}{2\pi}\int_{\mathbb{R}}
		e^{ix\xi}a(x,\xi)\widehat f(\xi)\,d\xi.
	\end{equation*}
	We write $T_a\in\operatorname{Op}(\mathcal{S}^{m,q})$.
	The following standard properties will be used below. 
	\begin{lem}\label{LemSGcalculus}
		Let $m,m_1,m_2,q,q_1,q_2\in\mathbb{R}$.
		\begin{enumerate}
			\item[(i)] If $m_1\leq m_2$ and $q_1\leq q_2$, then $\mathcal{S}^{m_1,q_1}\subseteq \mathcal{S}^{m_2,q_2}$.
			
			\item[(ii)] If $a\in\mathcal{S}^{0,0}$, then
			$T_a$ extends uniquely to a bounded operator on
			$L^p(\mathbb{R})$ for every $1<p<\infty$.
			
			\item[(iii)] If $1<p<\infty$ and $a\in\mathcal{S}^{m,q}$, then
			\begin{equation*}
				\|T_af\|_{p} \leq  c_a \|\langle x\rangle^qJ^m f\|_{p},   \qquad f\in\mathcal{S}(\mathbb{R}).
			\end{equation*}
			
			\item[(iv)] If $a\in\mathcal{S}^{m_1,q_1}$ and
			$b\in\mathcal{S}^{m_2,q_2}$, then there exists $c=a\# b\in\mathcal{S}^{m_1+m_2,q_1+q_2}$ such that
			\begin{equation*}T_aT_b=T_c.\end{equation*}
			Moreover,
			\begin{equation*}
				a\# b\sim
				\sum_{j=0}^{\infty}\frac{(-i)^j}{j!}
				(\partial_\xi^j a)(\partial_x^j b),
			\end{equation*}
			in the sense that, for every integer $N\geq1$,
			\begin{equation*}a\# b-\sum_{j=0}^{N-1}\frac{(-i)^j}{j!}(\partial_\xi^j a)(\partial_x^j b)\in\mathcal{S}^{m_1+m_2-N,q_1+q_2-N}.\end{equation*}
			
			\item[(v)] Under the assumptions in \emph{(iv)}, the commutator belongs
			to the class
			\begin{equation*}[T_a;T_b] \in\operatorname{Op} (\mathcal{S}^{m_1+m_2-1,q_1+q_2-1}).\end{equation*}
			More precisely, if $d=a\# b-b\# a$ denotes its symbol, then
			\begin{equation*}
				d\sim\sum_{j=1}^{\infty}\frac{(-i)^j}{j!}
				\left(
				\partial_\xi^j a\,\partial_x^j b
				-\partial_\xi^j b\,\partial_x^j a
				\right).
			\end{equation*}
			Thus, for every integer $N\geq2$,
			\begin{equation*}
				\begin{aligned}
					d-\sum_{j=1}^{N-1}\frac{(-i)^j}{j!}
					\left(
					\partial_\xi^j a\,\partial_x^j b
					-\partial_\xi^j b\,\partial_x^j a
					\right) \in \mathcal{S}^{m_1+m_2-N,q_1+q_2-N}.
				\end{aligned}
			\end{equation*}
			In particular, its principal symbol is
			\begin{equation*}
				d=\frac{1}{i}
				\bigl(\partial_\xi a\,\partial_x b
				-\partial_x a\,\partial_\xi b\bigr)
				\quad\bmod
				\mathcal{S}^{m_1+m_2-2,q_1+q_2-2}.
			\end{equation*}
		\end{enumerate}
	\end{lem}
	\begin{proof}
		The properties follow from the classical $S^0_{1,0}$ theory. For the $L^p$ boundedness in (iii), we refer to \cite{LPILW}, see also, \cite{NP,stein2,MendezRiano2025,KPV98}.  
	\end{proof} 
	
	\begin{rem}
		If a symbol $a=a(x)$ is independent of $\xi$, then
		$T_a$ is simply multiplication by $a(x)$. In this case, we will use the same notation for the symbol and the corresponding operator $T_a$.
	\end{rem}
	As a consequence of the above result, we have: 
	\medskip
	\begin{cor}\label{pseudo-cor1} Let $1<p<\infty$, and $a,b>0$.
		If
		\begin{equation*}
			\la x\ra^b |D_x|^a f\in L^p(\mathbb{R}) \hskip15pt\text{and}\hskip15pt \la x\ra^b  f\in L^2(\mathbb{R}),
		\end{equation*}
		then
		\begin{equation*}
			\la x\ra^b J^a f\in L^p(\mathbb{R})
		\end{equation*}
		with
		\begin{equation*}
			\|\la x\ra^b J^a f\|_p
			\leq
			c\big(\|\la x\ra^b |D_x|^a f\|_p+\|\la x\ra^b  f\|_{p}).
		\end{equation*}
	\end{cor}
	\begin{proof}
		Let $\chi\in C_0^\infty(\mathbb{R})$ with $\chi(x)=1$ if $|x|<1$
		and $\chi(x)=0$ if $|x|>2$. Then
		\begin{equation*}
			J^a f
			=
			J^a\chi(D_x)f
			+
			J^a(1-\chi(D_x))f.
		\end{equation*}
		Since $(1+\xi^2)^{a/2}\chi(\xi)\in C_0^\infty(\mathbb{R})$ one has that
		\begin{equation*}
			J^a\chi(D_x)f=K_a*f,
			\qquad
			\text{with }\;\;\;
			K_a(x)=\big((1+\xi^2)^{a/2}\chi(\xi)\big)^{\vee}(x).
		\end{equation*}
		Hence, $K_a\in\mathcal{S}(\mathbb{R})$,
		$\la x\ra^bK_a\in L^1(\mathbb{R})$ and
		\begin{equation*}
			\|\la x\ra^b(K_a*f)\|_p
			\leq
			c\,\|\la x\ra^bK_a\|_1\,\|\la x\ra^bf\|_p,
		\end{equation*}
		which follows from $\la x\ra^b\leq c_b\la x-y\ra^b\la y\ra^b$. Also, since
		\begin{equation*}
			J^a(1-\chi(D_x))|D_x|^{-a}\in \operatorname{Op}(\mathcal{S}^{0,0}),
		\end{equation*}
		by (iii) in Lemma \ref{LemSGcalculus}, it follows that
		\begin{equation*}
			\|\la x\ra^bJ^a(1-\chi(D_x))f\|_p
			\leq
			c\,\|\la x\ra^b|D_x|^af\|_p,
		\end{equation*}
		which yields the desired result.
	\end{proof}
	\medskip
	The next result is a version of Corollary \ref{pseudo-cor1} for negative powers of the weight.
	\begin{cor}\label{pseudo-cor2} If $b<0$ and $a\in(0,2]$, then 
		\begin{equation*}
			\|\la x\ra^b J^a f\|_2
			\leq
			c\left(
			\|\la x\ra^b |D_x|^a f\|_2+\| f\|_2
			\right).
		\end{equation*}
	\end{cor}
	\begin{proof}
		Since 
		\begin{equation*}
			J^af=|D_x|^af+(J^a-|D_x|^a)f=|D_x|^af+T_{m_a}f,
		\end{equation*}
		with
		\begin{equation*}
			m_a(\xi)=(1+\xi^2)^{a/2}-|\xi|^a\in L^{\infty}(\mathbb{R}),\;\;\;\;a\in(0,2],
		\end{equation*}
		one has that
		\begin{equation*}
			\begin{split}
				\|\la x\ra^b J^a f\|_2&\leq c(\|\la x\ra^b |D_x|^a f\|_2+\|\la x\ra^b T_{m_a} f\|_2)\\
				&\leq c(\|\la x\ra^b |D_x|^a f\|_2+\| T_{m_a}f\|_2)\\
				&\leq c(\|\la x\ra^b |D_x|^a f\|_2+\| f\|_2).
			\end{split}
		\end{equation*}
		The proof is complete.
	\end{proof}
	We now estimate the commutator between the weight $\langle x\rangle^{b}$ and the dispersive operator $\partial_x|D_x|^{a}$. The lack of smoothness of the symbol $i\xi|\xi|^{a}$ at the origin limits the size of $b$, and, as in Theorem \ref{AA2}, the zero mean condition allows us to go beyond the threshold $a+3/2$.
	\begin{cor}\label{pseudo-cor3}  Let $a\in(0,2]$. 
		
		(i) If $b\in (0,a +3/2)$, then
		\begin{equation}\label{ee1}
			\|[\langle x\rangle^b;\partial_x|D_x|^a]f\|_2
			\leq c_{a,b} (
			\|\langle x\rangle^{b-1}J^af\|_2+\|\langle x\rangle^bf\|_2).
		\end{equation}
		
		(ii) If $\int_{\mathbb{R}} f(x)dx=0$, then the inequality \eqref{ee1} extends to $b\in[a+3/2,a+5/2)$.
	\end{cor}
	\begin{proof}
		Let  $\chi\in C_0^\infty(\mathbb{R})$ as in the proof of Corollary \ref{pseudo-cor1}.
		We write 
		\begin{equation*}
			\partial_x|D_x|^a=\partial_x|D_x|^a\chi(D_x)+\partial_x|D_x|^a(1-\chi(D_x))=T_0+T_\infty,
		\end{equation*}
		that is
		\begin{equation*}
			\widehat{T_0f}(\xi)=i\xi|\xi|^a\chi(\xi)\widehat f(\xi), \;\;\text{and}\;\;\;
			\widehat{T_\infty f}(\xi)=i\xi|\xi|^a(1-\chi(\xi))\widehat f(\xi). 
		\end{equation*}
		Thus
		\begin{equation*}[\langle x\rangle^b;\partial_x|D_x|^a]
			=[\langle x\rangle^b;T_\infty]f+[\langle x\rangle^b;T_0]f.
		\end{equation*}
		Since $[\langle x \rangle^b;T_{\infty}]\in \operatorname{Op}(\mathcal{S}^{a,b-1})$, using (iii) in Lemma \ref{LemSGcalculus}, one has that
		\begin{equation}\label{ee2a}
			\|[\langle x\rangle^b;T_\infty]f\|_2\le c_{a,b} \|\langle x\rangle^{b-1}J^af\|_2. 
		\end{equation}
		Next, we observe that 
		\begin{equation*}
			T_0f=K_a*f\hskip15pt\text{with}\hskip10pt K_a(x)=(i\xi|\xi|^a\chi(\xi))^\vee(x),
		\end{equation*}
		with the kernel $K_a$ being compact supported and smooth except at $\xi=0$, where  it behaves like $i\xi|\xi|^a$.
		The inverse Fourier transform of this homogeneous expression has degree
		$-(a+2)$. Thus, 
		\begin{equation}\label{ee2}
			|K_a(x)|\le c_a\langle x\rangle^{-a-2}\quad \text{and} \quad
			|K_a'(x)|\le c_a\langle x\rangle^{-a-3}.
		\end{equation}
		Hence,
		\begin{equation*}
			K_a\in L^1(\mathbb{R}),
		\end{equation*}
		and
		\begin{equation}\label{ee3}
			[\langle x\rangle^b;T_0]f=\langle x\rangle^b(K_a*f)-K_a*(\langle x\rangle^bf). 
		\end{equation}
		By Young's inequality,
		\begin{equation*}
			\|K_a*(\langle x\rangle^bf)\|_2\le c_a\|\langle x\rangle^bf\|_2.
		\end{equation*}
		\vskip.2in 
		Thus, it only remains to show  
		\begin{equation}\label{ee5}
			\|\langle x\rangle^b(K_a*f)\|_2\leq c_{a,b} \|\langle x\rangle^b f\|_2.
		\end{equation}
		Notice that until this point we have not used the hypothesis on $b$.
		\vskip.2in 
		First, we assume : \underline{$0<b<a+1$}. Since
		\begin{equation*}
			\langle x\rangle^b\le c_b\langle x-y\rangle^b\langle y\rangle^b. 
		\end{equation*}
		it follows that
		\begin{equation*}
			\langle x\rangle^b|(K_a*f)(x)|
			\le c_b\int_{\mathbb{R}} \langle x-y\rangle^b|K_a(x-y)|
			\langle y\rangle^b|f(y)|\,dy.
		\end{equation*}
		and consequently
		\begin{equation*}
			\|\langle x\rangle^bK_a*f\|_2\le
			c_b\|\langle x\rangle^bK_a\|_1\|\langle x\rangle^bf\|_2. 
		\end{equation*}
		By \eqref{ee2} 
		\begin{equation*}
			\langle x\rangle^b|K_a(x)|\leq c\langle x\rangle^{b-a-2},
		\end{equation*}
		which is integrable when $b<a+1$. This yields \eqref{ee5}  in the case $b\in(0,a+1)$.
		\vskip.1in
		Now, assume : \underline{$a+1\le b<a+3/2$}.  Since in this case  $b>1/2$, one has that
		\begin{equation*}
			\|f\|_1
			\le c_b\| \langle x\rangle^{b}f\|_2.
		\end{equation*}
		For $|x|\le2$,
		\begin{equation*}
			\|\langle x\rangle^{b}K_a*f\|_{L^2(|x|\le2)}
			\le c\|K_a*f\|_2\le c\|f\|_2\le c\|\langle x\rangle^{b}f\|_2. 
		\end{equation*}
		For $|x|>2$, we write 
		\begin{equation*}
			K_a*f(x)=I_1(x)+I_2(x),
		\end{equation*}
		with 
		\begin{equation*}
			I_1(x)=\int_{|y|\le|x|/2}K_a(x-y)f(y)\,dy\;\;\;\;\text{and}\;\;\;\;
			I_2(x)=\int_{|y|>|x|/2}K_a(x-y)f(y)\,dy.
		\end{equation*}
		For $I_1$, since $|x|>2$ it follows that  $|x-y|\ge|x|/2$ and from \eqref{ee2}
		\begin{equation*}
			|I_1(x)|\le c_a|x|^{-a-2}\|f\|_1.
		\end{equation*}
		Thus
		\begin{equation*}
			\langle x\rangle^{b}|I_1(x)|\le
			c_{a,b}|x|^{b-a-2}\|\langle x\rangle^{b}f\|_2. 
		\end{equation*}
		The function $|x|^{b-a-2}$ belongs to $L^2(|x|>2)$ exactly when
		\begin{equation*}
			b<a+3/2. 
		\end{equation*}
		Therefore
		\begin{equation*}
			\|\langle x\rangle^{b}I_1\|_2\le c_{a,b}\|\langle x\rangle^{b}f\|_2.
		\end{equation*}
		For $I_2$, since $|y|>|x|/2$, we have
		$\langle x\rangle^b\le c_b\langle y\rangle^b$, so
		\begin{equation*}
			\langle x\rangle^{b}|I_2(x)|
			\le c_b\bigl(|K_a|*(\langle x\rangle^{b}|f|)\bigr)(x).
		\end{equation*}
		Since $K_a\in L^1$,
		\begin{equation*}
			\|\langle x\rangle^{b}I_2\|_2\le c_{a,b}\|\langle x\rangle^{b}f\|_2.
		\end{equation*}
		A combination of  the above estimates yields \eqref{ee5} in the case
		$0<b<a+3/2$. Hence, we have shown \eqref{ee1} whenever $0<b<a+3/2$, with no condition on $\int_{\mathbb{R}} f\,dx$.
		
		\vskip.2in
		
		Next, we turn to the proof of part (ii) of Corollary \ref{pseudo-cor2}, i.e.,  for : \underline{$b\in[a+3/2, a+5/2)$} with $\int_{\mathbb{R}} f\,dx=0$.
		
		The cancellation
		\begin{equation*}
			\int_{\mathbb{R}}f(y)\,dy=0 
		\end{equation*}
		allows us to write
		\begin{equation*}
			K_a*f(x)=\int_{\mathbb{R}}\bigl(K_a(x-y)-K_a(x)\bigr)f(y)\,dy.
		\end{equation*}
		For $|x|\le2$,
		\begin{equation*}
			\|\langle x\rangle^{b}K_a*f\|_{L^2(|x|\le2)}\le c\|\langle x\rangle^{b}f\|_2. 
		\end{equation*}
		For $|x|>2$, we split $K_a*f(x)=I_1(x)+I_2(x)$, where
		\begin{equation*}
			I_1(x)=\int_{|y|\le|x|/2}
			(K_a(x-y)-K_a(x))f(y)\,dy,
		\end{equation*}
		and
		\begin{equation*}
			I_2(x)=\int_{|y|>|x|/2}
			(K_a(x-y)-K_a(x))f(y)\,dy. 
		\end{equation*}
		First, we shall estimate $I_1$.  For $|y|\le|x|/2$, the mean-value theorem gives
		\begin{equation*}
			K_a(x-y)-K_a(x)
			=-y\int_0^1K_a'(x-\theta y)\,d\theta.
		\end{equation*}
		Since $|x-\theta y|\ge|x|/2$, \eqref{ee2} yields
		\begin{equation*}
			|K_a(x-y)-K_a(x)|
			\le c_a|y||x|^{-a-3}. 
		\end{equation*}
		Thus
		\begin{equation*}
			|I_1(x)|\le c_a|x|^{-a-3}\int_{\mathbb{R}} |y||f(y)|\,dy. 
		\end{equation*}
		Since $b\ge a+3/2>3/2$,
		\begin{equation*}
			\int_{\mathbb{R}} |y||f(y)|\,dy
			\le c_b\| \langle x\rangle^{b}f\|_2.
		\end{equation*}
		Hence
		\begin{equation*}
			\langle x\rangle^{b}|I_1(x)|
			\le c_{a,b}|x|^{b-a-3}\|\langle x\rangle^{-b}f\|_2. 
		\end{equation*}
		This belongs to $L^2(|x|>2)$ if $b<a+\frac52$. Therefore
		\begin{equation*}
			\|\langle x\rangle^{b}I_1\|_2\le c_{a,b}\|\langle x\rangle^{b}f\|_2.
		\end{equation*}
		Finally, we estimate $I_2(x)$ writing it as 
		\begin{equation*}
			I_2(x)=I_{21}(x)-I_{22}(x),
		\end{equation*}
		where
		\begin{equation*}
			I_{21}(x)=\int_{|y|>|x|/2}K_a(x-y)f(y)\,dy, 
		\end{equation*}
		and
		\begin{equation*}
			I_{22}(x)=K_a(x)\int_{|y|>|x|/2}f(y)\,dy. 
		\end{equation*}
		For $I_{21}$, since  $|y|>|x|/2$,
		\begin{equation*}
			\langle x\rangle^{b}\le c_b \langle y\rangle^{b}
		\end{equation*} 
		and consequently 
		\begin{equation*}
			\|\langle x\rangle^{b}I_{21}\|_2\le c_{a,b}\|\langle x\rangle^{b}f\|_2.
		\end{equation*}
		For $I_{22}$, Cauchy-Schwarz gives
		\begin{equation*}
			|\int_{|y|>|x|/2}f(y)\,dy|
			\le
			(\int_{|y|>|x|/2}\langle y\rangle^{-2b}\,dy)^{1/2}
			\|\langle x\rangle^{b}f\|_2
			\le c_b|x|^{1/2-b}\|\langle x\rangle^{b}f\|_2.
		\end{equation*}
		Hence, combining \eqref{ee2} and the above estimate one gets
		\begin{equation*}
			\langle x\rangle^{b}|I_{22}(x)|
			\le c_{a,b}|x|^{-a-3/2}\|\langle x\rangle^{b}f\|_2,
		\end{equation*}
		which is square integrable on $|x|>2$. Hence
		\begin{equation*}
			\|\langle x\rangle^{b}I_{22}\|_2\le c_{a,b}\|\langle x\rangle^{b}f\|_2.
		\end{equation*}
		Gathering the above information we obtain 
		\begin{equation*}
			\|\langle x\rangle^{b}K_a*f\|_2\le c_{a,b}\|\langle x\rangle^{b}f\|_2,
		\end{equation*}
		and
		\begin{equation*}
			\|[\langle x\rangle^{b};T_0]f\|_2\le c_{a,b}\|\langle x\rangle^{b}f\|_2.
		\end{equation*}
		Combining the above  estimate with \eqref{ee2}  we obtain the desired result \eqref{ee1}
		for
		\begin{equation*}
			a+3/2\le b<a+5/2,\;\;\;\;\text{with}\;\;\;\;\int_{\mathbb{R}}f\, dx=0.
		\end{equation*}
	\end{proof}
	The following proposition allows us to recover non-homogeneous derivatives from homogeneous ones in the presence of polynomial weights. It will be used in several steps of the iterative argument in Section 3.
	\begin{prop}\label{JandDprop}
		Let $1<p<\infty$, $m\in \mathbb{Z}^{+}$, $\theta\in \mathbb{R}$, $0<\beta<1$. If
		\begin{equation*}
			\, \langle x\rangle^{m+\theta}v,\,  \langle x\rangle^{\theta}|D_x|^{\beta}(x^m v)\in L^p(\mathbb{R}),
		\end{equation*}
		then $\langle x\rangle^{\theta} x^m  J^{\beta} v\in L^p(\mathbb{R})$. Moreover,
		\begin{equation*}
			\|\langle x\rangle^{\theta} x^m J^{\beta}v\|_{p}\leq c_{_{m,\theta,\beta}} \big(\|\langle x\rangle^{m+\theta} v\|_{p}+\|\langle x\rangle^{\theta}|D_x|^{\beta}(x^m v)\|_{p}\big).
		\end{equation*} 
	\end{prop}
	\begin{proof}
		Let $\chi\in C_0^\infty(\mathbb{R})$ such that $0\leq \chi\leq 1$ and $\chi(x)=1$, whenever $|x|\leq 1$. We write
		\begin{equation*}
			\begin{aligned}
				\langle x\rangle^{\theta} x^m J^{\beta}v=&\langle x\rangle^{\theta}x^m \chi(D_x)J^{\beta}v+\langle x\rangle^{\theta}x^m (1-\chi(D_x))J^{\beta}v  \\
				=&[\langle x\rangle^{\theta}x^m;\chi(D_x)J^{\beta}]v +\chi(D_x)J^{\beta}\big(\langle x\rangle^{\theta}x^m v\big)+  \langle x\rangle^{\theta}[x^m; (1-\chi(D_x))J^{\beta}]v\\
				&+\langle x\rangle^{\theta}(1-\chi(D_x))J^{\beta}(x^m v).
			\end{aligned}
		\end{equation*}
		Since $\langle x \rangle^{\theta} x^m$ $\in \operatorname{Op}(\mathcal{S}^{0,m+\theta})$, $\chi(D_x) J^{\beta}\in \operatorname{Op}(\mathcal{S}^{0,0})$, by Lemma \ref{LemSGcalculus}, we have that $[\langle x\rangle^{\theta}x^m;\chi(D_x)J^{\beta}]\in \operatorname{Op}(\mathcal{S}^{0,m+\theta-1})$, and thus
		\begin{equation*}
			\begin{aligned}
				\| [\langle x\rangle^{\theta}x^m;\chi(D_x)J^{\beta}]v\|_{p} +\|\chi(D_x)J^{\beta}\big(\langle x\rangle^{\theta}x^m v\big)\|_{p}\leq c\|\langle x \rangle^{m+\theta} v\|_{p}. 
			\end{aligned}    
		\end{equation*}
		Using that $(1-\chi(D_x))J^{\beta}\in \operatorname{Op}(\mathcal{S}^{\beta,0})$, that $ \langle x\rangle^{\theta}[x^m; (1-\chi(D_x))J^{\beta}]\in \operatorname{Op}(\mathcal{S}^{\beta-1,m+\theta-1})\subseteq \operatorname{Op}(\mathcal{S}^{0,m+\theta})$, it follows
		\begin{equation*}
			\| \langle x\rangle^{\theta}[x^m; (1-\chi(D_x))J^{\beta}]v\|_{p}
			\leq c \|\langle x \rangle^{m+\theta} v\|_{p}.
		\end{equation*}
		Next, we write
		\begin{equation*}
			\begin{aligned}
				\langle x\rangle^{\theta}(1-\chi(D_x))J^{\beta}(x^m v)=&\langle x\rangle^{\theta}|D_x|^{-\beta}(1-\chi(D_x))J^{\beta}|D_x|^{\beta}(x^m v)\\
				=&[\langle x\rangle^{\theta};|D_x|^{-\beta}(1-\chi(D_x))J^{\beta}]|D_x|^{\beta}(x^m v)\\
				&+|D_x|^{-\beta}(1-\chi(D_x))J^{\beta}\big(\langle x\rangle^{\theta}|D_x|^{\beta}(x^m v)\big).
			\end{aligned}    
		\end{equation*}
		Using that $|D_x|^{-\beta}(1-\chi(D_x))\in \operatorname{Op}(\mathcal{S}^{-\beta,0})$, we have $[\langle x\rangle^{\theta};|D_x|^{-\beta}(1-\chi(D_x))J^{\beta}]\in \operatorname{Op}(\mathcal{S}^{-1,\theta-1})$. We deduce from the continuity of pseudo-differential operators that
		\begin{equation*}
			\begin{aligned}
				\|[\langle x\rangle^{\theta};|D_x|^{-\beta}(1-\chi(D_x))J^{\beta}]|D_x|^{\beta}(x^m v)\|_{p}+&\||D_x|^{-\beta}(1-\chi(D_x))J^{\beta}\big(\langle x\rangle^{\theta}|D_x|^{\beta}(x^m v)\big)\|_{p}\\
				\leq &  c\|\langle x\rangle^{\theta-1} |D_x|^{\beta}(x^m v)\|_{p}+\|\langle x\rangle^{\theta}|D_x|^{\beta}(x^m v)\|_{p}\\
				\leq & c \|\langle x\rangle^{\theta}|D_x|^{\beta}(x^m v)\|_{p}.
			\end{aligned}    
		\end{equation*}
		The proof is complete.
	\end{proof}
	\subsection{Further commutator Expansions }\label{sectionginibrevelo}
	
	In this part, we first present several tools obtained by Ginibre and Velo \cite{GV1}, \cite{GV2}. The basic problem is to handle the non-local operator $|D_{x}|^{s}$ for non-integer $s$ and, in particular, to obtain representations of its commutator with multiplication operators by functions that exhibit as much locality as possible.
	
	Let $b=2\mu+1>1,$ let $n$ be a non-negative integer, and let $h$ be a smooth function such that its derivatives have suitable decay at infinity.

	We define the operator
	\begin{equation}\label{eq8}
		{{R}}_{n}(b)=-\left[\mathcal{H}\, |D_{x}|^{b}; h\right]-\frac{1}{2}\left(\mathcal{P}_{n}(b)-\mathcal{H}\,\mathcal{P}_{n}(b)\,\mathcal{H}\right).
	\end{equation}
	where $\mathcal H$ stands for the Hilbert transform, see \eqref{hita}, and 
	
	\begin{equation}\label{eq105}
		\mathcal{P}_{n}(b)=b\sum_{0\leq j\leq n} c_{2j+1}(-1)^{j}4^{-j} |D_{x}|^{\mu-j}\left(h^{(2j+1)}|D_{x}|^{\mu-j}\right).
	\end{equation}
	where
	\begin{equation*}
		c_{1}=1,\hskip10pt\text{and}\hskip10pt c_{2j+1}=\frac{1}{(2j+1)!}\prod_{0\leq k<j}\left(b^{2}-(2k+1)^{2}\right).
	\end{equation*}
	It was shown in \cite{GV1} that the operator ${{R}}_{n}(b)$ can be represented in terms of anti-commutators\footnote{For any two operators $P$ and $Q$ we denote the \emph{anti-commutator} by $[P;Q]_{+}:=PQ+QP.$} as follows
	\begin{equation}\label{aj1}
		{R}_{n}(b)=-\frac{1}{2}\left(\left[\mathcal{H}; \mathcal{Q}_{n}(b)\right]_{+}
		+
		\left[|D_{x}|^{b};\left[\mathcal{H}; h\right]\right]_{+}\right).
	\end{equation}
	where the operator $	\mathcal{Q}_{n}(a)$ is represented in the Fourier space variables by the integral kernel
	\begin{equation}\label{aj2}
		\mathcal{Q}_{n}(b)\longrightarrow (2\pi)^{\frac{1}{2}}
		\widehat{h}(\xi-\xi')|\xi\xi'|^{\frac{b}{2}}2bq_{n}(b,t),
	\end{equation}
	with $|\xi|=|\xi'|e ^{2t}$ and
	\begin{equation*}
		q_{n}(b,t)=\frac{1}{b}(b^{2}-(2n+1)^{2})c_{2n+1}\int_{0}^{t} \sinh^{2n+1}\tau\,\sinh((b(t-\tau)))\,d\tau.
	\end{equation*}
	Based on \eqref{aj1} and \eqref{aj2}, Ginibre and Velo \cite{GV2}
	obtain the following boundedness and compactness properties of the operator ${{R}}_{n}(b).$
	\begin{prop}\label{propo2}
		Let $n$ be a non-negative integer, $b\geq 1,\,$ and $\sigma\geq 0,$ be such that
		\begin{equation}\label{eq21}
			2n+1\leq b+2\sigma\leq2n+3.
		\end{equation}
		Then
		\begin{itemize}
			\item[(a)] The operator $|D_x|^{\sigma}	{{R}}_{n}(b)|D_x|^{\sigma}$ is bounded in $L^{2}$ with norm
			\begin{equation}\label{eq98}
				\left\||D_{x}|^{\sigma}	{{R}}_{n}(b)|D_{x}|^{\sigma}f\right\|_{2}\leq c(2\pi)^{-1/2}\left\|\widehat{(|D_{x}|^{b+2\sigma}h)}\right\|_{L^{1}_{\xi}}\|f\|_{2}.
			\end{equation}
			If $b\geq 2n+1,$ one can take $c=1.$
			\item[(b)] Assume in addition that
			\begin{equation*}
				2n+1\leq b+2\sigma<2n+3.
			\end{equation*}
			Then the operator ${\displaystyle |D_{x}|^{\sigma} {{R}}_{n}(b)|D_{x}|^{\sigma}}$ is compact in $L^{2}(\mathbb{R}).$
		\end{itemize}
	\end{prop}
	\begin{proof}
		See Proposition 2.2 in \cite{GV2}.
	\end{proof}
	In fact, Proposition \ref{propo2} is a generalization of a previous result in which the derivatives of the operator $\tilde{R}_{n}(b)$ are not considered (cf. Proposition 1 in \cite{GV1}). 
	
	The estimate (\ref{eq98}) yields the following identity for the localization of derivatives.
	\begin{lem}\label{prop1}
		Assume $0<\kappa<1.$ Let $\varphi\in C^{\infty}(\mathbb{R})$ with $\varphi'\in C_0^{\infty}(\mathbb{R})$. Then,
		\begin{equation}\label{eq1ab}
			\begin{split}
				\int_{\mathbb{R}} \varphi f\, |D_{x}|^{\kappa+1}\partial_{x}f\,dx&= -\left(\frac{\kappa+2}{4}\right)\int_{\mathbb{R}} \left(\left||D_{x}|^{\frac{\kappa+1}{2}}f\right|^{2}+\left||D_{x}|^{\frac{\kappa+1}{2}}\mathcal{H}f\right|^{2}\right)\varphi' \,dx\\
				&\quad-\frac{1}{2}\int_{\mathbb{R}}f\,{{R}}_{0}(\kappa+2)f\,dx.
			\end{split}
		\end{equation}
	\end{lem}
	\begin{proof}
		First notice that
		\begin{equation}\label{eq2}
			\left[|D_{x}|^{\kappa+1}\partial_{x}; \varphi\right]f=|D_{x}|^{\kappa+1}\partial_{x}(\varphi f)-\varphi|D_{x}|^{\kappa+1}\partial_{x}f.
		\end{equation}
		Using identity (\ref{eq2}) and integration by parts, we obtain
		\begin{equation}\label{eq3}
			\begin{split}
				\int_{\mathbb{R}}\varphi\,f\,|D_{x}|^{\kappa+1}(\partial_{x}f)\,dx&=-\int_{\mathbb{R}}f\left[|D_{x}|^{\kappa+1}\partial_{x}; \varphi\right]f\,dx +\int_{\mathbb{R}}f\,|D_{x}|^{\kappa+1}\partial_{x} \left(\varphi\,f\right)\,dx\\
				&=-\int_{\mathbb{R}}f\left[|D_{x}|^{\kappa+1}\partial_{x}; \varphi\right]f\,dx-\int_{\mathbb{R}}\partial_{x}f\,|D_{x}|^{\kappa+1} \left(\varphi\,f\right)\,dx.
			\end{split}
		\end{equation}
		By the symmetry of the operator $|D_x|^{\kappa+1}$, we find  
		\begin{equation}\label{eq4}
			\int_{\mathbb{R}}\partial_{x}f\,|D_{x}|^{\kappa+1} \left(\varphi\,f\right)\,dx=\int_{\mathbb{R}} \varphi\,f\,|D_{x}|^{\kappa+1}\partial_{x}f\, \,dx.
		\end{equation}
		Thus, gathering (\ref{eq3}) and (\ref{eq4}) we obtain 
		\begin{equation}\label{eq92}
			\begin{split}
				\int_{\mathbb{R}}\varphi\,f\,|D_{x}|^{\kappa+1}(\partial_{x}f)\,dx &=-\frac{1}{2}\int_{\mathbb{R}}f\left[|D_{x}|^{\kappa+1}\partial_{x}; \varphi\right]f\,dx\\
				&=\frac{1}{2} \int_{\mathbb{R}}f\left[\mathcal{H}|D_{x}|^{\kappa+2}; \varphi\right]f\,dx.
			\end{split}
		\end{equation}
		Now, we apply \eqref{eq8}, with ${\displaystyle \mu=\frac{\kappa+1}{2}}$ and $n=0$ to obtain
		\begin{equation}\label{eq93}
			\begin{split}
				[\mathcal{H}|D_{x}|^{\kappa+2};\varphi] &=-\frac{1}{2}\mathcal{P}_{0}(\kappa+2)-\frac{1}{2}\mathcal{H}\mathcal{P}_{0}(\kappa+2)\mathcal{H}+\tilde{R}_{0}(\kappa+2)\\
				&=-\left(\frac{\kappa+2}{2}\right)|D_{x}|^{\frac{\kappa+1}{2}}\varphi' |D_{x}|^{\frac{\kappa+1}{2}}-\left(\frac{\kappa+2}{2}\right)\mathcal{H}|D_{x}|^{\frac{\kappa+1}{2}}\varphi' |D_{x}|^{\frac{\kappa+1}{2}}\mathcal{H}-{{R}}_{0}(\kappa+2).
			\end{split}
		\end{equation}
		Hence, inserting (\ref{eq93}) into (\ref{eq92}) and using Plancherel's identity, we obtain 
		\begin{equation}\label{eq94}
			\begin{split}
				\int_{\mathbb{R}} f[\mathcal{H}|D_{x}|^{\kappa+2};\varphi]f\,dx&= -\Big(\frac{\kappa+2}{2}\Big)\int_{\mathbb{R}}||D_{x}|^{\frac{\kappa+1}{2}}f|^{2}\varphi'(x)\,dx\\
				&\hskip15pt-\Big(\frac{\kappa+2}{2}\Big)\int_{\mathbb{R}}||D_{x}|^{\frac{\kappa+1}{2}}\mathcal{H}f|^{2}\varphi'(x)\,dx \\
				&\hskip15pt-\int_{\mathbb{R}}
				f{{R}}_{0}(\kappa+2)f\,dx .
			\end{split}
		\end{equation}
		The identity \eqref{eq1ab} follows from \eqref{eq92} and \eqref{eq94}.
	\end{proof}
	\begin{lem}\label{prop1a}
		Assume $0<\kappa<1.$ Let $\varphi\in C^{\infty}(\mathbb{R})$ with $\varphi'\in C_0^{\infty}(\mathbb{R})$. Then,
		\begin{equation}\label{eq1a}
			\begin{aligned}
				\int_{\mathbb{R}} \varphi f\, |D_{x}|^{\kappa+1}\partial_{x}f\,dx&= -\left(\frac{\kappa+2}{2}\right)\int_{\mathbb{R}} (|D_{x}|^{\frac{\kappa+1}{2}}f)^{2}\,\varphi' \,dx\\
				&\hskip12pt-\frac{1}{2}\int_{\mathbb{R}} f\,{{R}}_{0}(\kappa+2)f\,dx+\int_{\mathbb{R}} f E_0(\kappa+2)f dx
			\end{aligned}
		\end{equation}
		with
		\begin{equation*}
			\Big|\int_{\mathbb{R}} f \,E_0(\kappa+2)f \,dx\Big|     \leq c\| \widehat{(|D_x|^{\kappa+1}\varphi')}\|_{L^1_{\xi}}\,\| f\|_2^2.
		\end{equation*}
		Moreover, with the parameters as in \eqref{eq98}, i.e., $\kappa+1+2\sigma\leq 3$
		\begin{equation}\label{eq98b}
			\left\||D_{x}|^{\sigma}	E_{0}(\kappa+2)|D_{x}|^{\sigma}f\right\|_{2}\leq c \left\|\widehat{(|D_{x}|^{{\kappa}+1+2\sigma}{\varphi})}\right\|_{L^{1}_{\xi}}\|f\|_{2}.
		\end{equation}
	\end{lem}
	\begin{proof}
		We write
		\begin{equation}\label{abc1}
			\begin{split}
				\int_{\mathbb{R}} |D_x|^{\frac{\kappa+1}{2}}\mathcal Hf\,|D_x|^{\frac{\kappa+1}{2}}\mathcal Hf\,\varphi'\,dx&=-\int_{\mathbb{R}} |D_x|^{\frac{\kappa+1}{2}}f\,\mathcal H(|D_x|^{\frac{\kappa+1}{2}}\mathcal Hf\,\varphi')\,dx\\
				&=\int_{\mathbb{R}} (|D_x|^{\frac{\kappa+1}{2}}f)^2\,\varphi'\,dx-\int_{\mathbb{R}} |D_x|^{\frac{\kappa+1}{2}}f\,[\mathcal H;\varphi'](|D_x|^{\frac{\kappa+1}{2}}\mathcal Hf)\,dx.
			\end{split}
		\end{equation}
		Thus, defining
		\begin{equation*}
			\int_{\mathbb{R}} fE_0(\kappa+{2})f dx=-\int_{\mathbb{R}} |D_x|^{\frac{\kappa+1}{2}}f\,[\mathcal H;\varphi'](|D_x|^{\frac{\kappa+1}{2}}\mathcal Hf)\,dx
			=-\int_{\mathbb{R}} f\,|D_x|^{\frac{\kappa+1}{2}}[\mathcal H;\varphi'] (|D_x|^{\frac{\kappa+1}{2}}\mathcal Hf)\,dx
		\end{equation*}
		we see by \eqref{calderon2} that
		\begin{equation*}
			\Big|\int_{\mathbb{R}} f\,E_0(\kappa+2)f\,dx\Big|\leq c\,\|\widehat { (|D_x|^{\kappa+{2}}\varphi')}\|_{L^1_{\xi}}\,\| f\|_2^2,
		\end{equation*}
		and  \eqref{eq98} also hold. Thus, combining these estimates and \eqref{eq1ab}, we obtain the desired result. 
	\end{proof}
	\begin{rem}
		Since the operators, ${{R}}_{0}(\kappa+2)$ and $E_0(\kappa+2)$ satisfy the same estimate. we define
		\begin{equation*}
			\tilde{R}_{0}(\kappa+2):=-\frac{1}{2}\tilde{R}_{0}(\kappa+2)+E_0(\kappa+2),
		\end{equation*}
		and use identity \eqref{eq1a} as
		\begin{equation}\label{eq1abc}
			\int_{\mathbb{R}} \varphi f\, |D_{x}|^{\kappa+1}\partial_{x}f\,dx= -\left(\frac{\kappa+2}{2}\right)\int_{\mathbb{R}} (|D_{x}|^{\frac{\kappa+1}{2}}f)^{2}\,\varphi' \,dx
			+\int_{\mathbb{R}} f\,\tilde{R}_{0}(\kappa+2)f\,dx,
		\end{equation}
		with $\tilde{R}_{0}(\kappa+2)$ satisfying \eqref{eq98}.
	\end{rem}
	\begin{rem} It should be clear that Lemma \ref{prop1a} extends to $\varphi \in C^{\infty}(\mathbb{R})$ with 
		$\widehat{(|D_{x}|^{\kappa+1+2\sigma}\varphi})\in L^1(\mathbb{R})$.
	\end{rem}
	In our applications we take $\kappa=a-1\in(0,1)$, so that $|D_x|^{\kappa+1}\partial_x=\partial_x|D_x|^{a}$ and $\tilde{R}_0(\kappa+2)=\tilde{R}_0(a+1)$. Identity \eqref{eq1abc} is the key tool to extract the local smoothing effect from the dispersive term in the weighted energy estimates. To handle the commutators between fractional derivatives and the truncated weights, we also need the following expansion, due to Ginibre and Velo \cite{GV2}.
	Let $b>0$, let $\delta\in\{0,1\}$, and let $n\in\mathbb{N}_{0}$. 
	Let $h$ be a smooth function with sufficient decay at infinity; for instance, assume that 
	$\partial_{x}h\in\mathcal{S}(\mathbb{R})$.
	
	We define the operator
	\begin{equation*}
		\mathcal{S}_{\delta,n}(b)
		:=
		\left[(-\mathcal{H})^{\delta}|D_{x}|^{b};\, h\right]
		+
		\sum_{j=1}^{n}
		\binom{b}{j}
		(-\mathcal{H})^{j+\delta}
		|D_{x}|^{b-j}
		\big(\,\cdot\, \partial_{x}^{j} h\big),
	\end{equation*}
	where the summation term is omitted if $n=0$. Here
	\begin{equation*}
		\binom{b}{j}
		:=
		\frac{b(b-1)(b-2)\cdots(b-j+1)}{j!},
	\end{equation*}
	and $\partial_{x}^{j}h$ denotes the $j$-th derivative of $h$, viewed as a multiplication operator.
	\begin{thm}\label{gv2}
		Let $b>0$, and let $\alpha,\beta\ge 0$ satisfy
		\begin{equation}\label{condition}
			n \le b+\alpha \le b+\alpha+\beta \le n+1.
		\end{equation}
		Then the operator
		\begin{equation*}
			|D_{x}|^{\alpha}\,\mathcal{S}_{\delta,n}(b)\,|D_{x}|^{\beta}
		\end{equation*}
		extends to a bounded operator on $L^{2}(\mathbb{R})$. More precisely, for all 
		$f\in L^{2}(\mathbb{R})$,
		\begin{equation}\label{condicons}
			\left\|
			|D_{x}|^{\alpha}\,\mathcal{S}_{\delta,n}(b)\,|D_{x}|^{\beta} f
			\right\|_{2}
			\le
			\frac{c}{2\sqrt{\pi}}
			\left\|
			\widehat{|D_{x}|^{b+\alpha+\beta} h}
			\right\|_{L^{1}_{\xi}}
			\,
			\|f\|_{2},
		\end{equation}
		where
		\begin{equation}\label{condicons1}
			c	=	1	+\sum_{j=0}^{n}\left|
			\binom{b}{j}
			\right|.
		\end{equation}
		Assume in addition that
		\begin{equation*}
			b+\alpha+\beta < n+1.
		\end{equation*}
		Then the operator
		\begin{equation*}
			|D_{x}|^{\alpha}\,\mathcal{S}_{\delta,n}(b)\,|D_{x}|^{\beta}
		\end{equation*}
		is compact on $L^{2}(\mathbb{R})$.
	\end{thm}
	\begin{proof}
		See the work of Ginibre and Velo \cite{GV2}.
	\end{proof}
	\subsection{Nonlinear estimates}
	
	As an additional consequence of the pseudo-differential calculus, we shall derive weighted nonlinear estimates. This result exploits the pseudo-differential nature of the operators $\langle x\rangle$ and $J$. Before doing so, let us first recall the following bounded approximation of $\langle x \rangle$. Given $\theta>0$, $N\in \mathbb{N}$, we define the truncated weights $\langle \cdot \rangle_{N}^{\theta} : \mathbb{R} \rightarrow \mathbb{R}$ as
	\begin{equation}\label{WNwightApp}
		\langle x \rangle_{N}^{\theta} =\left\{\begin{aligned} 
			& (1+ |x|^2)^{\frac{\theta}{2}}, \hspace{0.3cm} \text{ if } |x|\leq N, \\
			&\Big(\frac{3N}{2}\Big)^{\theta}, \hspace{1cm}\text{if } |x|\geq 2N,
		\end{aligned}\right.
	\end{equation}
	in such a way that $\langle x \rangle_{N}$ is smooth, even and non-decreasing in $|x|$ with $|\partial_x\langle x \rangle_{N}| \leq 1$,  and there exist constants $c_{l}>0$ independent of $N$ such that $|\partial_x^l\langle x \rangle_{N}| \leq c_l|\partial_x^l\langle x \rangle|$, for each $l=2,3,4,5$. 
	\begin{lem}\label{Nonlinear estimates}
		Let $\sigma>\frac{1}{2}$, $0<\theta\leq \frac{1}{2}$, $m,l\in \mathbb{Z}^{+}\cup\{0\}$. Then 
		\begin{equation}\label{nonlinearE1}
			\Big|\int_{\mathbb{R}} \big(v\partial_x v\big) \big(x^{2m+1}\langle x\rangle^{2\theta}_N \langle x\rangle^{-l} v\big)\, dx \Big|\leq c_{m,l,\theta} 
			\left\{\begin{aligned}
				&\|v\|_{3}^3, \hspace{2.8cm} \text{ if } m+\theta-\frac{l}{2} \leq 0, \\
				&\|v\|_{\infty}\|\langle x \rangle^{m+\theta-\frac{l}{2}} v\|_{2}^2, \hspace{0.5cm} \text{ if } m+\theta-\frac{l}{2} > 0.
			\end{aligned}\right.  
		\end{equation}
		Moreover,
		\begin{equation}\label{nonlinearE2} 
			\begin{aligned}
				\Big|\int_{\mathbb{R}} & J^{\sigma}\big(v\partial_x v\big) \big(x^{2m+1}\langle x\rangle^{2\theta}_N \langle x\rangle^{-l} J^{\sigma}v\big)\, dx\Big|\\
				&\leq c_{m,l,\theta}
				\left\{\begin{aligned}
					&\big(\|J^{\sigma}v\|_{2}+\|\partial_x v\|_{\infty}\big)\|\langle x\rangle^{m+\frac{1}{2}+\theta-\frac{l}{2}} J^{\sigma} v\|_{2}^2, \hspace{1cm} \text{ if } \frac{1}{2}< \sigma \leq 1, \\
					&\big(\|J^{3/2}v\|_{2}+\|\partial_x v\|_{\infty}\big)\|\langle x\rangle^{m+\frac{1}{2}+\theta-\frac{l}{2}} J^{\sigma} v\|_{2}^2, \hspace{0.7cm} \text{ if } 1 <\sigma < 3/2, \\
					&\|J^{3/2^{+}}v\|_{2}\|\langle x\rangle^{m+\frac{1}{2}+\theta-\frac{l}{2}} J^{\sigma} v\|_{2}^2, \hspace{2.45cm} \text{ if } \sigma= 3/2, \\
					&\|J^{\sigma}v\|_{2}\|\langle x\rangle^{m+\frac{1}{2}+\theta-\frac{l}{2}} J^{\sigma} v\|_{2}^2, \hspace{2.9cm} \text{ if } \sigma>3/2,
				\end{aligned}\right. 
			\end{aligned}    
		\end{equation}
		where all the implicit constants above are independent of $N$. Same estimates hold for \eqref{nonlinearE1} and \eqref{nonlinearE2} with $\langle x\rangle_N$ replaced by $\langle x\rangle$.
	\end{lem}
	\begin{proof}
		Using that $v^2\partial_x v=\frac{1}{3}\partial_x (v^3)$ and that $\Big|\Big(x^{2m+1}\langle x\rangle_N^{2\theta}\langle x \rangle^{-l}\Big)'\Big|\leq c_{m,l,\theta}\langle x \rangle^{2m+2\theta-l}$, where $c_{m,l,\theta}>0$ independent of $N$, we have that \eqref{nonlinearE1} follows from integration by parts.
		
		We turn to the deduction of \eqref{nonlinearE2}. We write 
		\begin{equation*}
			\begin{aligned}
				&\int_{\mathbb{R}} J^{\sigma}\big(v\partial_x v\big) \big(x^{2m+1}\langle x\rangle^{2\theta}_N \langle x\rangle^{-l} J^{\sigma}v\big)\, dx\\
				\quad &= \int_{\mathbb{R}} \langle x\rangle^{m+\frac{1}{2}+\theta-\frac{l}{2}}J^{\sigma}\big(v\partial_x v\big) \big( \langle x\rangle^{m+\frac{1}{2}+\theta-\frac{l}{2}}J^{\sigma}v\big)\, \frac{x^{2m+1}\langle x\rangle^{2\theta}_N}{\langle x\rangle^{2m+1+2\theta}}\, dx  \\
				\quad &=:\mathcal{N}_1,
			\end{aligned}    
		\end{equation*}
		where the fact that $\langle x\rangle_{N}\leq \langle x\rangle$ shows $\Big|\frac{x^{2m+1}\langle x\rangle^{2\theta}_N}{\langle x\rangle^{2m+1+2\theta}}\Big|\leq 1$. Note also that this will justify why the desired result is valid when replacing $\langle x \rangle_N$ with $\langle x\rangle$.
		
		Using commutators, we divide our arguments as follows
		\begin{equation*}
			\begin{aligned}
				\mathcal{N}_1=&\int_{\mathbb{R}} [\langle x\rangle^{m+\frac{1}{2}+\theta-\frac{l}{2}};J^{\sigma}](v\partial_x v)\big(\langle x\rangle^{m+\frac{1}{2}+\theta-\frac{l}{2}}J^{\sigma} v\big)\frac{x^{2m+1}\langle x\rangle^{2\theta}_N}{\langle x\rangle^{2m+1+2\theta}}\, dx\\
				&+\int_{\mathbb{R}} [J^{\sigma};\langle x\rangle^{m+\frac{1}{2}+\theta-\frac{l}{2}} v]\partial_x v \big(\langle x\rangle^{m+\frac{1}{2}+\theta-\frac{l}{2}}J^{\sigma}v\big)\frac{x^{2m+1}\langle x\rangle^{2\theta}_N}{\langle x\rangle^{2m+1+2\theta}}\, dx \\
				&+\int_{\mathbb{R}} \langle x\rangle^{m+\frac{1}{2}+\theta-\frac{l}{2}}v \partial_x J^{\sigma}v  \big(\langle x\rangle^{m+\frac{1}{2}+\theta-\frac{l}{2}} J^{\sigma}v\big)\frac{x^{2m+1}\langle x\rangle^{2\theta}_N}{\langle x\rangle^{2m+1+2\theta}}\, dx\\
				=:& \mathcal{N}_{1,1}+\mathcal{N}_{1,2}+\mathcal{N}_{1,3}.
			\end{aligned}    
		\end{equation*}
		By means of the asymptotic expansion (v) in Lemma \ref{LemSGcalculus}, it is not difficult to deduce that given $M\in \mathbb{Z}^{+}$ with $M\geq \max\{m+\frac{1}{2}+\theta-\frac{l}{2};\sigma \}$, there exist $a_1,\dots,a_{M},a_{M+1}\in \mathcal{S}^{0,0}$ such that
		\begin{equation}\label{Pseudodecomp}
			\begin{aligned}
				\,[\langle x\rangle^{m+\frac{1}{2}+\theta-\frac{l}{2}};J^{\sigma}]\partial_x (v^2)=&\sum_{j=1}^M T_{a_j}\big(J^{\sigma+1-j}\big(\langle x\rangle^{m+\frac{1}{2}+\theta-\frac{l}{2}-j} v^2\big)\big)+T_{a_{M+1}}(v^2).    
			\end{aligned}    
		\end{equation}
		To estimate the terms in the sum above, in the cases where $\sigma+1-j>0$, we apply Lemma \ref{lem1} to get
		\begin{equation}\label{Pseudodecomp1}
			\begin{aligned}
				\|J^{\sigma+1-j}\big(\langle x\rangle^{m+\frac{1}{2}+\theta-\frac{l}{2}-j} v^2\big)\|_{2}\leq c &\big(\|v\|_{\infty} \|J^{\sigma+1-j}\big(\langle x\rangle^{m+\frac{1}{2}+\theta-\frac{l}{2}-j} v\big)\|_{2}\\
				&+\|J^{\sigma+1-j}v\|_{2} \|\langle x\rangle^{m+\frac{1}{2}+\theta-\frac{l}{2}-j} v\|_{\infty}\big).
			\end{aligned}   
		\end{equation}
		Given that $\sigma>\frac{1}{2}$, using that $\langle x\rangle^{m+\frac{1}{2}+\theta-\frac{l}{2}-j}\leq \langle x\rangle^{m+\frac{1}{2}+\theta-\frac{l}{2}}$, Sobolev embedding and that $[J^{\sigma};\langle x\rangle^{m+\frac{1}{2}+\theta-\frac{l}{2}}]\in \operatorname{Op}(\mathcal{S}^{\sigma-1,m+\frac{1}{2}+\theta-\frac{l}{2}-1})\subseteq \operatorname{Op}(\mathcal{S}^{\sigma,m+\frac{1}{2}+\theta-\frac{l}{2}})$, we get from (iii) in Lemma \ref{LemSGcalculus} that
		\begin{equation*}
			\begin{aligned}
				\|\langle x\rangle^{m+\frac{1}{2}+\theta-\frac{l}{2}-j} v\|_{\infty} \leq &  c  \|J^{\sigma}\big(\langle x\rangle^{m+\frac{1}{2}+\theta-\frac{l}{2}} v\big)\|_{2} \\
				\leq &  c\big( \|[J^{\sigma};\langle x\rangle^{m+\frac{1}{2}+\theta-\frac{l}{2}}]v\|_{2}+\|\langle x\rangle^{m+\frac{1}{2}+\theta-\frac{l}{2}}J^{\sigma}v\|_{2}\big)\\
				\leq &  c \|\langle x\rangle^{m+\frac{1}{2}+\theta-\frac{l}{2}}J^{\sigma}v\|_{2}.
			\end{aligned}    
		\end{equation*}
		To complete the estimate of \eqref{Pseudodecomp1}, we observe
		\begin{equation*}
			\begin{aligned}
				J^{\sigma+1-j}\big(\langle x\rangle^{m+\frac{1}{2}+\theta-\frac{l}{2}-j} v\big)=&J^{1-j}[J^{\sigma};\langle x\rangle^{m+\frac{1}{2}+\theta-\frac{l}{2}-j}]v\\
				&+J^{1-j}\big(\langle x\rangle^{-j} \langle x\rangle^{m+\frac{1}{2}+\theta-\frac{l}{2}}J^{\sigma}v\big).  
			\end{aligned}
		\end{equation*}
		Hence, since $[J^{\sigma};\langle x\rangle^{m+\frac{1}{2}+\theta-\frac{l}{2}-j}]\in \operatorname{Op}(\mathcal{S}^{\sigma-1,m+\frac{1}{2}+\theta-\frac{l}{2}-j-1})$, we deduce from the previous decomposition and the continuity of the pseudo-differential operators (iii) in Lemma \ref{LemSGcalculus} that
		\begin{equation*}
			\begin{aligned}
				\|J^{\sigma+1-j}\big(\langle x\rangle^{m+\frac{1}{2}+\theta-\frac{l}{2}-j} v\big)\|_{2}\leq c   \|\langle x\rangle^{m+\frac{1}{2}+\theta-\frac{l}{2}}J^{\sigma}v\|_{2}. 
			\end{aligned}    
		\end{equation*}
		Then, we get from \eqref{Pseudodecomp} that
		\begin{equation*}
			\begin{aligned}
				\|[\langle x\rangle^{m+\frac{1}{2}+\theta-\frac{l}{2}},J^{\sigma}]\partial_x (v^2)\|_{2}\leq c  \|J^{\sigma} v\|_{2}\|\langle x\rangle^{m+\frac{1}{2}+\theta-\frac{l}{2}}J^{\sigma}v\|_{2}.
			\end{aligned}    
		\end{equation*}
		The cases $\sigma+1-j\leq 0$ are easier to control, and similarly establish the same bound. As a consequence, we arrive at
		\begin{equation*}
			\begin{aligned}
				|\mathcal{N}_{11}|\leq c \|J^{\sigma} v\|_{2}\|\langle x\rangle^{m+\frac{1}{2}+\theta-\frac{l}{2}}J^{\sigma}v\|_{2}^2.
			\end{aligned}    
		\end{equation*}
		Next, if $\sigma>3/2$, we apply \eqref{kp-comm}, Sobolev embedding, and arguments analogous to those used before to get
		\begin{equation*}
			\begin{aligned}
				|\mathcal{N}_{12}|\leq c &\big(\|\partial_x\big(\langle x\rangle^{m+\frac{1}{2}+\theta-\frac{l}{2}} v\big)\|_{\infty}\|J^{\sigma}v\|_{2}\\
				&\quad +\|J^{\sigma}\big(\langle x\rangle^{m+\frac{1}{2}+\theta-\frac{l}{2}} v\big)\|_{2}\|\partial_x v\|_{\infty}\big)\|\langle x\rangle^{m+\frac{1}{2}+\theta-\frac{l}{2}}J^{\sigma}v\|_{2}\\
				\leq c &\|J^{\sigma}\big(\langle x\rangle^{m+\frac{1}{2}+\theta-\frac{l}{2}} v\big)\|_{2}\|J^{\sigma}v\|_{2}\|\langle x\rangle^{m+\frac{1}{2}+\theta-\frac{l}{2}}J^{\sigma}v\|_{2}\\
				\leq c & \|J^{\sigma}v\|_{2}\|\langle x\rangle^{m+\frac{1}{2}+\theta-\frac{l}{2}}J^{\sigma}v\|_{2}^2. 
			\end{aligned}   
		\end{equation*}
		If $1< \sigma \leq 3/2$, a further application of \eqref{kp-comm} and previous ideas establish
		\begin{equation*}
			\begin{aligned}
				|\mathcal{N}_{12}|\leq c &\big(\|\partial_x\big(\langle x\rangle^{m+\frac{1}{2}+\theta-\frac{l}{2}} v\big)\|_{p_1}\|J^{\sigma}v\|_{p_2}\\
				&\quad +\|J^{\sigma}\big(\langle x\rangle^{m+\frac{1}{2}+\theta-\frac{l}{2}} v\big)\|_{2}\|\partial_x v\|_{\infty}\big)\|\langle x\rangle^{m+\frac{1}{2}+\theta-\frac{l}{2}}J^{\sigma}v\|_{2}\\
				\leq c &\|\partial_x\big(\langle x\rangle^{m+\frac{1}{2}+\theta-\frac{l}{2}} v\big)\|_{p_1}\|J^{\sigma}v\|_{p_2}\|\langle x\rangle^{m+\frac{1}{2}+\theta-\frac{l}{2}}J^{\sigma}v\|_{2}\\
				&\quad +c\|\partial_x v\|_{\infty}\|\langle x\rangle^{m+\frac{1}{2}+\theta-\frac{l}{2}}J^{\sigma}v\|_{2}^2,
			\end{aligned}   
		\end{equation*}
		where $\frac{1}{p_1}+\frac{1}{p_2}=\frac{1}{2}$. According to $0<\sigma-1\leq \frac{1}{2}$, let us choose $p_1,p_2$. If $0<\sigma-1<\frac{1}{2}$, we set $\beta_1=\sigma-1$, $\beta_2=3/2-\sigma$, and then $\frac{1}{p_1}=\frac{1}{2}-\beta_1$, $\frac{1}{p_2}=\frac{1}{2}-\beta_2$. Applying Hardy–Littlewood–Sobolev inequality, we deduce
		\begin{equation*}
			\begin{aligned}
				\|\partial_x\big(\langle x\rangle^{m+\frac{1}{2}+\theta-\frac{l}{2}} v\big)\|_{p_1}\|J^{\sigma}v\|_{p_2}\leq  c & \|J^{\sigma}\big(\langle x\rangle^{m+\frac{1}{2}+\theta-\frac{l}{2}} v\big)\|_{2}\|J^{3/2}v\|_{2}\\ 
				\leq c &\|J^{3/2}v\|_{2}\|\langle x\rangle^{m+\frac{1}{2}+\theta-\frac{l}{2}} J^{\sigma} v \|_{2}.
			\end{aligned}    
		\end{equation*}
		Now, if $\sigma-1=\frac{1}{2}$, i.e., $\sigma=3/2$, for a given $0<\epsilon\ll 1$, we set $\beta_1=\frac{1}{2}-\epsilon$, $\beta_2=\epsilon$, and $\frac{1}{p_1}=\frac{1}{2}-\beta_1$, $\frac{1}{p_2}=\frac{1}{2}-\beta_2$. Once again by Hardy–Littlewood–Sobolev inequality, it is seen that
		\begin{equation*}
			\begin{aligned}
				\|\partial_x\big(\langle x\rangle^{m+\frac{1}{2}+\theta-\frac{l}{2}} v\big)\|_{p_1}\|J^{\sigma}v\|_{p_2}\leq c & \|J^{3/2-\epsilon}\big(\langle x\rangle^{m+\frac{1}{2}+\theta-\frac{l}{2}} v\big)\|_{2}\|J^{3/2+\epsilon}v\|_{2}\\
				\leq c & \|\langle x\rangle^{m+\frac{1}{2}+\theta-\frac{l}{2}} J^{\sigma}v\|_{2}\|J^{3/2+\epsilon}v\|_{2}.
			\end{aligned}    
		\end{equation*}
		Summarizing, if $1< \sigma\leq 3/2$, we have
		\begin{equation*}
			\begin{aligned}
				|\mathcal{N}_{12}|&\leq c
				\left\{\begin{aligned}
					&\big(\|J^{3/2} v\|_{2}+\|\partial_x v\|_{\infty}\big)\|\langle x\rangle^{m+\frac{1}{2}+\theta-\frac{l}{2}} J^{\sigma} v\|_{2}^2, \hspace{0.5cm} \text{ if } 1<\sigma< 3/2, \\
					&\|J^{3/2^{+}} v\|_{2}\|\langle x\rangle^{m+\frac{1}{2}+\theta-\frac{l}{2}} J^{\sigma} v\|_{2}^2, \hspace{2.4cm} \text{ if } \sigma=3/2.
				\end{aligned}\right. 
			\end{aligned}    
		\end{equation*}
		Next, we consider the case $\frac{1}{2}<\sigma\leq 1$. Using \eqref{kp-comm2} and preceding arguments, we obtain
		\begin{equation*}
			\begin{aligned}
				|\mathcal{N}_{12}|\leq c & \|\partial_x v\|_{\infty}\|J^{\sigma-1}\partial_x\big(\langle x\rangle^{m+\frac{1}{2}+\theta-\frac{l}{2}} v\big)\|_{2}\|\langle x\rangle^{m+\frac{1}{2}+\theta-\frac{l}{2}} J^{\sigma}v\|_{2}\\
				\leq c & \|\partial_x v\|_{\infty}\|\langle x\rangle^{m+\frac{1}{2}+\theta-\frac{l}{2}} J^{\sigma}v\|_{2}^2.
			\end{aligned}   
		\end{equation*}
		This completes the estimate of $\mathcal{N}_{12}$, for any $\sigma>\frac{1}{2}$.  Finally, by integration by parts we deduce
		\begin{equation*}
			\begin{aligned}
				\mathcal{N}_{13}=&-\frac{1}{2}
				\int_{\mathbb{R}} \big(J^{\sigma} v\big)^2\Big(x^{2m+1}\langle x \rangle^{2\theta}_N\langle x\rangle^{-l}v \Big)'\, dx\\
				=&-\frac{1}{2}
				\int_{\mathbb{R}} \big(\langle x\rangle^{m+\frac{1}{2}+\theta-\frac{l}{2}} J^{\sigma} v\big)^2 \langle x\rangle^{-2m-1-2\theta+l}\Big(x^{2m+1}\langle x \rangle^{2\theta}_N\langle x\rangle^{-l}v \Big)'\, dx.   
			\end{aligned}    
		\end{equation*} 
		From which it follows
		\begin{equation*}
			\begin{aligned}
				\big| \mathcal{N}_{13}\big|\leq c &\big(\|v\|_{\infty}+\|\partial_x v\|_{\infty}\big)\|\langle x\rangle^{m+\frac{1}{2}+\theta-\frac{l}{2}} J^{\sigma} v\|_{2}^2.
			\end{aligned}    
		\end{equation*} 
		The proof is complete.
	\end{proof}	
	\section{Proof of Theorem \ref{main}}
	In this section we prove Theorem \ref{main}. The argument is based on an iterative process. At each step, the decay of the solution at two different times, combined with a weighted energy estimate, yields a local smoothing effect of order $a/2$ at the price of half a power of the weight. Since this information holds for almost every time in the corresponding time interval, it provides a new two-time condition, which is the starting point of the next step. After finitely many steps, the interpolation inequalities in Section 2 give the regularity $H^{a\alpha}(\mathbb{R})$ at some time, and the local theory extends it to the whole interval $[-T,T]$.
	Part (i) follows by energy estimates and approximating by smooth solutions, using continuous dependence, with initial data in the Schwartz class (where Theorem \ref{AA2} holds). For the sake of brevity, and because this case is contained in our arguments below, we omit further details here. However, it is worth emphasizing that in this case, the estimate of the nonlinear term produces a factor $\|u\|_{3}^3$, which is controlled by the embedding $H^{1/6}(\mathbb{R})\hookrightarrow L^{3}(\mathbb{R})$. This justifies the required regularity  $s\geq 1/6$ in (i). 
	
	\medskip

	Let us deduce (ii) in  Theorem \ref{main}. We split the proof into several steps according to the size of the weight parameter $\alpha\in (s/a,5/2+a)$. We will also frequently use the approximation $\langle x\rangle_N$ in \eqref{WNwightApp}.
	The cases $\alpha\le 3/2$ contain the main ideas of the argument. In the subsequent cases, the weight is decomposed as a power of $x$ times a truncated weight $\langle x\rangle_N^{2\theta}$, the identity \eqref{generalconmmutator} produces additional dispersive terms, and these terms are controlled by means of the results obtained in the previous cases.
	\\ \\
	\noindent \underline{\bf Case $\alpha\in (s/a,1]$}. We have $3/2<s<a\alpha<2$. By the persistence result in weighted spaces in Theorem \ref{AA2}, the fact that $|x|^{\alpha}u(t_j)\in L^2(\mathbb{R})$, $j=1,2$ implies
	\begin{equation}\label{classn1}
		u\in C([-T,T]:H^s(\mathbb{R}))\cap C([-T,T]:L^2(|x|^{2s/a}\, dx)).    
	\end{equation}
	We first show that the weighted energy estimate produces the local smoothing effect $\langle x\rangle^{\alpha-1/2}|D_x|^{a/2}u\in L^{2}([t_1,t_2]:L^2(\mathbb{R}))$.
	We multiply the equation in \eqref{dgbo} by $x\langle x\rangle_N^{2\alpha-1}\, u$ and integrate over the space variable to get the identity
	\begin{equation}\label{energyN11}
		\begin{split}
			&\frac{1}{2}\frac{d}{dt}
			\int_{\mathbb{R}} u^{2}\,x \langle x\rangle_N^{2\alpha-1} \,dx-
			\int_{\mathbb{R}} |D_x|^a\partial_xu \,u \,x \langle x\rangle_N^{2\alpha-1} \,dx \\
			&+\int_{\mathbb{R}} u\,\partial_{x}u\,u\,x \langle x\rangle_N^{2\alpha-1} \,dx
			= 0.
		\end{split}
	\end{equation}
	Notice that since $s/a>1/2$, the above identity is justified by \eqref{classn1} and $\langle x\rangle_N\in L^{\infty}(\mathbb{R})$. By integration by parts (or \eqref{nonlinearE1}), we get
	\begin{equation*}	
		\begin{split}
			\Big|\int_{\mathbb{R}} u^{2}\,\partial_{x}u\,x\langle x\rangle^{2\alpha-1}_N \,dx\Big|\leq c\int_{\mathbb{R}}|u|^3\langle x\rangle^{2\alpha-1}_N\, dx\leq & c\|u\|_{\infty}\|\langle x \rangle^{\alpha-1/2}_N u\|_2^2\\
			\leq c &\|u\|_{\infty}\|\langle x \rangle^{s/a} u\|_2^2,
		\end{split}
	\end{equation*}
	which is controlled by \eqref{classn1} with constant independent of $N$.  To bound the second term in \eqref{energyN11}, we use Lemma \ref{prop1a} to get
	\begin{equation}\label{energyM1}
		\begin{split}	
			\int_{\mathbb{R}} |D_x|^a\partial_xu \,u \,x\langle x\rangle_N^{2\alpha-1} \,dx
			&=-(\frac{a+1}{2}\Big)\int_{\mathbb{R}}\Big||D_x|^{a/2}u\Big|^2\,(x\langle x\rangle_N^{2\alpha-1})'\,dx \\
			&\hskip15pt-\frac{1}{2} \int_{\mathbb{R}}|D_x|^{\sigma}|D_x|^{-\sigma}u\,\tilde{R}_0(a+1)\,|D_x|^{\sigma}|D_x|^{-\sigma}u\,dx,
		\end{split}
	\end{equation}
	with  $\sigma=1-a/2\in(0,1/2)$. To bound the last term above, we apply \eqref{eq98b} to get
	\begin{equation}
		\label{energyM2}
		\begin{split}
			|\int_{\mathbb{R}}|D_x|^{\sigma}|D_x|^{-\sigma}u\,\tilde{R}_0(a+1)\,|D_x|^{\sigma}|D_x|^{-\sigma}u\,dx|
			&= |\int_{\mathbb{R}}|D_x|^{-\sigma}u\,|D_x|^{\sigma}\tilde{R}_0(a+1)\,|D_x|^{\sigma}|D_x|^{-\sigma}u\,dx|\\
			&\leq 
			\||D_x|^{-\sigma}u\|_2\;
			\||D_x|^{\sigma}\tilde{R}_0(a+1)|D_x|^{\sigma}|D_x|^{-\sigma}u\|_2\\
			&\leq c \|\widehat{|D_x|^{a+1+2\sigma}(x\langle x\rangle_N^{2\alpha-1})}\|_{L^1_{\xi}}\,\||D_x|^{-\sigma}u\|^2_2.
		\end{split}
	\end{equation}
	Let $0<\sigma<\beta<1/2$, $\phi\in C^{\infty}_0(\mathbb{R})$ such that $\phi=1$ on an neighborhood of the origin. Using Plancherel's identity, H\"older's inequality, and the Hardy-Littlewood-Sobolev inequality, we obtain 
	\begin{equation}\label{energyInequalityM}
		\begin{aligned}
			\||D_x|^{-\sigma}u\|_2\leq & c\||\xi|^{-\sigma }\phi\widehat{u}\|_{2}+\||\xi|^{-\sigma }(1-\phi)\widehat{u}\|_{2}\\
			\leq & c\||\xi|^{-\sigma }\phi\|_{\frac{1}{\beta}}\|\widehat{u}\|_{\frac{2}{1-2\beta}}+\||\xi|^{-\sigma}(1-\phi)\|_{\infty}\|\widehat{u}\|_{2}\\
			\leq & c\|J_{\xi}^{\beta}\widehat{u}\|_2\\
			\leq &c \|\langle x\rangle^{1/2^{-}}u\|_2 <\infty,    
		\end{aligned}
	\end{equation}
	which holds provided \eqref{classn1} and that $s/a>1/2$. Next, since $\sigma=1-a/2$ one has
	\begin{equation*}
		\|\widehat{|D_x|^{a+1+2\sigma}(x\langle x\rangle_N^{2\alpha-1})}\|_{L^1_{\xi}}= \|\widehat{|D_x|^{3}(x\langle x\rangle_N^{2\alpha-1})}\|_{L^1_{\xi}}
		= \|\widehat{\partial_x^{3}(x\langle x\rangle_N^{2\alpha-1})}\|_{L^1_{\xi}}.
	\end{equation*}
	To bound the above expression for $\alpha<3/2$, with a bound independent of $N$, we see that if 
	\begin{equation}
		\label{G1}
		f_{\alpha}(x)=\partial_x^3(x\langle x\rangle_N^{2\alpha-1}),
	\end{equation}
	then 
	\begin{equation}
		\label{G2}
		f_{\alpha}(x)\sim |x|^{2\alpha-3},\;\;\;|x|\uparrow \infty,
	\end{equation}
	with 
	\begin{equation}
		\label{G3}
		\partial_x^2f_{\alpha}\in L^1(\mathbb{R})\;\;\;\;\;\;\;\;\;\text{so}\;\;\;\;\;\;\;\;\;|\widehat{f_{\alpha}}(\xi)|\leq c\frac{\|\partial_x^2f_{\alpha}\|_1}{|\xi|^2}
	\end{equation}
	and 
	\begin{equation}
		\label{G4}
		f_{\alpha}(x)\sim |x|^{2\alpha-3}\;\;\;\;\;\text{so}\;\;\;\;\;\widehat{f_{\alpha}}(\xi)\sim |\xi|^{2-2\alpha},\;\;\;\xi\sim 0.
	\end{equation}
	Hence, $\widehat{f}_{\alpha}\in L^1(\mathbb{R})$ if $2\alpha-2<1$, i.e. $\alpha<3/2$, with a bound independent of $N$.

	On the other hand, since $\langle x \rangle_N$ is an even, non-decreasing function of $x\geq 0$, we have that $x\big(\partial_x\langle x\rangle_N\big)\geq 0$, then
	\begin{equation*}
		\begin{aligned}
			\big(\langle x\rangle_N^{2\alpha-1} x\big)'=\langle x\rangle_N^{2\alpha-1}+(2\alpha-1) \langle x\rangle_N^{2\alpha-1-1}\, x\partial_x\big(\langle x\rangle_N\big)\geq \langle x\rangle_{N}^{2\alpha-1}.    
		\end{aligned}    
	\end{equation*}
	Therefore, integrating in the time interval $[t_1,t_2]$ identity \eqref{energyN11}, using previous estimates and the two time conditions $|x|^{\alpha}u(t_1),|x|^{\alpha}u(t_2)\in L^2(\mathbb{R})$, it follows
	\begin{equation*}
		\int_{t_1}^{t_2}\,\int_{\mathbb{R}} \langle x\rangle_N^{2\alpha-1}\,\Big||D_x|^{a/2}u(x,t)\Big|^2\,dx\,dt\leq c<\infty,
	\end{equation*}
	where $c>0$ is independent of $N$. Thus, taking $N\to \infty$, we deduce
	\begin{equation*}
		\langle x\rangle^{\alpha-1/2}\,|D_x|^{a/2}u(x,t)\in L^{2}([t_1,t_2]:L^2(\mathbb{R})).    
	\end{equation*}
	Using the previous result, and reapplying the argument implemented above, with the weight $\langle x\rangle^{2\alpha}_N$, we conclude that
	\begin{equation}\label{z}
		\langle x\rangle^{\alpha}\,u(x,t)\in L^{\infty}([t_1,t_2]:L^2(\mathbb{R})).
	\end{equation}
	Hence, by Corollary \ref{pseudo-cor1}, we get
	\begin{equation}\label{Jconcl1}
		\langle x\rangle^{\alpha-1/2}\,J^{a/2}u(x,t)\in L^{2}([t_1,t_2]:L^2(\mathbb{R})).
	\end{equation}
	In particular, there exist two times $t_1\leq t_3 <t_4\leq t_2$ such that
	\begin{equation}\label{newtwotimes}
		\langle x\rangle^{\alpha-1/2}\,J^{a/2}u(x,t_j)\in L^2(\mathbb{R}), \qquad j=3,4.
	\end{equation}
	We now repeat the previous argument at the level of $J^{a/2}u$, with the times $t_3$ and $t_4$ in \eqref{newtwotimes} in place of $t_1$ and $t_2$.
	Next, we apply $J^{a/2}$ to the equation in \eqref{dgbo}, multiplying then by $   x\langle x\rangle^{2\alpha-2} J^{a/2}u$ and integrate in the spatial variable, we get
	\begin{equation}\label{energyN12}
		\begin{split}
			&\frac{1}{2}\frac{d}{dt}
			\int_{\mathbb{R}} (J^{a/2}u)^{2}\,x \langle x\rangle^{2\alpha-2} \,dx-
			\int_{\mathbb{R}} |D_x|^a\partial_x J^{a/2}u \,J^{a/2}u \,x \langle x\rangle^{2\alpha-2} \,dx \\
			&+\int_{\mathbb{R}} J^{a/2}(u\,\partial_{x}u)\,J^{a/2}u\,x\langle x\rangle_N^{2\alpha-2} \,dx
			= 0.
		\end{split}
	\end{equation}
	For this energy estimate, we do not use an approximation argument with weight $\langle x\rangle_N$. Instead, \eqref{energyN12} can be justified from \eqref{z} and \eqref{Jconcl1}, and Theorem \ref{AA2}, together with approximation by functions with initial data in the Schwartz class.

	Since $1/2<a/2<1$, we apply \eqref{nonlinearE2} with $\sigma=a/2$, $\theta=\alpha-1/2$, $m=0$, $l=1$ to deduce
	\begin{equation}\label{energy12}	
		\begin{split}
			\Big|\int_{\mathbb{R}} J^{a/2}(u\,\partial_{x}u)\,J^{a/2}u\,x\langle x\rangle^{2\alpha-2} \,dx\Big|\leq & c(\|J^{a/2} u\|_{2}+\|\partial_x u\|_{\infty})\|\langle x \rangle^{\alpha-1/2}J^{a/2}u\|_2^2\\
			\leq & c\|J^{s} u\|_{2}\|\langle x \rangle^{\alpha-1/2}J^{a/2}u\|_2^2
		\end{split}
	\end{equation}
	which, upon integration between $t_3$ and $t_4$, gives the desired control in terms of \eqref{z} and the $H^s$-norm of $u$ with $s\geq a/2$ and $s>3/2$.

	To bound the second term in \eqref{energyN12}, we first assume that $\alpha<1$ (i.e., $\alpha-1/2<1/2$) and use Lemma \ref{prop1a} to get
	\begin{equation}\label{energy12b}
		\begin{split}	
			\int_{\mathbb{R}} |D_x|^a\partial_x J^{a/2}u \, J^{a/2}u \,x \langle x\rangle^{2\alpha-2} \,dx
			&=-\Big(\frac{a+1}{2}\Big)\int_{\mathbb{R}}\Big||D_x|^{a/2}J^{a/2}u\Big|^2\,(x\langle x\rangle^{2\alpha-2})'\,dx \\
			&\hskip15pt-\frac{1}{2} \int_{\mathbb{R}} J^{a/2}u\,\tilde{R}_0(a+1)\,J^{a/2}u\,dx, 
		\end{split}
	\end{equation}
	we have
	\begin{equation}\label{EstimateNOWeight}
		\|\tilde{R}_0(a+1)u\|_2\leq c \|\widehat{|D_x|^{a+1}(x\langle x\rangle^{2\alpha-2})}\|_{L^1_{\xi}}\,\|J^{a/2}u\|_2.    
	\end{equation}
	Using that $\|\widehat{f}\|_{1}\leq c \|\langle \xi\rangle^{\frac{1}{2}^{+}}\widehat{f}\|_{2}\leq c \|f\|_{H^1}$, we write $|D_x|^{a+1}=-\mathcal{H}|D_x|^{a-2}\partial_x^3$ to get
	\begin{equation}\label{L1estimate1}
		\begin{aligned}
			\|\widehat{|D_x|^{a+1}(x\langle x\rangle^{2\alpha-2})}\|_{L^1_{\xi}}
			\leq c & \||D_x|^{a-2}\partial_x^3(x\langle x\rangle^{2\alpha-2})\|_{2}+c\||D_x|^{a-2}\partial_x^4(x\langle x\rangle^{2\alpha-2})\|_{2}.
		\end{aligned}    
	\end{equation}
	Set $g(x)=\partial_x^3( x\langle x\rangle^{2\alpha-2})$. Notice that $|g(x)|\leq c \langle x\rangle^{2\alpha-4}$, with $c>0$ independent of $N$, and $\widehat{g}(0)=0$. Then, using that $\widehat{g}(\xi)=\xi\int_0^1 \frac{d}{d\xi}\widehat{g}(\sigma\xi)\, d\sigma $, and a function $\chi\in C^{\infty}_0(\mathbb{R})$ such that $\chi=1$ on an neighborhood of the origin, it follows that
	\begin{equation}\label{L1estimate2}
		\begin{aligned}
			\||D_x|^{a-2}\partial_x^3(x\langle x\rangle^{2\alpha-2})\|_{2}\leq c & \||\xi|^{a-2}\chi\widehat{g}\|_{2}+c\||\xi|^{a-2}(1-\chi)\widehat{g}\|_{2} \\
			\leq & c\int_0^1\||\xi|^{a-1}\chi  \frac{d}{d\xi}\widehat{g}(\sigma\xi)\|_{2}\, d\sigma+c\||\xi|^{a-2}(1-\chi)\|_{\infty}\|\widehat{g}\|_{2}\\
			\leq c&\int_0^1 \sigma^{-\frac{1}{2}}\||\xi|^{a-1}\chi\|_{\infty}\|\frac{d}{d\xi}\widehat{g}\|_{2}\, d\sigma +c\|g\|_{2}\\
			\leq c & \|\langle x\rangle g\|_{2}\\
			\leq c & \|\langle x\rangle^{2\alpha-3}\|_{2}<\infty,
		\end{aligned}    
	\end{equation}
	which holds for $\alpha<5/4$. Similarly, $ \||D_x|^{a-2}\partial_x^4(x\langle x\rangle^{2\alpha-2})\|_{2}\leq c \|\langle x\rangle^{2\alpha-4}\|_{2}$, with constant $c>0$ independent of $N$. This completes the estimate of \eqref{EstimateNOWeight}.
	
	On the other hand, if $\alpha=1$, using that $\partial_x |D_x|^a$ is a skew-symmetric operator, by \eqref{generalconmmutator} with $m=1$, it follows
	\begin{equation*}
		\begin{aligned}
			\int_{\mathbb{R}} |D_x|^a\partial_x J^{a/2}u \,J^{a/2}u \,x  \,dx =&-\frac{1}{2}\int_{\mathbb{R}}  J^{a/2}u \,[\partial_x|D_x|^{a};x]J^{a/2}u\, \,dx\\
			=&-\Big(\frac{a+1}{2}\Big)\int_{\mathbb{R}}  J^{a/2}u \, |D_x|^aJ^{a/2}u\, \,dx\\
			=&-\Big(\frac{a+1}{2}\Big)\int_{\mathbb{R}}  \Big||D_x|^{a/2}J^{a/2}u\Big|^2 \,dx.
		\end{aligned}    
	\end{equation*}
	Additionally, we also have that for $0<\alpha-1/2\leq 1/2$, there exists $c_{\alpha}>0$ such that
	\begin{equation*}
		(x\langle x\rangle^{2\alpha-2})'\geq c_{\alpha} \langle x \rangle^{2\alpha-2}.   
	\end{equation*}
	Integrating in the time interval $[t_3,t_4]$ the identity \eqref{energyN12}, using \eqref{newtwotimes} and the previous estimates, we find
	\begin{equation}\label{energy13b}
		\int_{t_3}^{t_4}\,\int_{\mathbb{R}} \langle x\rangle^{2\alpha-2}\,\Big||D_x|^{a/2}J^{a/2}u(x,t)\Big|^2\,dx\,dt<\infty.
	\end{equation}
	In other words, for $\frac{s}{a}-\frac{1}{2}<\alpha-1/2\leq \frac{1}{2}$, it follows
	\begin{equation*}
		\langle x\rangle^{\alpha-1}\,|D_x|^{a/2}J^{a/2}u(x,t)\in L^{2}([t_3,t_4]:L^2(\mathbb{R})).    
	\end{equation*}
	By Corollary \ref{pseudo-cor2} and \eqref{classn1}, we get
	\begin{equation*}
		\langle x\rangle^{\alpha-1}\,J^{a}u(x,t)\in L^{2}([t_3,t_4]:L^2(\mathbb{R})).
	\end{equation*}
	The interpolation inequality \eqref{interpol-complete} between the above result and \eqref{z} yields
	\begin{equation*}
		J^{a \alpha}u(x,t)\in L^{2}([t_3,t_4]:L^2(\mathbb{R})).
	\end{equation*}
	Thus, there exists a time $t^{\ast}\in [t_3,t_4]$ such that $J^{a \alpha}u(x,t^{\ast})\in L^2(\mathbb{R})$. Hence, by the local well-posedness in $H^s$-spaces, we get  $u\in C([-T,T]:H^{a\alpha}(\mathbb{R}))$.
	\\ \\
	Summarizing, in the case $\alpha\in(s/a,1]$ a few applications of the weighted energy estimate are enough to reach the regularity $H^{a\alpha}(\mathbb{R})$. For larger weights the same scheme is iterated, and the number of steps increases with $\alpha$.
	
	\underline{\bf Case $\alpha\in (1,3/2]$}. We have $3/2<s<a\alpha$, $1<a<2$. We first consider $\alpha\in(1,3/2)$. By the result in the previous case and Theorem \ref{AA2}, we can assume 
	\begin{equation}\label{class2}
		u\in C([-T,T]:H^{\max\{s,\, a\}}(\mathbb{R}))\cap C([-T,T]:L^2(|x|^{2}\, dx)).    
	\end{equation}
	and at two different times $t_1,t_2\in [-T,T]$, $|x|^{\alpha}u(t_j)\in L^2(\mathbb{R})$.
	
	We consider the energy estimate 
	\begin{equation}\label{energy11}
		\begin{split}
			&\frac{1}{2}\frac{d}{dt}
			\int_{\mathbb{R}} u^{2}\,x \langle x\rangle_N^{2\alpha-1} \,dx-
			\int_{\mathbb{R}} |D_x|^a\partial_xu \,u \,x \langle x\rangle_N^{2\alpha-1} \,dx \\
			&+\int_{\mathbb{R}} u\,\partial_{x}u\,u\,x \langle x\rangle_N^{2\alpha-1} \,dx
			= 0.
		\end{split}
	\end{equation}
	By the identity \eqref{eq1abc}, we have
	\begin{equation}\label{energy12c}
		\begin{split}	
			\int_{\mathbb{R}} |D_x|^a\partial_xu \,u \,x\langle x\rangle_N^{2\alpha-1} \,dx
			&=-(\frac{a+1}{2}\Big)\int_{\mathbb{R}}\Big||D_x|^{a/2}u\Big|^2\,(x\langle x\rangle_N^{2\alpha-1})'\,dx \\
			&\hskip15pt-\frac{1}{2} \int_{\mathbb{R}}|D_x|^{\sigma}|D_x|^{-\sigma}u\,\tilde{R}_0(a+1)\,|D_x|^{\sigma}|D_x|^{-\sigma}u\,dx,
		\end{split}
	\end{equation}
	with  $\sigma=1-a/2\in(0,1/2)$. To bound the last term in \eqref{energy12c} we write
	\begin{equation}
		\label{energy15}
		\begin{split}
			|\int_{\mathbb{R}}|D_x|^{\sigma}|D_x|^{-\sigma}u\,\tilde{R}_0(a+1)\,|D_x|^{\sigma}|D_x|^{-\sigma}u\,dx|
			&= |\int_{\mathbb{R}}|D_x|^{-\sigma}u\,|D_x|^{\sigma}\tilde{R}_0(a+1)\,|D_x|^{\sigma}|D_x|^{-\sigma}u\,dx|\\
			&\leq 
			\||D_x|^{-\sigma}u\|_2\;
			\||D_x|^{\sigma}\tilde{R}_0(a+1)|D_x|^{\sigma}|D_x|^{-\sigma}u\|_2\\
			&\leq c \|\widehat{|D_x|^{a+1+2\sigma}(x\langle x\rangle_N^{2\alpha-1})}\|_{L^1_{\xi}}\,\||D_x|^{-\sigma}u\|^2_2.
		\end{split}
	\end{equation}
	By \eqref{energyInequalityM}, we have
	\begin{equation*}
		\||D_x|^{-\sigma}u\|_2\leq c \|\langle x\rangle^{1/2}u\|_2 <\infty,
	\end{equation*}
	which is justified by \eqref{class2}. By the ideas in \eqref{G1}-\eqref{G4}, it follows $\|\widehat{|D_x|^{a+1+2\sigma}(x\langle x\rangle_N^{2\alpha-1})}\|_{L^1_{\xi}}\leq c <\infty$ with constant independent of $N$.

	To control the nonlinear term in \eqref{energy11}, by integration by parts 
	\begin{equation*}
		\begin{aligned}
			\Big|\int_{\mathbb{R}} u\,\partial_{x}u\,u\,x \langle x\rangle_N^{2\alpha-1} \,dx\Big|\leq  c \|\langle x\rangle_N^{\alpha-1/2} u\|_2^2\,\|u\|_{\infty}\leq  c \|\langle x\rangle u\|_2^2\,\|u\|_{\infty}
		\end{aligned}    
	\end{equation*}
	which is bounded by \eqref{class2}. Since $(x\langle x \rangle_N^{2\alpha-1})'\geq c\langle x\rangle_N^{2\alpha-1}$ with $c>0$ independet of $N$, integrating \eqref{energy11} over $[t_1,t_2]$, we collect the above estimates to take $N\to \infty$, deducing that 
	\begin{equation*}
		\langle x\rangle^{\alpha-1/2}\,|D_x|^{a/2}u(x,t)\in L^{2}([t_1,t_2]:L^2(\mathbb{R})).    
	\end{equation*}   
	Using the above results and replacing the energy estimate with weight $\langle x \rangle^{2\alpha}$, we arrive at
	\begin{equation}\label{alphadecayI}
		\langle x\rangle^{\alpha}\, u(x,t)\in L^{\infty}([t_1,t_2]:L^2(\mathbb{R})). 
	\end{equation}   
	Hence, Corollary \ref{pseudo-cor1} establishes
	\begin{equation*}
		\langle x\rangle^{\alpha-1/2}\,J^{a/2}u(x,t)\in L^{2}([t_1,t_2]:L^2(\mathbb{R})).    
	\end{equation*}   
	Thus, there exist two times $t_1\leq t_3 <t_4\leq t_2$ such that
	\begin{equation*}
		\langle x\rangle^{\alpha-1/2}\,J^{a/2}u(x,t_j)\in L^2(\mathbb{R}), \qquad j=3,4.
	\end{equation*}
	Since $\alpha-1/2\in (1/2,1)$, we can now return to the previous case of weights in $(s/a,1]$ and repeat the argument, with the regularity increased by $J^{a/2}$ and the decay exponent replaced by $\alpha-1$. For example, the first energy estimate replaces $u$ with $J^{a/2}u$ and the weight with $\alpha-1$. Since the new estimates involve powers of $J$, the dispersive contribution remains unchanged, whereas the nonlinear contribution can be controlled by Lemma \ref{Nonlinear estimates} with $m=0$, $l=0$, $\theta=\alpha-1\in (0,1/2)$ as follows
	\begin{equation*}
		\begin{aligned}
			\Big|\int_{\mathbb{R}}  J^{a/2}\big(u\partial_x u\big) \big(x \langle x\rangle^{2\alpha-1}_N J^{a/2}u\big)\, dx\Big|\leq & c\big(\|J^{a/2}u\|_{2}+\|\partial_x u\|_{\infty}\big)\|\langle x\rangle^{1/2+\alpha-1} J^{a/2} u\|_{2}^2\\
			\leq &  c\|J^{(3/2)^{+}}u\|_{2}\|\langle x\rangle^{\alpha-1/2} J^{a/2} u\|_{2}^2.
		\end{aligned} 
	\end{equation*}
	Summarizing, one deduces that
	\begin{equation*}
		\langle x\rangle^{\alpha-1}|D_x|^{a/2}\,J^{a/2}u(x,t)\in L^{2}([t_3,t_4]:L^2(\mathbb{R})).  
	\end{equation*}  
	and
	\begin{equation*}
		\langle x\rangle^{\alpha-1/2}\,J^{a/2}u(x,t)\in L^{\infty}([t_3,t_4]:L^2(\mathbb{R})). 
	\end{equation*}  
	Hence, by Corollary \ref{pseudo-cor1}
	\begin{equation*}
		\langle x\rangle^{\alpha-1}\,J^{a}u(x,t)\in L^{2}([t_3,t_4]:L^2(\mathbb{R})).  
	\end{equation*}  
	Moreover, taking two times $t_3\leq t_5<t_6\leq t_4$ such that $\langle x\rangle^{\alpha-1}\,J^{a}u(x,t_j)\in L^2(\mathbb{R})$, $j=5,6$, it follows
	\begin{equation*}
		\langle x\rangle^{\alpha-3/2}\,|D_x|^{a/2}J^{a}u(x,t)\in L^{2}([t_5,t_6]:L^2(\mathbb{R}))
	\end{equation*}  
	and
	\begin{equation*}
		\langle x\rangle^{\alpha-1}\,J^{a}u(x,t)\in L^{\infty}([t_5,t_6]:L^2(\mathbb{R})).  
	\end{equation*}  
	Thus, Corollary \ref{pseudo-cor2} yields
	\begin{equation*}
		\langle x\rangle^{\alpha-3/2}\,J^{3a/2}u(x,t)\in L^{2}([t_5,t_6]:L^2(\mathbb{R})).
	\end{equation*}  
	The interpolation inequality \eqref{interpol-complete} between the above result and \eqref{alphadecayI} establishes 
	\begin{equation*}
		J^{a \alpha}u(x,t)\in L^{2}([t_5,t_6]:L^2(\mathbb{R})).
	\end{equation*}
	Consequently, there exists a time $t^{\ast}\in [t_5,t_6]$ such that $J^{a \alpha}u(x,t^{\ast})\in L^2(\mathbb{R})$. The LWP theory yields $u\in C([-T,T]:H^{a\alpha}(\mathbb{R}))$.
	\\ \\
	The case $\alpha=3/2$ is not covered by the previous argument, since the bound for $\widehat{\partial_x^{3}(x\langle x\rangle_N^{2\alpha-1})}$ in $L^{1}(\mathbb{R})$ requires $\alpha<3/2$, see \eqref{G1}-\eqref{G4}. In this case we use instead the polynomial weight $x^{3}$ and the commutator identity \eqref{I1}.
	
	Next, we consider the case $\alpha=3/2$. Formally, we multiply the equation in \eqref{dgbo} by $x^3 u$ and integrate over the spatial variable. As in the previous cases, we just need to consider the linear terms:
	\begin{equation}
		\label{3/2}
		A_{0} = -\int_{\mathbb{R}} \partial_{x} |D_{x}|^{a} u \cdot x^{3} u \, dx
	\end{equation}
	Using the identity (see Proposition \ref{GeneralConmutator}) 
	\begin{equation} \label{I1}
		x|D_x|^a\partial_xu=|D_x|^a\partial_x(xu)-(a+1)|D_x|^au,
	\end{equation}
	we see that:
	\begin{equation*}
		\begin{aligned}
			A_{0} &= -\int_{\mathbb{R}} \partial_{x} |D_{x}|^{a} (xu) \cdot xu \cdot x \, dx + (a + 1) \int_{\mathbb{R}} |D_{x}|^{a} u \cdot (xu)\, x \, dx =: A_{1} + A_{2}.
		\end{aligned}
	\end{equation*}
	Now, by integration by parts
	\begin{equation*}
		A_{1} = \int_{\mathbb{R}} |D_{x}|^{a} (xu) \partial_{x} (xu) x \, dx + \int_{\mathbb{R}} |D_{x}|^{a} (xu) xu \, dx,
	\end{equation*}
	and then using the identity
	\begin{equation}\label{I2}
		x|D_x|^af=|D_x|^a(xf)-a|D_x|^{a-1}\mathcal Hf,
	\end{equation}
	it follows that:
	\begin{equation*}
		\begin{split}
			A_1&= \int_{\mathbb{R}} xu |D_{x}|^{a} \partial_{x} (xu) x \, dx - a \int_{\mathbb{R}} xu |D_{x}|^{a} (xu) \, dx + \int_{\mathbb{R}} \left(|D_{x}|^{a/2} (xu)\right)^{2} \, dx \\
			&= -A_{1} - a \int_{\mathbb{R}} xu |D_{x}|^{a} (xu) \, dx + \int_{\mathbb{R}} \left(|D_{x}|^{a/2} (xu)\right)^{2} \, dx.
		\end{split}
	\end{equation*}
	Thus,
	\begin{equation*}
		A_{1} = \frac{1 - a}{2} \int_{\mathbb{R}} \left(|D_{x}|^{a/2} (xu)\right)^{2} \, dx,
	\end{equation*}
	and by \eqref{I2}
	\begin{equation*}
		\begin{aligned}
			A_{2} &= (a + 1) \int_{\mathbb{R}} |D_{x}|^{a} (xu) \cdot xu \, dx - (a + 1)a \int_{\mathbb{R}} |D_{x}|^{a - 1} \mathcal{H}u \cdot xu \, dx \\
			&= (a + 1) \int_{\mathbb{R}} \left(|D_{x}|^{a/2} (xu)\right)^{2} \, dx - (a + 1)a \int_{\mathbb{R}} |D_{x}|^{a - 1} \mathcal{H}u \cdot xu \, dx.
		\end{aligned}
	\end{equation*}
	Gathering the above information, we conclude that:
	\begin{equation}
		\label{3/2/b}
		A_{0} = \left(\frac{a}{2} + \frac{3}{2}\right) \int_{\mathbb{R}} \left(|D_{x}|^{a/2} (xu)\right)^{2} \, dx - (a + 1)a \int_{\mathbb{R}} |D_{x}|^{a - 1} \mathcal{H}u \cdot xu \, dx.
	\end{equation}
	The last term above can be bounded by the previous result, so using a familiar argument.  We remark that  the above formal proof can be made rigorously justified by using an approximation argument and the commutator estimates found in D.~Li \cite[Corollary 1.4]{Dongli}, which are extensions of those described in Theorem \ref{kprefi} common to the homogeneous operators $|D_{x}|^{p},\;\mathcal H |D_x|^p$, $p > 0$. Thus, integrating between $t_1$ and $t_2$ one obtains
	\begin{equation*}
		|D_x|^{a/2}(x u)\in L^2([t_1,t_2]:L^2(\mathbb{R})).   
	\end{equation*}
	Applying Proposition \ref{JandDprop} with $m=1$, $\theta=0$, and $\beta=a/2$, and using \eqref{class2} it follows
	\begin{equation*}
		\langle x\rangle J^{a/2}u \in L^2([t_1,t_2]:L^2(\mathbb{R})).  
	\end{equation*}
	We conclude the proof for the case $\alpha=3/2$. 
	\\ \\
	With the case $\alpha=3/2$ at hand, the function $xu$ plays the role of $u$ in the next cases, and the additional decay is measured by means of the truncated weights $\langle x\rangle_N^{2\theta}$.
	
	\underline{\bf Case $\alpha\in (3/2,2]$}. Together  with the result from the previous case, Theorem \ref{AA2} allows us to assume
	\begin{equation}\label{class3}
		u\in C([-T,T]:H^{3a/2}(\mathbb{R}))\cap C([-T,T]:L^2(|x|^{3}\, dx)).    
	\end{equation}
	and at two different times $t_1,t_2\in [-T,T]$, $|x|^{\alpha}u(t_j)\in L^2(\mathbb{R})$. We also assume 
	\begin{equation}\label{Extraclass}
		|D_x|^{a/2}(x u)\in L^{\infty}([t_1,t_2]:L^2(\mathbb{R})).   
	\end{equation}
	We first consider $\alpha=(3/2,2)$ and write $\alpha=3/2+\theta$, $\theta\in (0,1/2)$. Multiply the equation by $x^3u\langle x\rangle_N^{2\theta}$ with $\theta\in(0,1/2)$, use the identity \eqref{generalconmmutator} with $m=1$ and integrate in the space variable to get
	\begin{equation}
		\label{EE1}
		\begin{aligned}
			& \frac{1}{2} \frac{d\;}{dt}\int_{\mathbb{R}} (xu)^2x\langle x \rangle_N^{2\theta}\,dx -\int_{\mathbb{R}} |D_x|^a\partial_x(xu)\,xu\,x\langle x\rangle_N^{2\theta} dx\\
			&+(a+1) \int_{\mathbb{R}} |D_x|^au\, xu\,x\langle x\rangle_N^{2\theta} dx + \int_{\mathbb{R}} xu\partial_xu \,xu\,x\langle x\rangle_N^{2\theta} dx=0.
	\end{aligned}\end{equation}
	By integration by parts and \eqref{class3} we have 
	\begin{equation*}
		\begin{aligned}
			|\int_{\mathbb{R}} xu\partial_xu \,xu\,x\langle x\rangle_N^{2\theta} dx|\leq c \int_{\mathbb{R}} |u|^3\langle x\rangle_N^{2+2\theta}dx \leq & c \| \langle x\rangle_N^{1+\theta}u\|^2_2\,\|u\|_{\infty}\\
			\leq &c \| \langle x\rangle^{3/2}u\|^2_2\,\|u\|_{\infty}.
		\end{aligned}    
	\end{equation*}
	To handle the second term in \eqref{EE1}, we use \eqref{eq1abc} to write
	\begin{equation}\label{EE2}
		\begin{split} 
			\int_{\mathbb{R}} |D_x|^a\partial_x(xu)\,xu\,x\langle x\rangle_N^{2\theta} dx
			&= -\Big( \frac{a+1}{2}\Big)\int_{\mathbb{R}} ||D_x|^{a/2}(xu)|^2\,(x\langle x\rangle_N^{2\theta})'\,dx\\
			& \hskip12pt -\frac{1}{2} \int_{\mathbb{R}} (xu)\,\tilde{R}_0(a+1)(xu)\,dx.
		\end{split}
	\end{equation}
	To bound the last term in \eqref{EE2}, we apply \eqref{eq98} with $\sigma=1-a/2\in(0,1/2)$ so that $a+1+2\sigma=3$  to deduce
	\begin{equation*}
		\begin{aligned} 
			|\int_{\mathbb{R}} (xu)\,\tilde{R}_0(a+1)(xu)\,dx|&\leq c\| |D_x|^{-\sigma}(xu)\|_2^2\,\| \widehat{(|D_x|^{a+1+2\sigma}(x\langle x\rangle_N^{2\theta})}\|_{L^1_{\xi}}\\
			&\leq c \|\langle x\rangle^{3/2}u\|_2
		\end{aligned}
	\end{equation*}
	where similar arguments in \eqref{G1}-\eqref{G4} show the validity of the previous inequality for $\theta<1$. 
	
	Next, we consider the third term in \eqref{EE1}. Using the identity \eqref{I2} one sees that
	\begin{equation}\label{EE4}
		\begin{split} 
			(a+1) \int_{\mathbb{R}} |D_x|^au\, xu\,x\langle x\rangle_N^{2\theta} dx
			=&(a+1) \int_{\mathbb{R}} |D_x|^a(xu)\, xu\,\langle x\rangle_N^{2\theta} dx\\
			&-(a+1)a\int_{\mathbb{R}} |D_x|^{a-1}\mathcal{H}(u)\, xu \, \langle x\rangle_N^{2\theta}\, dx\\ 
			=:&B_1+B_2.
		\end{split}
	\end{equation}
	First, we have that
	\begin{equation}\label{EE5}
		\begin{split} 
			B_1&=(a+1) \int_{\mathbb{R}} |D_x|^a(xu)\, xu\,\langle x\rangle_N^{2\theta} dx\\
			&=(a+1) \int_{\mathbb{R}} ||D_x|^{a/2}(xu)|^2\,\langle x\rangle_N^{2\theta} dx+(a+1)\int_{\mathbb{R}} |D_x|^{a/2}(xu) \big[|D_x|^{a/2};\langle x\rangle_N^{2\theta}\big](xu)\,dx\\
			&=(a+1) \int_{\mathbb{R}} ||D_x|^{a/2}(xu)|^2\,\langle x\rangle_N^{2\theta} dx+(a+1)\int_{\mathbb{R}} |D_x|^{a-1}(xu)\,|D_x|^{1-a/2}\,\big[|D_x|^{a/2};\langle x\rangle_N^{2\theta}\big](xu)\,dx\\
			&=:B_{11}+B_{12}.
		\end{split}
	\end{equation}
	By Lemma \ref{lowcommu}, we have
	\begin{equation*}
		|B_{12}|\leq \| |D_x|^{a-1}(xu)\|_2 \|\partial_x(\langle x\rangle^{2\theta})\|_{\infty} \|xu\|_2,
	\end{equation*}
	where by Proposition \ref{Dint} 
	\begin{equation}\label{EE5a} 
		\| |D_x|^{a-1}(xu)\|_2\leq c \| |D_x|^{a/2}(xu)\|_2^{2-2/a}\| xu\|_2^{2/a-1},
	\end{equation}
	which, by Young's inequality and the previous result $\alpha=3/2$ (see \eqref{Extraclass}) is bounded.

	It remains to consider the term $B_2$ in \eqref{EE4}
	\begin{equation}\label{EE7}
		\begin{split} 
			|B_2|&\leq c\int_{\mathbb{R}} |\langle x\rangle^{\gamma}|D_x|^{a-1}\mathcal H(u)\,xu \langle x\rangle^{2\theta-\gamma}|\,dx\\
			&\leq c\|\langle x\rangle^{\gamma}\mathcal H|D_x|^{a-1}(u)\|_2\|xu \langle x\rangle_N^{2\theta-\gamma}\|_2=:B_{21}\,B_{22}.
		\end{split}
	\end{equation}
	
	From \eqref{class3} one has that the term $B_{22}$ is bounded if $1+2\theta-\gamma\leq 3/2$, i.e., $2\theta\leq \gamma+1/2$. To control the factor $B_{21}$ in \eqref{EE7} we follow the argument in \cite{FONLP} (or the arguments in the proof of Theorem \ref{AA2}). Thus, we write
	\begin{equation} \label{EE8}
		\begin{split} 
			B_{21}&=\|\langle x\rangle^{\gamma}\mathcal H|D_x|^{a-1}u\|_2=c\| J^{\gamma}(|\xi|^{a-1}\,\operatorname{sgn}(\xi)\widehat{u})||_2\\
			&\leq c(\| J^{\gamma}(|\xi|^{a-1}\,\operatorname{sgn}(\xi)\chi(\xi)\widehat{u})||_2+\| J^{\gamma}(|\xi|^{a-1}\,\operatorname{sgn}(\xi)(1-\chi(\xi))\widehat{u})||_2)\\
			&=:B_{211}+B_{212},
		\end{split}
	\end{equation}
	where $\chi\in C^{\infty}_0(\mathbb{R})$ with $\chi(\xi)=1,\;|\xi|<1$ and $\chi(\xi)=0,\;|\xi|>2$. We observe, see  \cite{FONLP} that
	\begin{equation}\label{EE8a}
		|\xi|^{a-1}\,\operatorname{sgn}(\xi)\chi(\cdot)\in H^{(a-1/2)^{-}}(\mathbb{R})\;\;\;\;\;\;\text{and}\;\;\;\;\;\;\widehat{u}\in H^{3/2^-}(\mathbb{R}).
	\end{equation}
	Since $a-1/2>1/2$ it follows that $B_{211}<\infty$ if $\gamma <a-1/2$, which holds for $\gamma\leq 1/2$. Thus, we have bounded $B_{211}$. To control the contribution of $B_{212}$ one has that for $\gamma\leq 1$ (see \cite{FONLP})
	\begin{equation}\label{EE8b}
		\begin{aligned} 
			B_{212} &\leq c \|\langle x\rangle^{\gamma} J^{a-1}(u)\|_2\leq c(\|J^{a-1}u\|_2+\|xJ^{a-1}u\|_2)\\
			&\leq c( \|u\|_2 +\| |D_x|^{a/2}u\|_2+\|J^{a-1}(xu)\|_2)\\
			&\leq c( \|u\|_2 +\| |D_x|^{a/2}u\|_2+\|xu\|_2+\||D_x|^{a-1}(xu)\|_2) \\
			& \leq c( \|u\|_2 +\| |D_x|^{a/2}u\|_2+ \|xu\|_2+\| |D_x|^{a/2}(xu)\|_2^{2-2/a}\| xu\|_2^{2/a-1}),
		\end{aligned}
	\end{equation}
	which is bounded by \eqref{Extraclass}.
	
	\medskip

	Finally, after integration  in the time interval $[t_1,t_2]$ using the hypothesis and inserting \eqref{EE5}, \eqref{EE8a}, and \eqref{EE8b} in \eqref{EE2}, we get that
	\begin{equation}\label{EE9}
		\frac{a+1}{2}\int_{t_1}^{t_2}\int_{\mathbb{R}}  ||D_x|^{a/2}(xu)|^2\Big((x\langle x\rangle_N^{2\theta})'+2\langle x\rangle_N^{2\theta}\Big)dx\,dt<\infty.
	\end{equation}
	Since $\theta\in (0,1)$, we have
	\begin{equation*}
		\partial_x(x\langle x\rangle^{2\theta}_N)+2\langle x\rangle^{2\theta}_N\geq c_{\theta}
		\langle x\rangle^{2\theta},
	\end{equation*}
	we get from \eqref{EE9} and taking $N\to \infty$ that 
	\begin{equation}\label{xuestimatewithD}
		\langle x\rangle^{\theta}|D_x|^{a/2}(x u)\in L^{2}([t_1,t_2]:L^{2}(\mathbb{R})). 
	\end{equation}
	And reapplying the same ideas, one deduces
	\begin{equation}
		\langle x\rangle^{3/2+\theta} u\in L^{\infty}([t_1,t_2]:L^{2}(\mathbb{R})). 
	\end{equation}
	Consequently, by Proposition \ref{JandDprop} we deduce $x\langle x \rangle^{\theta}J^{a/2}u \in L^2([t_1,t_2]:L^2(\mathbb{R}))$, and since $J^{a/2}u\in C([-T,T]: L^2(\mathbb{R}))$, it follows
	\begin{equation}\label{conlcusioncase3/2}
		\langle x\rangle^{1+\theta}J^{a/2} u \in L^{2}([t_1,t_2]:L^{2}(\mathbb{R})). 
	\end{equation}
	Thus, for two different times $t_1\leq t_3<t_4\leq t_2$, $  \langle x\rangle^{1+\theta}J^{a/2} u(t_j)\in L^2(\mathbb{R})$, $j=3,4$. In view of \eqref{conlcusioncase3/2} and since $\alpha-1/2=1+\theta\in (1,3/2)$, the proof proceeds along the same lines as the previous case, with the regularity increased by $a/2$ through the application of $J^{a/2}$. The dispersive estimate remains unchanged, while the nonlinear term is controlled using Lemma \ref{Nonlinear estimates}. More precisely, in the nonlinear case, assume that the following estimate holds on some time interval \(I\):
	\begin{equation*}
		\langle x \rangle^{\alpha-k/2}J^{ak/2}u \in L^{2}(I:L^2(\mathbb{R})),    
	\end{equation*}
	for some $k\in \mathbb{N}_0$, $k \leq 2\alpha$. Then, in each step, the energy estimate is given by $\alpha-k/2=m+1/2+\theta$, $\theta\in(0,1)$, $m\in \mathbb{N}_0$. It follows from Lemma \ref{Nonlinear estimates} that
	\begin{equation*}
		\begin{aligned}
			\Big|\int_{\mathbb{R}}  J^{ak/2}\big(u\partial_x u\big) \big(x^{2m+1} \langle x\rangle^{2\theta}_N J^{ak/2}u\big)\, dx\Big|\leq & c\|J^{\max\{(3/2)^{+}; ak/2\}}u\|_{2}\|\langle x\rangle^{\alpha-k/2} J^{ak/2} u\|_{2}^2.
		\end{aligned} 
	\end{equation*}
	Since this estimate is independent of the step, it applies at every stage of the iterative scheme, and we do not need to repeat its proof.
	
	\medskip
	For the case $\alpha=2$, we can formally apply the energy estimate by multiplying the equation in \eqref{dgbo} by $x^4 u$ and proceed similarly to the previous case. This case can be made rigorously justified by using the commutator estimates found in D.~Li (Corollary 1.4) and an approximation argument.

	\medskip
	
	At this point, we emphasize that, thanks to the preceding results and the iterative structure of the energy estimate, we only need to study the linear term in the following energy estimates.
	
	\medskip
	
	\noindent \underline{\bf Case $\alpha\in (2,5/2]$}. By the previous case, we can assume
	\begin{equation}\label{class4}
		u\in C([-T,T]:H^{2a}(\mathbb{R}))\cap C([-T,T]:L^2(|x|^{4}\, dx)),    
	\end{equation}
	and at two different times $t_1,t_2\in [-T,T]$, $|x|^{\alpha}u(t_j)\in L^2(\mathbb{R})$. We first consider the case $\alpha\in (3/2,2)$ and use the same energy estimate as in the previous case. However, here we show that the above results extend to $\theta \in [1/2,1)$.

	We first observe that the arguments in \eqref{G1}-\eqref{G4} imply that for $0<\theta<1$, $\|\widehat{\partial_x^3(x\langle x\rangle_N^{2\theta})}\|_{L^1}<c$ with $c>0$ independent of $N$. Thus, the estimate \eqref{EE2} remains the same.  
	
	To estimate $B_{12}$ in \eqref{EE5}, we use the identity
	\begin{equation}\label{ExtraExpansion}
		[|D_x|^{a/2}; \langle x\rangle^{2\theta}_N]=\mathcal{S}_{0,1}(a/2)+\frac{a}{2}\mathcal{H}|D_x|^{a/2-1}(\cdot \, (\langle x\rangle_{N}^{2\theta})').    
	\end{equation}
	Then
	\begin{equation*}
		\begin{aligned}
			B_{12}=&(a+1)\int_{\mathbb{R}} |D_x|^{a/2}(xu)\,\mathcal{S}_{0,1}(a/2)(xu)\,dx\\
			&+\frac{(a+1)a}{2}\int_{\mathbb{R}} |D_x|^{a-1}(xu)\mathcal{H}\big((\langle x\rangle_N^{2\theta})'(xu)\big)\,dx\\
			=:& B_{121}+B_{122}.
		\end{aligned}    
	\end{equation*}
	By Theorem \ref{gv2}
	\begin{equation*}
		\begin{aligned}
			|B_{121}|=&\Big|\int_{\mathbb{R}}  |D_x|^{a/2}(x u)\mathcal{S}_{0,1}(a/2)(x u)\,dx\Big|\\
			\leq & \||D_x|^{-\sigma_1+a/2}(x u)\|_2\, \||D_x|^{-\sigma_2}(x u)\|_2 \,\||D_x|^{\sigma_1}\mathcal{S}_{0,1}(a/2)\,|D_x|^{\sigma_2}\|_{L^{2}\to L^{2}},
		\end{aligned}
	\end{equation*}
	where we set $\sigma_1+\sigma_2+a/2=2$, $\sigma_1>a/2$. It follows from the ideas in \eqref{energyInequalityM} that
	\begin{equation*}
		\||D_x|^{-\sigma_1+a/2}(x u)\|_2, \, \||D_x|^{-\sigma_2}(x u)\|_2\leq c \| \langle x\rangle^{3/2} u\|_2<\infty,
	\end{equation*}
	and if $\theta <1$, by previous arguments
	\begin{equation*}
		\||D_x|^{\sigma_1}\mathcal{S}_{0,1}(a/2)\,|D_x|^{\sigma_2}\|_{L^{2}\to L^{2}}\leq c \| \widehat{\partial_x^2(\langle x\rangle^{2\theta}_N)}\|_1<\infty.
	\end{equation*}
	Next, since $|(\langle x\rangle_N^{2\theta})'x|\leq c \langle x\rangle^{2\theta}\leq \langle x\rangle^2$, with constant $c>0$ independent of $N$, we have
	\begin{equation*}
		\begin{aligned}
			|B_{122}|\leq c \||D_{x}|^{a-1}(xu)\|_2\|\langle x \rangle^2 u\|_2 \leq c  \||D_{x}|^{a/2}(xu)\|_2\|\langle x \rangle^2 u\|_2, 
		\end{aligned}    
	\end{equation*}
	which is controlled by previous cases and \eqref{class4}.
	
	By the arguments in \eqref{EE7}-\eqref{EE8b}, we have
	\begin{equation}
		\begin{split} 
			|B_2|&\leq c\|\langle x\rangle^{\gamma}\mathcal H|D_x|^{a-1}(u)\|_2\|xu \langle x\rangle_N^{2\theta-\gamma}\|_2=B_{21}\,B_{22},
		\end{split}
	\end{equation}
	where $B_{21}<\infty$ if $\gamma<a-1/2$, and by \eqref{class4}, $B_{22}$ is controlled if $1+2\theta-\gamma\leq 2$. Setting $\gamma=(a-1/2)^{-}$, we have that $B_{22}<\infty$ if $\theta<1/2+(a-1/2)/2=(a+1/2)/2$.

	Consequently, if $1/2\leq \theta<\min\{1,\, (a+1/2)/2\}$, we have
	\begin{equation}\label{NewconA}
		\langle x\rangle^{\theta}\,|D_x|^{a/2}(x u)\in L^2([t_1,t_2]:L^2(\mathbb{R}))\hskip10pt\text{and}\hskip10pt |x|^{3/2+\theta}u \in L^{\infty}([t_1,t_2]:L^2(\mathbb{R})).
	\end{equation}\
	With the above result, we can repeat the argument and progressively increase the weight $\theta<1$. The only modification is the estimate for $B_2$, for which taking $\gamma=(a-1/2)^{-}$ and  using \eqref{NewconA}, we have that $B_{22}<\infty$ if $\theta<1/2+(a-1/2)(1/2+1/4)$. Thus, we extend \eqref{NewconA} to $\theta<\min\{1,\, 1/2+(a-1/2)(1/2+1/4)\}$. Iterating this argument, we extend \eqref{NewconA} to weights $\theta<1/2+(a-1/2)\Big(\sum_{j=1}^M \frac{1}{2^j}\Big)$, $\theta<1$. Taking the number of iteration $M\to \infty$, we deduce \eqref{NewconA} to $\theta<1$.

	Finally, by Proposition \ref{JandDprop} we arrive at
	\begin{equation}\label{weightJ3/2a}
		\langle x \rangle^{1+\theta}J^{a/2}u \in L^2([t_1,t_2]:L^2(\mathbb{R})). 
	\end{equation}
	for $\theta\in [1/2,1)$. Since $1+\theta\in [3/2,2)$, we can reapply the same argument in the previous case, completing the proof.

	\medskip
	
	This completes the case $\alpha\in(2,5/2)$. As in the case $\alpha=3/2$, the endpoint $\alpha=5/2$ is treated with a polynomial weight, now $x^{5}$.
	
	Next, we consider the case $\alpha=5/2$. Thus, formally we multiply the equation in \eqref{dgbo} by $x^5 u$ and integrate the result in the space variable. From the hypothesis and previous result, we just need to consider the contribution of the linear terms
	\begin{equation*}
		A_{0} = -\int_{\mathbb{R}} \partial_{x} |D_{x}|^{a} u \cdot x^{5} u \, dx.
	\end{equation*}
	We get that
	\begin{equation}\label{case5a}
		\begin{split}
			A_{0} &= -\int_{\mathbb{R}} \partial_{x} |D_{x}|^{a}(x^{2}u) x^{2}u \cdot x \, dx + 2(a+1)\int_{\mathbb{R}} |D_{x}|^{a}(xu) x\, x^{2}u \, dx\\
			&\hskip12pt + a(a+1)\int_{\mathbb{R}} |D_{x}|^{a-1}\mathcal{H}u \,x\, x^{2}u \, dx =: A_{1} + A_{2} + A_{3}.
		\end{split}
	\end{equation}
	Thus, combining integration by parts and the identity \eqref{I2} we have:
	\begin{equation*}
		\begin{split}
			A_{1} &= \int_{\mathbb{R}} |D_{x}|^{a}(x^{2}u) \partial_{x}(x^{2}u) x \, dx + \int_{\mathbb{R}} |D_{x}|^{a}(x^{2}u) x^{2}u \, dx \\
			&= \int_{\mathbb{R}} (|D_{x}|^{a/2}(x^{2}u))^{2} dx - A_{1} + a \int_{\mathbb{R}} x^{2}u |D_{x}|^{a}(x^{2}u) \, dx
		\end{split}
	\end{equation*}
	Therefore,
	\begin{equation*}
		A_{1} = \frac{a+1}{2} \int_{\mathbb{R}} (|D_{x}|^{a/2}(x^{2}u))^{2} dx.
	\end{equation*}
	Also, by \eqref{I2}
	\begin{equation*}
		\begin{aligned}
			A_{2} &= 2(a+1)\int_{\mathbb{R}} |D_{x}|^{a}(x^{2}u)\,x^{2}u \, dx - 2a(a+1)\int_{\mathbb{R}} |D_{x}|^{a-1}\mathcal{H}(xu) x^{2}u \, dx\\
			&=:2(a+1)\int_{\mathbb{R}} (|D_{x}|^{a/2}(x^{2}u))^{2} dx+A_{21}
		\end{aligned}
	\end{equation*}
	with $A_{21} $ satisfying 
	\begin{equation*}
		|A_{21}| \le 2a(a+1) \| |D_{x}|^{a-1}(xu) \|_{2} \|x^{2}u\|_{2} \le c \||D_x|^{1/2}(xu)\|_{2} \|x^{2}u\|_{2}<\infty,
	\end{equation*}
	from previous results. Next, we have that 
	
	\begin{equation*}
		|A_3| \le c \| |\langle x\rangle^{1/2}|D_{x}|^{a-1}\mathcal Hu \|_{2} \|\langle x\rangle^{2+1/2}u\|_{2} =:A_{31}\,A_{32}.
	\end{equation*}
	Since  $A_{32} = \| \langle x \rangle^{2+1/2} u \|_{2}$ is bounded by previous argument,  it only remains to consider $A_{31}$. An already applied argument shows that
	\begin{equation*}
		\begin{split}
			A_{31} &\le c \| J^{1/2+} ( |\xi|^{a-1} \operatorname{sgn}(\xi) \widehat{u}(\xi) ) \|_{2} \\
			&\le c \| J^{1/2+} ( |\xi|^{a-1} \operatorname{sgn}(\xi) \chi(\xi) \widehat{u}(\xi) ) \|_{2} + c \| \langle x \rangle^{1/2} J^{a-1} u \|_{2} \\
			&=: A_{311} + A_{312}.
		\end{split}
	\end{equation*}
	Since $|\xi|^{a-1} \operatorname{sgn}(\xi) \chi(\xi) \in H^{a-1/2}(\mathbb{R})$ (with $a-1/2 > 1/2$) and $\widehat{u} \in H^{2}(\mathbb{R})$, we have that 
	$A_{311} < \infty$. 
	
	Finally,
	\begin{equation*}
		\begin{split}
			A_{312} &= \| \langle x \rangle^{1/2} J^{a-1} u \|_{2} \le c(\|J^{a-1}u\|_{2} + \|x J^{a-1}u\|_{2}) \\
			&\le c(\|u\|_{H^{a/2}} + \|J^{a-1}(xu)\|_{2}) \\
			&\le c(\|u\|_{H^{a/2}} + \|xu\|_{2} + \| |D_{x}|^{a/2}(xu) \|_{2}) < \infty
		\end{split}
	\end{equation*}
	by previous results. As before, we remark that the above formal proof can be made rigorously justified by an approximation argument.
	
	\medskip

	For weights larger than $5/2$, the function $x^{2}u$ plays the role of $xu$, and we use the identity \eqref{I3} below, which follows from \eqref{generalconmmutator} with $m=2$. In this range the terms in \eqref{I3} whose symbols are less regular at the origin become relevant, and they are responsible for the restriction $\alpha<a+3/2$ in the absence of the zero mean condition.
	
	\noindent \underline{\bf Case $\alpha\in(5/2,\nu)$ with $\nu=\min\{a+3/2,3\}$}.  We set $\alpha=5/2+\theta$, where $0<\theta<a-1$, if $1<a\leq 3/2$, and $0<\theta<1/2$, if $3/2<a<2$. By the previous case, we may assume
	\begin{equation}\label{class5}
		u\in C([-T,T]:H^{5a/2}(\mathbb{R}))\cap C([-T,T]:L^2(|x|^{5}\, dx)),    
	\end{equation}
	and at two different times $t_1,t_2\in [-T,T]$, $|x|^{\alpha}u(t_j)\in L^2(\mathbb{R})$. Furthermore, it holds
	\begin{equation*}
		|D_x|^{a/2}(x^2 u)\in L^{\infty}([t_1,t_2]:L^2(\mathbb{R})).   
	\end{equation*}
	For the moment, we assume that $0<\theta<(a-1)/2$.  We multiply the equation in \eqref{dgbo} by $x^4u\langle x\rangle^{2\theta}_Nx\;$ and integrate in the space variable. As before, we are reduced to estimating the term coming from the dispersion in \eqref{dgbo}, namely
	\begin{equation*}
		A_0= -\int_{\mathbb{R}} x^2\,|D_x|^a\partial_x u\,x^2u\,\langle x\rangle^{2\theta}_N x\,dx.
	\end{equation*}
	By Proposition \ref{GeneralConmutator}, we shall use the identity
	\begin{equation}\label{I3}
		x^{2}\partial_{x}|D_{x}|^{a}u = \partial_{x}|D_{x}|^{a}(x^{2}u) - 2(a+1)|D_{x}|^{a}(xu) + a(a+1)|D_{x}|^{a-1}\mathcal{H}u.
	\end{equation}
	One has that
	\begin{equation}\label{NA1}
		\begin{split}
			A_0&=-\int_{\mathbb{R}} |D_x|^a\partial_x (x^2u)\,x^2u\,\langle x\rangle^{2\theta}_{N}x\,dx+2(a+1)\int_{\mathbb{R}} |D_x|^a(xu)\,x^2u\,\langle x\rangle^{2\theta}_{N}x\,dx\\
			&\hskip12pt-(a+1)a\int_{\mathbb{R}} |D_x|^{a-1}\mathcal H(u)\,x^2u\,\langle x\rangle^{2\theta}_{N}x\,dx=:A_1+A_2+A_3.
		\end{split}
	\end{equation}
	To handle $A_1$, we use \eqref{eq1abc} to get
	\begin{equation}\label{NA_1}
		\begin{split}
			A_1 &= \frac{a+1}{2} \int_{\mathbb{R}} \left( (|D_x|^{a/2} (x^2 u))^2 \right) (\langle x\rangle_{N}^{2\theta}x)' \, dx \\
			&\hskip12pt + \frac{1}{2} \int_{\mathbb{R}} (x^2 u) \tilde{R}_0(a+1)(x^2 u) \, dx,
		\end{split}
	\end{equation}
	where setting $\sigma=1+a/2$, by \eqref{energyInequalityM} we get
	\begin{equation*}
		\begin{aligned}
			|A_{12}| &\le c \| |D_x|^{-\sigma} (x^2 u) \|_2^2 \,\ \| \widehat{\partial^3_x (\langle x\rangle_{N}^{2\theta}x)} \|_{L^1_{\xi}}\\
			&\le c \| \langle x \rangle^{5/2-} u \|_2^2 \| \widehat{\partial^3_x (\langle x\rangle_{N}^{2\theta}x)} \|_{L^1_{\xi}}.
		\end{aligned}
	\end{equation*}
	The above is bounded by the arguments in \eqref{G1}-\eqref{G4}. To estimate $A_2$ in \eqref{NA1}, we use the identity \eqref{I2}  to get
	\begin{equation}
		\begin{aligned}
			\label{NA_2}
			A_2 &= 2(a+1) \int_{\mathbb{R}} |D_x|^a (x^2 u) x^2 u\langle x\rangle_N^{2\theta} \, dx - 2(a+1)a \int_{\mathbb{R}} |D_x|^{a-1} \mathcal{H}(xu) x^2 u\langle x\rangle_N^{2\theta} \, dx\\
			&=: A_{21} + A_{22}.
		\end{aligned}
	\end{equation}
	As in previous cases, we write
	\begin{equation}\label{NA_21}
		\begin{aligned}
			A_{21} &= 2(a+1) \int_{\mathbb{R}} (|D_x|^{a/2} (x^2 u))^2 \langle x\rangle^{2\theta}_N \, dx + 2(a+1) \int_{\mathbb{R}} |D_x|^{a/2} (x^2 u)[|D_x|^{a/2}; \langle x\rangle^{2\theta}_N]  (x^2 u) \, dx \\
			&=: A_{211} + A_{212}. 
		\end{aligned}
	\end{equation}
	To handle $A_{212}$, we write
	\begin{equation*}\label{NA212}
		\begin{aligned}
			\int_{\mathbb{R}} |D_x|^{a/2} (x^2 u)[|D_x|^{a/2}; \langle x\rangle^{2\theta}_N]  (x^2 u) \, dx=\int_{\mathbb{R}} |D_x|^{a-1} (x^2 u)|D_x|^{1-a/2}[|D_x|^{a/2}; \langle x\rangle^{2\theta}_N]  (x^2 u) \, dx.  
		\end{aligned}    
	\end{equation*}
	By Lemma \ref{lowcommu}, we have
	\begin{equation*}
		\begin{aligned}
			|A_{212}|\leq &c  \|\partial_x(\langle x\rangle^{2\theta}_N)\|_{\infty}\||D_x|^{a-1}(x^2 u)\|_2\|x^2 u\|_2^2\leq c \||D_x|^{a/2}(x^2 u)\|_2\|x^2 u\|_2^2.  
		\end{aligned}    
	\end{equation*}
	This completes the bound for $A_{21}$. It remains to handle $A_{22}$ in \eqref{NA_2} and $A_3$  in \eqref{NA1}. We first observe that from the arguments in \eqref{EE7}-\eqref{EE8b}, for $0\leq \gamma<a-1/2$ (in particular, $\gamma\leq 1/2$ is admissible) it follows
	\begin{equation}\label{Hilbertineqconclu}
		\|\langle x \rangle^{\gamma}|D_x|^{a-1}\mathcal{H}f\|_2\leq c\big(\|f\|_2+\||D_x|^{a/2}f\|_2+\|xf\|_2+\||D_x|^{a/2}(xf)\|_2^{2-2/a}\|xf\|_2^{2/a-1}\big).   
	\end{equation}
	Hence, for any $\max\{2\theta-1/2,0\}\leq \gamma<a-1/2$, we get
	\begin{equation}\label{NA_4}
		\begin{aligned}
			|A_{22}|&\leq c\| \langle x\rangle^{\gamma}|D_x|^{a-1}\mathcal{H} (xu)\|_2\,\| \langle x\rangle^{2(1+\theta)-\gamma}u\|_2\\
			&\leq c \big(\|\langle x\rangle^2 u\|_2+\||D_x|^{a/2}(xu)\|_2+\||D_x|^{a/2}(x^2u)\|_2^{2-2/a}\|\langle x\rangle^2 u\|_2^{2/a-1}\big)\| \langle x\rangle^{5/2}u\|_2,
		\end{aligned}
	\end{equation}
	where the condition $2\theta-1/2\leq \gamma$ implies that $2(1+\theta)-\gamma\leq 5/2$, and it follows that $2\theta-1/2<a-1/2$, whenever $0<\theta<1/2$.  This completes the estimate for $A_2$.

	We observe that the results established so far for $A_1$ and $A_2$ hold for every $0<\theta<1/2$, $1<a<2$.
	
	Next, we consider the term $A_3$ in \eqref{NA1}. We follow the previous argument for $A_{22}$, thus we have
	\begin{equation}
		\begin{aligned}
			\label{NA_7}
			|A_{3}|&\leq c\| \langle x\rangle^{\gamma}|D_x|^{a-1}\mathcal{H} u\|_2\,\| \langle x\rangle^{2\theta+3-\gamma}u\|_2=:A_{31}\cdot A_{32},
		\end{aligned}
	\end{equation}
	where  $A_{31}<\infty$ if $\gamma<a-1/2$ and $A_{32}<\infty$ if $\theta\leq (\gamma-1/2)/2$ (this last condition is due to \eqref{class5}). Combining these conditions, we get that \eqref{NA_7} is controlled if $\theta<(a-1)/2$. We observe that this is the only step so far where we require this condition on $\theta$.
	\medskip
	
	Finally, we gather the first terms in \eqref{NA_1} and \eqref{NA_21} to get
	\begin{equation}\label{EE9a}
		\frac{a+1}{2}\int_{\mathbb{R}} ||D_x|^{a/2}(x^2u)|^2\Big((x\langle x\rangle^{2\theta}_N)'+2\langle x\rangle^{2\theta}_N\Big)dx.
	\end{equation}
	When $\theta\in (0,1/2)$, there exists $c_{\theta}>0$ such that
	\begin{equation*}
		\partial_x(x\langle x\rangle^{2\theta}_N)+2\langle x\rangle^{2\theta}_N\geq c_{\theta} \langle x\rangle^{2\theta}_N.   
	\end{equation*}
	Collecting all the results above, we conclude that for $\alpha=5/2+\theta$, $0<\theta<(a-1)/2$, one has that 
	\begin{equation}\label{A_8}
		\langle x\rangle^{\theta}\,|D_x|^{a/2}(x^2u)\in L^2([t_1,t_2]:L^2(\mathbb{R}))\hskip10pt\text{and}\hskip10pt |x|^{5/2+\theta}u \in L^{\infty}([t_1,t_2]:L^2(\mathbb{R})).
	\end{equation}
	With this result at hand, we can iterate the argument to increase the value of $\theta$. The only modification occurs in the estimate for $A_3$ in \eqref{NA_7}. For instance, if \eqref{A_8} holds, then $A_{32}$ is controlled provided $2\theta+3-\gamma<5/2+(a-1)/2$ with $\gamma<a-1/2$. Thus, by taking $\gamma=(a-1/2)^{-}$, \eqref{A_8} implies that $A_{3,2}<\infty$, if $(a-1)/2\leq \theta< 3(a-1)/4$ with $\theta<1/2$. As a consequence, \eqref{A_8} is now valid for $0<\theta<\min\{3(a-1)/4,\, 1/2\}$. Iterating this argument, we obtain that the recursive scheme is given by $\theta<\min\{(a-1)\Big(\sum_{j=1}^M \frac{1}{2^j}\Big);\, 1/2\}$, and, upon letting the number of iterations $M\to \infty$, we extend \eqref{A_8} to the desired range $0<\theta<\min\{a-1,1/2\}$.

	By \eqref{A_8} and Proposition \ref{JandDprop} we conclude
	\begin{equation}\label{Conclusioncasea3/2}
		\langle x\rangle^{2+\theta}\,J^{a/2}u \in L^2([t_1,t_2]:L^2(\mathbb{R})).    
	\end{equation}
	Since $2+\theta\in (0,\nu_1)$ with $\nu_1=\min\{a+1,5/2\}$, we can apply the estimates in the above case with weights $\alpha\in (2, 5/2)$ to obtain the claimed result.

	\medskip 
	
	Once the weight reaches the threshold $3/2+a$, part (ii) of Theorem \ref{AA2} shows that the solution must have zero mean, and this cancellation allows us to continue the iteration.

	\noindent \underline{\bf Case $\alpha\in [\nu,7/2]$ with $\nu=\min\{3/2+a,3\}$}. Let us first work on the case $\alpha\in [\nu,7/2)$. If $a\in (3/2,2)$, by the above case, we can assume that
	\begin{equation}\label{class5E}
		u\in C([-T,T]:H^{3^{-}}(\mathbb{R}))\cap C([-T,T]:L^2(|x|^{6^{-}}\, dx)),    
	\end{equation}
	and at two different times $t_1,t_2\in [-T,T]$, $|x|^{\alpha}u(t_j)\in L^2(\mathbb{R})$.  
	
	\noindent  If $a\in (1,3/2]$, we assume that
	\begin{equation}\label{class5.1}
		u\in C([-T,T]:H^{a(3/2+a)^{-}}(\mathbb{R}))\cap C([-T,T]:L^2(|x|^{2(3/2+a)^{-}}\, dx)),    
	\end{equation}
	and at two different times $t_1,t_2\in [-T,T]$, $|x|^{\alpha}u(t_j)\in L^2(\mathbb{R})$. Moreover, for $a\in (1,3/2]$, by \eqref{class5.1} and (ii) in Theorem \ref{AA2}, we have that $\widehat{u_0}(0)=0$, and by the $L^1$-conservation law of \eqref{dgbo} it must follow
	\begin{equation}\label{zeromean}
		\widehat{u}(0,t)=\int_{\mathbb{R}} u(x,t)\, dx=0, \qquad t\in [-T,T].    
	\end{equation}
	We set $\alpha=5/2+\theta$, and use the same energy estimate as in the previous case, i.e., with the factor $x^5\langle x\rangle_{N}^{2\theta}$ applied to the equation in \eqref{dgbo}. The difference is that we increase the values of $\theta$ to the following intervals:
	\begin{equation*}
		\left\{\begin{aligned}
			&a-1\leq \theta < 1, \qquad \quad \text{if $1<a\leq 3/2$,  and $\widehat{u}(0)=0$}, \\
			&1/2\leq \theta < a-1, \qquad \text{if $3/2<a<2$, no condition on $\widehat{u}(0)$}, \\
			&a-1\leq \theta < 1, \qquad \quad \text{if $3/2<a< 2$, and $\widehat{u}(0)=0$}.
		\end{aligned}\right.    
	\end{equation*}
	We need to study the factors $A_1$, $A_2$ and $A_3$ in \eqref{NA1}. The estimate for $A_1$ is unchanged and remains valid for all $\theta<1$. To control \eqref{NA_21} for $0<\theta<1$, by Theorem \ref{gv2} one has that 
	\begin{equation*}
		A_{21}=2(a+1)\int_{\mathbb{R}} (|D_x|^{a}(x^2u))^2\,\langle x\rangle^{2\theta}_{N}\,dx+2(a+1)\int_{\mathbb{R}} |D_x|^{a/2}(x^2u)\big[|D_x|^{a/2};\langle x \rangle^{2\theta}_N\big](x^2u)\,dx.
	\end{equation*}
	Using \eqref{ExtraExpansion}, we have
	\begin{equation*}
		\begin{aligned}
			2(a+1)\int_{\mathbb{R}} |D_x|^{a/2}(x^2u)\big[|D_x|^{a/2};\langle x \rangle^{2\theta}_N\big](x^2u)\,dx=&2(a+1)\int_{\mathbb{R}} |D_x|^{a/2}(x^2u)\mathcal{S}_{0,1}(a/2)(x^2u)\,dx\\
			&+(a+1)a \int_{\mathbb{R}} |D_x|^{a-1}(x^2u)\mathcal{H}\big((\langle x \rangle^{2\theta}_N)'x^2u\big)\,dx      
		\end{aligned}
	\end{equation*}
	with
	\begin{equation*}
		\begin{aligned}
			\Big|\int_{\mathbb{R}} & |D_x|^{a/2}(x^2u)\mathcal{S}_{0,1}(a/2)(x^2u)\,dx\Big|\\
			&\leq \||D_x|^{-\sigma_1+a/2}(x^2u)\|_2\, \||D_x|^{-\sigma_2}(x^2u)\|_2 \,\||D_x|^{\sigma_1} \mathcal{S}_{0,1}(a/2)\,|D_x|^{\sigma_2}\|_{L^{2}\to L^{2}},
		\end{aligned}
	\end{equation*}
	where we set $\sigma_1+\sigma_2+a/2=2$, $\sigma_1>a/2$. By \eqref{energyInequalityM}. we have
	\begin{equation*}
		\||D_x|^{-\sigma_1+a/2}(x^2u)\|_2\, \||D_x|^{-\sigma_2}(x^2u)\|_2\leq c \| \langle x\rangle^{5/2} u\|_2<\infty,
	\end{equation*}
	and if $\theta <1$, by previous argument
	\begin{equation*}
		\||D_x|^{\sigma_1}\mathcal{S}_{0,1}(a/2)\,|D_x|^{\sigma_2}\|_{L^{2}\to L^{2}}\leq c \| \widehat{\partial_x^2(\langle x\rangle^{2\theta}_N)}\|_1<\infty.
	\end{equation*}
	On the other hand, we have
	\begin{equation}\label{EstimateA21}
		\begin{aligned}
			\Big|\int_{\mathbb{R}} |D_x|^{a-1}(x^2u)\mathcal{H}\big((\langle x \rangle^{2\theta}_N)'x^2u\big)\,dx \Big|\leq & c\||D_x|^{a-1}(x^2u)\|_2\|\langle x\rangle^{1+2\theta} u\|_2\\
			\leq & c \||D_x|^{a/2}(x^2u)\|_2\|\langle x\rangle^{1+2\theta} u\|_2,
		\end{aligned}    
	\end{equation}
	which by \eqref{class5E} and \eqref{class5.1} is bounded if $1+2\theta<\nu$, $\nu=\min\{3/2+a,3\}$. In particular, if $3/2<a<2$, $\theta<1$.

	Now, to control $A_{22}$ in \eqref{NA_2} we write
	\begin{equation}
		\begin{aligned}
			|A_{22}|&\leq c\|\langle x\rangle^{\gamma}|D_x|^{a-1}\mathcal H(xu)\|_2 \|\langle x \rangle^{2\theta-\gamma}x^2u\|_2=:A_{221}A_{222}.
		\end{aligned}
	\end{equation}
	If $\gamma<a-1/2$, by \eqref{Hilbertineqconclu} and \eqref{A_8}, it follows that $A_{221}<\infty$. While by \eqref{class5} and \eqref{class5.1}, we have that $A_{222}<\infty$ if $2\theta-\gamma+2<\nu$, where $\nu=\min\{3/2+a,\, 3\}$. This completes the study of $A_2$.

	Next, we consider $A_3$ in \eqref{NA1}. Thus,
	\begin{equation*}
		|A_3|\leq c\|\,\langle x\rangle^{\gamma_1}\,|D_x|^{a-1}\mathcal H u\|_2\, \|\langle x \rangle^{2\theta+2-\gamma_1} u\|_2=:A_{311}A_{312}.
	\end{equation*}
	By familiar arguments one has $A_{311}<\infty$, if $\gamma_1<a-1/2$, and by  Lemma \ref{Hardy-type} if $\widehat{u}(0)=0$, this extends to $\gamma_1<a+1/2$. By \eqref{class5} and \eqref{class5.1},  $A_{312}<\infty$, if $2\theta-\gamma_1+2<\nu$, where $\nu=\min\{3/2+a,\, 3\}$. This concludes the estimate for $A_3$.
	
	Summarizing, in the case $1<a\leq 3/2$, taking $\gamma=\gamma_1=(a-1/2)^{-}$ and $a-1\leq\theta<a-1/2$, we deduce 
	\begin{equation}\label{NewA_8}
		\langle x\rangle^{\theta}\,|D_x|^{a/2}(x^2u)\in L^2([t_1,t_2]:L^2(\mathbb{R}))\hskip10pt\text{and}\hskip10pt |x|^{5/2+\theta}u \in L^{\infty}([t_1,t_2]:L^2(\mathbb{R})).
	\end{equation}
	The key point is that thanks to \eqref{NewA_8}, we can now increase the value of $\theta$ in the previous estimates for \eqref{EstimateA21}, $A_{222}$ and $A_{312}$, and hence, by iteration, extend \eqref{NewA_8} to the full desired range $a-1\leq \theta <1$. We remark that due to the zero-mean condition, the estimate for $A_{311}$ remains valid over a wider range of parameters $\gamma_1$.

	Similarly, if $3/2<a<2$, we can iterate as in the previous case and obtain that \eqref{NewA_8} holds for $1/2\leq \theta<a-1$, without any restriction on $\widehat{u}(0)$. Once \eqref{NewA_8} holds for $\theta<a-1$, i.e., $\alpha=5/2+\theta<a+3/2$, we can apply (ii) in Theorem \ref{AA2} to obtain that the zero-mean condition \eqref{zeromean} must hold. Using this and iterations, we extend \eqref{NewA_8} to $a-1\leq \theta<1$. Once again, $\widehat{u}(0)=0$ allows the estimate $A_{3111}$ to be extended beyond its original range of validity.
	
	By Proposition \ref{JandDprop} and \eqref{NewA_8}, we get
	\begin{equation*}
		\langle x \rangle^{2+\theta}J^{a/2}u\in L^{2}([t_1,t_2]:L^2(\mathbb{R})),   
	\end{equation*}
	$0<\theta<1/2$. As noted above, this estimate allows us to iterate the argument from the previous case obtaining the desired result.

	\medskip
	
	As before, the endpoint $\alpha=7/2$ is treated separately, now with the weight $x^{7}$ and the identity \eqref{id3}, which follows from \eqref{generalconmmutator} with $m=3$.
	
	Next, we consider the case $\alpha=7/2$. By the previous result, we can assume that
	\begin{equation}\label{class6}
		u\in C([-T,T]:H^{(7a/2)^{-}}(\mathbb{R}))\cap C([-T,T]:L^2(|x|^{7^{-}}\, dx)),
	\end{equation}
	\eqref{zeromean} holds, and the first result in \eqref{NewA_8} is valid for $\theta<1$. Thus, formally we multiply the equation in \eqref{dgbo} by $x^7 u$ and integrate the result in the space variable. From the hypothesis and previous results, we just need to consider the contribution of the linear terms:
	\begin{equation*}
		A_{0} = -\int_{\mathbb{R}} \partial_{x} |D_{x}|^{a} u \cdot x^{7} u \, dx
	\end{equation*}
	Using the identity:
	\begin{equation}
		\label{id3}
		\begin{aligned}
			x^{3}\partial_{x}|D_{x}|^{a}u &= \partial_{x}|D_{x}|^{a}(x^{3}u) +3(a+1)|D_{x}|^{a}(x^2u) \\
			&\hskip12pt+3a(a+1)|D_{x}|^{a-1}\mathcal{H}(xu)-a(a+1)(a-1)|D_x|^{a-2}u,
		\end{aligned}
	\end{equation}
	we get that:
	\begin{equation}\label{a71}
		\begin{split}
			A_{0} &= -\int_{\mathbb{R}} \partial_{x} |D_{x}|^{a}(x^{3}u) x^{3}u \cdot x \, dx + 3(a+1)\int_{\mathbb{R}} |D_{x}|^{a}(x^2u) x\, x^{3}u \, dx \\
			&\hskip12pt +3 a(a+1)\int_{\mathbb{R}} |D_{x}|^{a-1}\mathcal{H}(xu) \,x\, x^{3}u \, dx -a(a^2-1)\int_{\mathbb{R}} |D_x|^{a-2}u \,x\, x^{3}u \, dx\\
			&=: A_{1} + A_{2} + A_{3}+A_{4}.
		\end{split}
	\end{equation}
	Thus, combining integration by parts and the identity \eqref{I2} we have:
	\begin{equation*}
		\begin{aligned}
			A_{1} &= \int_{\mathbb{R}} |D_{x}|^{a}(x^{3}u) \partial_{x}(x^{3}u) x \, dx + \int_{\mathbb{R}} |D_{x}|^{a}(x^{3}u) x^{3}u \, dx \\
			&= \int_{\mathbb{R}} (|D_{x}|^{a/2}(x^{3}u))^{2} dx - A_{1} + a \int_{\mathbb{R}} x^{3}u |D_{x}|^{a}(x^{3}u) \, dx.
		\end{aligned}
	\end{equation*}
	Therefore,
	\begin{equation*}
		A_{1} = \frac{a+1}{2} \int_{\mathbb{R}} (|D_{x}|^{a/2}(x^{3}u))^{2} dx.
	\end{equation*}
	Also,
	\begin{equation*}
		\begin{aligned}
			A_{2} &= 3(a+1)\int_{\mathbb{R}} |D_{x}|^{a}(x^{3}u)\,x^{3}u \, dx - 3a(a+1)\int_{\mathbb{R}} |D_{x}|^{a-1}\mathcal{H}(x^2u) x^{3}u \, dx\\
			&=3(a+1)\int_{\mathbb{R}} (|D_{x}|^{a/2}(x^{3}u))^{2} dx+A_{21},
		\end{aligned}
	\end{equation*}
	with  
	\begin{equation*}
		|A_{21}| \le 3a(a+1) \| |D_{x}|^{a-1}(x^2u) \|_{2} \|x^{3}u\|_{2} \le c \||D_x|^{a/2}(x^2u)\|_{2} \|x^{3}u\|_{2}<\infty,
	\end{equation*}
	by previous results. Next, we have 
	\begin{equation*}
		|A_3|\leq c\| \langle x\rangle^{1/2+}|D_x|^{a-1}\,\mathcal H(xu)\|_2\|\langle x\rangle^{1/2-}x^3u\|_2=:A_{31} A_{32},
	\end{equation*}
	with $A_{32}<\infty$ from previous results. It only remains to consider $A_{31}$. A familiar argument shows that
	\begin{equation*}
		\begin{aligned}
			A_{31} &\le c \| J^{1/2+} ( |\xi|^{a-1} \operatorname{sgn}(\xi) \widehat{xu}(\xi) ) \|_{2}+c \| \langle x \rangle^{1/2} J^{a-1} (xu) \|_{2}  \\
			&=:A_{311} + A_{312}.
		\end{aligned}
	\end{equation*}
	Since $|\xi|^{a-1} \operatorname{sgn}(\xi) \chi(\xi) \in H^{a-1/2}(\mathbb{R})$ (with $a-1/2 > 1/2$) and $\widehat{xu} \in H^{2}(\mathbb{R})$, we have that $A_{311} < \infty$. 
	
	Also,
	\begin{equation*}
		\begin{aligned}
			A_{312} &= \| \langle x \rangle^{1/2} J^{a-1} (xu) \|_{2} \le c(\|J^{a-1}(xu)\|_{2} + \|x J^{a-1}(xu)\|_{2}) \\
			&\le c(\|xu\|_{H^{a/2}} + \|J^{a-1}(x^2u)\|_{2}) \\
			&\le c(\|xu\|_{H^{a/2}} + \|x^2u\|_{2} + \| |D_{x}|^{a/2}(x^2u) \|_{2}) < \infty
		\end{aligned}
	\end{equation*}
	by previous results.

	Finally, we consider the term $A_4$. Then,
	\begin{equation*}
		|A_4|\leq c \| \langle x\rangle^{1/2+}|D_x|^{a-2}u\|_2\|\langle x\rangle^{1/2-}x^3u\|_2=:A_{41} A_{42}
	\end{equation*}
	where $A_{42}<\infty$ from previous results. To bound $A_{41}$ we write
	\begin{equation*}
		A_{41}\leq c(\| J^{1/2+}(|\xi|^{a-2}\,\chi(\xi)\,\widehat{u}(\xi))\|_2+\| \langle x\rangle^{1/2+}\,J^{a-2}u\|_2)=:A_{412}+A_{422}.
	\end{equation*}
	We see that
	\begin{equation*}
		A_{422}\leq c(\|u\|_2+\|xJ^{a-2}u\|_2)\leq c(\|u\|_2+\|xu\|_2), 
	\end{equation*}
	and by Lemma \ref{Hardy-type} using that $\widehat{u}(0)=0$ so $\widehat{u}(\xi)=\xi\,\widehat{v}(\xi)$ with $\widehat{v}\in H^2(\mathbb{R})$,
	\begin{equation*}
		A_{412}\leq c\|J^{1/2+}(|\xi|^{a-1}\operatorname{sgn}(\xi)\,\chi(\xi)\widehat{v} (\xi))\|_2.
	\end{equation*}
	Since $|\xi|^{a-1} \operatorname{sgn}(\xi) \chi(\xi) \in H^{a-1/2}(\mathbb{R})$ (with $a-1/2 > 1/2$) we complete the proof. 
	As before, we remark that the above formal proof can be made rigorously justified by using the commutator estimates found in D.~Li (Corollary 1.4) and an approximation argument.
	
	\medskip

	Finally, we consider the remaining range of weights $\alpha\in(7/2,a+5/2)$, where the function $x^{3}u$ plays the role of $x^{2}u$.
	
	\noindent \underline{\bf Case $\alpha\in (7/2,a+5/2)$}. By the previous case, we assume that 
	\begin{equation}\label{class7}
		u\in C([-T,T]:H^{7a/2}(\mathbb{R}))\cap C([-T,T]:L^2(|x|^{7}\, dx)),
	\end{equation}
	and $\widehat{u}(0,t)=0$ for all $t\in [-T,T]$. Moreover, it follows
	\begin{equation}\label{class6.1}
		\langle x\rangle|D_x|^{a/2}(x^2 u), \,  |D_x|^{a/2}(x^3 u)\,\in L^{\infty}([-T,T]:L^2(\mathbb{R})).  
	\end{equation}
	Here, we write $\alpha=7/2+\theta$ with $\theta<a-1$. We multiply the equation in \eqref{dgbo} by $x^7 \langle x\rangle^{2\theta}_N \,u$ and integrate the result in the space variable. From the hypothesis and previous results, we just need to consider the contribution of the linear term:
	\begin{equation*}
		A_{0} = -\int_{\mathbb{R}} \partial_{x} |D_{x}|^{a} u \cdot x^{7} \langle x\rangle^{2\theta}_N\,u \, dx. 
	\end{equation*}
	Using the identity \eqref{id3} we have
	\begin{equation}
		\label{b71}
		\begin{aligned}
			A_{0} &= -\int_{\mathbb{R}} \partial_{x} |D_{x}|^{a}(x^{3}u) x^{3}u\,\langle x\rangle^{2\theta}_N\, x \, dx + 3(a+1)\int_{\mathbb{R}} |D_{x}|^{a}(x^2u) x\, x^{3}u\,\langle x\rangle^{2\theta}_N \, dx \\
			&\hskip12pt +3 a(a+1)\int_{\mathbb{R}} |D_{x}|^{a-1}\mathcal{H}(xu) \,x\, x^{3}u \, \langle x\rangle^{2\theta}_N\,dx -a(a^2-1)\int_{\mathbb{R}} |D_x|^{a-2}u \,x\, x^{3}u \, \langle x\rangle^{2\theta}_N\,dx\\
			&=: A_{1} + A_{2} + A_{3}+A_{4}.
		\end{aligned}
	\end{equation}
	By \eqref{eq1abc}, we get 
	\begin{equation*}
		\begin{aligned}
			A_1&=\Big(\frac{a+1}{2}\Big)\int_{\mathbb{R}} (|D_x|^{a/2}(x^3u))^2)(x\langle x\rangle^{2\theta}_N)'\,dx\\
			&\hskip12pt- \int_{\mathbb{R}} (x^3u) \tilde{R}_0(a+1)(x^3u)\,dx,
		\end{aligned}
	\end{equation*}
	and 
	\begin{equation*}
		\begin{aligned}
			\Big| \int_{\mathbb{R}} (x^3u) \tilde{R}_0(a+1)(x^3u)\,dx\Big| &\leq c\||D_x|^{-\sigma}(x^3u)\|_2^2\,\||D_x|^{\sigma}\,\tilde{R}_0(a+1)\,|D_x|^{\sigma}\|_{L^{2}\to L^{2}}\\
			&\leq c \| \langle x\rangle^{7/2} u\|_2^2 \,\|\widehat{\partial_x^3(\langle x\rangle^{2\theta}_N\,x)}\|_{L^1_{\xi}}<\infty,
		\end{aligned}
	\end{equation*}
	if $\theta\in (0,1)$  by a previous result by taking $2\sigma+a+1=3$.
	
	Next, using \eqref{I2} we have
	\begin{equation}
		\begin{aligned}
			\label{N1}
			A_{2} &= 3(a+1)\int_{\mathbb{R}} |D_{x}|^{a}(x^{3}u)x^{3}u\langle x\rangle^{2\theta}_N\, dx - 3(a+1)a\int_{\mathbb{R}} |D_{x}|^{a-1}\mathcal{H}(x^{2}u)x^{3}u\langle x\rangle^{2\theta}_N\, dx \\
			&=:A_{21}+A_{22}.
		\end{aligned}
	\end{equation}
	We write
	\begin{align*}
		A_{21} &= 3(a+1)\int_{\mathbb{R}} (|D_{x}|^{a/2}(x^{3}u))^{2}\langle x\rangle^{2\theta}_N \, dx \\
		&\hskip10pt+ 3(a+1) \int_{\mathbb{R}} |D_x|^{a/2}(x^3u) [|D_{x}|^{a/2}; \langle x\rangle^{2\theta}_N](x^{3}u)\,dx \\
		&=: A_{211} + A_{212}.
	\end{align*}
	Using the identity \eqref{ExtraExpansion}, we have
	\begin{equation*}
		\begin{aligned}
			A_{211}=&3(a+1)\int_{\mathbb{R}} |D_x|^{a/2}(x^3 u)\mathcal{S}_{0,1}(a/2)(x^3u)\,dx\\
			&+\frac{3(a+1)a}{2} \int_{\mathbb{R}} |D_x|^{a-1}(x^3u)\mathcal{H}\big((\langle x \rangle^{2\theta}_N)'x^3 u\big)\,dx.      
		\end{aligned}
	\end{equation*}
	To bound the first term in the identity above, we use Theorem \ref{gv2} to write
	\begin{align*}
		\Big| \int_{\mathbb{R}} |D_{x}|^{a/2}&(x^{3}u)\,\mathcal{S}_{0,1}(a/2)(x^{3}u)dx \Big| \\
		&= \left| \int_{\mathbb{R}} |D_{x}|^{-\gamma_{1}}(|D_{x}|^{a/2}(x^{3}u)) |D_{x}|^{\gamma_{1}}\mathcal{S}_{0,1}(a/2)|D_{x}|^{\gamma_{2}} (|D_{x}|^{-\gamma_{2}}(x^{3}u))dx \right| \\
		&\le c \|\langle x\rangle^{7/2}u\|_{2}^{2} \,\||D_{x}|^{\gamma_{1}}\mathcal{S}_{0,1}(a/2)|D_{x}|^{\gamma_{2}}\|_{L^{2}\to L^{2}},
	\end{align*}
	which is bounded for $\theta<1$ by taking $\ a/2 + \gamma_{1} + \gamma_{2} = 2$ since from previous argument
	\begin{equation*}
		\||D_{x}|^{\gamma_{1}}\mathcal{S}_{0,1}(a/2)|D_{x}|^{\gamma_{2}}\|_{L^{2}\to L^{2}} \le c \|\widehat{\partial_{x}^{2}(\langle x\rangle^{2\theta}_N)}\|_{L^{1}_{\xi}} < \infty.
	\end{equation*}
	To control the last term in the identity $A_{211}$, 
	\begin{equation*}
		\begin{aligned}
			\Big|\int_{\mathbb{R}} |D_x|^{a-1}(x^3u)\mathcal{H}\big((\langle x \rangle^{2\theta}_N)'x^3 u\big)\,dx\Big|\leq & c\||D_x|^{a-1}(x^3 u)\|_2\|\langle x\rangle^{2+2\theta} u\|_2 \\ 
			\leq & c\||D_x|^{a/2}(x^3 u)\|_2\|\langle x\rangle^{2+2\theta} u\|_2, 
		\end{aligned}    
	\end{equation*}
	which by \eqref{class7} and \eqref{class6.1} is controlled when $2+2\theta\leq 7/2$, i.e., $\theta\leq 3/4$.
	
	It remains to consider $A_{22}$ in \eqref{N1}.  Then
	\begin{align*}
		|A_{22}| &\leq c \int_{\mathbb{R}} \big||D_{x}|^{a-1}\mathcal{H}(x^{2}u)\langle x\rangle^{3+2\theta}u \big| dx  \\
		&\le \|\langle x\rangle^{\nu}|D_{x}|^{a-1}\mathcal{H}(x^{2}u)\|_{2} \|\langle x\rangle^{3+2\theta-\nu}u\|_{2} =: E_{1} \cdot E_{2}.
	\end{align*}
	We have
	\begin{align*}
		E_{1} &\leq c\|J^{\nu}(|\xi|^{a-1}\operatorname{sgn}(\xi)\chi(\xi))(\widehat{x^{2}u})(\xi)\|_{2} + c\|\langle x\rangle^{\nu}J^{a-1}(x^{2}u)\|_{2} \\
		&=: E_{11} + E_{12}.
	\end{align*}
	Since $|\xi|^{a-1}\operatorname{sgn}(\xi)\chi \in H^{(a-1/2)}(\mathbb{R})$ and $\widehat{x^{2}u} \in H^{3/2}(\mathbb{R})$, one has that $E_{11} < \infty$ if $\nu < a - 1/2$. Also, if  $\nu\leq 1$
	\begin{align*}
		E_{12} &\le c \|J^{a-1}(x^{2}u)\|_{2} + c \|x J^{a-1}(x^{2}u)\|_{2} \\
		&\le c \|J^{a-1}(x^{2}u)\|_{ 2} + \||D_{x}|^{a/2}(x^{3}u)\|_{2},
	\end{align*}
	which is bounded after integration in time from the previous result.
	Also for $\nu < a - 1/2$ one has that $E_{2} < \infty$ if $3+2\theta-\nu \le 7/2$ implies that $\theta < a/2$.
	
	\medskip
	
	Thus, all the previous terms are collected and bounded for $\theta \le 3/4$ and $\theta < a/2$. Next we consider $A_{3}$ in \eqref{b71} and estimate
	\begin{equation*}
		|A_{3}| \leq c \| \langle x \rangle^{\nu}|D_{x}|^{a-1}\mathcal{H}(xu)\|_2\|\langle x\rangle^{4+2\theta-\nu}u\|_2,
	\end{equation*}
	where as before
	\begin{equation*}
		\begin{aligned}
			\|\langle x\rangle^{\nu}|D_{x}|^{a-1}H(xu)\|_{2} \le c \big(\|J^{\nu}(|\xi|^{a-1}\operatorname{sgn}(\xi)\chi(\xi)\widehat{xu}(\xi))\|_{2} + \|\langle x\rangle^{\nu}J^{a-1}(xu)\|_{2}\big) =: F_{1} + F_{2}.
		\end{aligned}
	\end{equation*}
	Using that
	\begin{equation*}
		\widehat{xu}(\xi) = c\partial_{\xi}\widehat{u}(\xi) = c\partial_{\xi}(\xi\widehat{v}(\xi)) = c\widehat{v}(\xi) + c\xi\partial_{\xi}\widehat{v}(\xi)
	\end{equation*}
	we see that with $\widehat{v} \in H^{5/2}(\mathbb{R})$
	\begin{equation*}
		|\xi|^{a-1}\operatorname{sgn}(\xi)\chi(\xi)\widehat{xu}(\xi) = c|\xi|^{a-1}\operatorname{sgn}(\xi)\chi(\xi)\widehat{v}(\xi) + c|\xi|^{a}\chi(\xi)\widehat{v}(\xi).
	\end{equation*}
	Thus $F_{1} < \infty$ if and only if $\nu < a - 1/2$ and $F_{2} < \infty$ if $\nu \le 2$. Hence $|A_{3}| < \infty$ if $4+2\theta-\nu \le 7/2$, i.e., 
	$\theta < (a-1)/2$.
	Finally, we consider $A_{4}$ in \eqref{N1}
	\begin{align*}
		|A_{4}| &\leq  c \|\langle x\rangle^{\nu}|D_{x}|^{a-2}u\|_{2} \|\langle x\rangle^{4+2\theta-\nu}u\|_{2}.
	\end{align*}
	Writing
	\begin{equation*}
		\|\langle x\rangle^{\nu}|D_{x}|^{a-2}u\|_{2} \le c \|J^{\nu}(|\xi|^{a-2}\widehat{u}(\xi)\chi(\xi))\|_{2} + \|\langle x\rangle^{\nu}J^{a-2}u\|_{2}
	\end{equation*}
	and using that $\widehat{u}(\xi) = \xi\widehat{v}(\xi)$ with $\widehat{v} \in H^{5/2}(\mathbb{R})$ we get
	\begin{equation*}
		\|J^{\nu}(|\xi|^{a-2}\widehat{u}(\xi)\chi(\xi))\|_{2} < \infty \iff \nu < a - 1/2
	\end{equation*}
	and 
	\begin{equation*}
		\|\langle x\rangle^{\nu}J^{a-2}u\|_{2} < \infty \iff \nu \leq 3.
	\end{equation*}
	For $A_{4}$ one has $|A_{4}| < \infty$, if $4+2\theta-\nu \le 7/2$, i.e., $\theta < (a-1)/2$.
	
	\medskip
	Collecting the above information we get that for $\theta < (a-1)/2$
	\begin{equation*}
		\langle x\rangle^{\theta}|D_{x}|^{a/2}(x^{3}u) \in L^{2}([t_{1}, t_{2}]:L^2(\mathbb{R}))
	\end{equation*}
	and 
	\begin{equation*}
		|x|^{7/2+\theta}u \in L^{\infty}([t_{1}, t_{2}]: L^{2}(\mathbb{R})), \;\;\;\;\;\text{for}\;\;\;\theta < \frac{a-1}{2}.
	\end{equation*}
	Reapplying this argument we extend the results above to (with $\nu < a - 1/2$)
	\begin{equation*}
		4+2\theta - \nu < \frac{7}{2} + \frac{a-1}{2} \implies \theta < 3(a-1)/4.
	\end{equation*}
	Iterating this process we extend it to $\theta < a-1$, i.e., for $\alpha < 7/2 + \theta = 5/2 + a,$ one has 
	\begin{equation*}
		\langle x\rangle^{\theta}|D_{x}|^{a/2}(x^{3}u) \in L^{2}([t_{1}, t_{2}]: L^{2}(\mathbb{R})),\;\;\; \quad |x|^{\alpha}u \in L^{\infty}([t_{1}, t_{2}]: L^{2}(\mathbb{R}))
	\end{equation*}
	which is the desired result.
	\section{Proof of Theorem \ref{AA2} }
	In this section we prove Theorem \ref{AA2}. We first deduce the persistence properties in parts (i) and (iii) by means of weighted energy estimates, where the main difficulty is to control the dispersive term with the sharp regularity $s\geq\max\{a/2,\,a\alpha\}$. Then we prove the unique continuation principles in parts (ii) and (iv), which are based on the integral formulation of \eqref{dgbo}.
	
	\subsection{Proof of Theorem \ref{AA2} parts (i) and (iii)}

	We consider $m\in \mathbb{N}_0$, $0<\theta\leq 1$. Multiplying the equation in \eqref{dgbo} by $x^{2m}\langle x\rangle_N^{2\theta} u$ and integrating in space, it follows that
	\begin{equation}\label{energyLWP}
		\begin{aligned}
			\frac{1}{2}\frac{d}{dt}\int_{\mathbb{R}} \big(x^m \langle x\rangle_N^{\theta} u\big)^2\, dx -&\underbrace{\int_{\mathbb{R}} \partial_x |D_x|^{a} u \big(x^{2m}\langle x \rangle_N^{2\theta} u\big)\, dx}_{\mathcal{B}_1} +\int_{\mathbb{R}} u\partial_x u \big(x^{2m}\langle x\rangle^{2\theta}_N u\big)\, dx=0.    
		\end{aligned}
	\end{equation}
	To deduce Theorem \ref{AA2} (i) and (iii), by the arguments in \cite{FONLP} and the estimate \eqref{nonlinearE1}, it suffices to control the dispersive term $\mathcal{B}_1$ while preserving the regularity at the condition $s\geq \max\{\frac{a}{2};\, a \alpha\}$, with $\alpha=m+2\theta$. It is also worth emphasizing that by the symmetry of the equation in \eqref{dgbo}, setting $w(x,t)=u(-x,-t)$, it is enough to study the energy estimates for positive times. 
	
	We summarize the estimate for $\mathcal{B}_1$ in the following lemma.
	\begin{lem}\label{DispersiveEstimate}
		Let $1<a<2$, $m\in\mathbb{N}_0$, $0<\theta\leq1$, and let
		$\mathcal{B}_1$ be the term defined in the weighted energy identity above. Then, there exists $c_{a,m,\theta}>0$ independent of $N$ such that
		\begin{equation*}
			\begin{aligned}
				|\mathcal{B}_1| \leq c_{a,m,\theta}
				\Big(\|J^{\max\{a/2,\, a(m+\theta)\}}u\|_{2}^2
				+\|\langle x\rangle^m u\|_{2}^2
				+\|x^m\langle x\rangle_N^\theta u\|_{2}^2\Big),
			\end{aligned}
		\end{equation*}
		provided that one of the following two conditions holds:
		\begin{equation}
			\left\{\begin{aligned}
				&m+\theta<a+3/2,\qquad\text{or}\\
				&a+3/2\leq m+\theta<a+5/2
				\quad\text{and}\quad \widehat u(0)=0.
			\end{aligned}\right.
		\end{equation}
	\end{lem}
	\begin{rem}
		The argument behind the proof of (i) and (iii) in Theorem \ref{AA2} is an iteration in the weight size. Starting with $m=0$, we write $\alpha=m+\theta$ and assume that persistence
		is already known at the integer weight $m$, i.e., $u\in C([-T,T]: L^2(|x|^{2m}\,dx)$. Lemma \ref{DispersiveEstimate}, with
		$s\geq \max\{a/2,\, \,a\alpha\}$, then controls the dispersive term $\mathcal{B}_1$ with weight $m+\theta$. This is why the fractional weight is first replaced by its bounded approximation $\langle x\rangle_N^\theta$ (the truncation justifies the weighted energy estimate). Gronwall's inequality then gives a bound independent of $N$, and thus taking $N\to\infty$, we recover the full weight $m+\theta$. Repeating this step for successive integers $m$ and proper $0\leq \theta\leq 1$ proves Theorem \ref{AA2}. 
	\end{rem}
	\begin{proof}[Proof of Lemma \ref{DispersiveEstimate}]
		We apply Proposition \ref{GeneralConmutator} to write     
		\begin{equation}\label{eqPersistence1}
			\begin{aligned}
				\int_{\mathbb{R}} \partial_x |D_x|^{a} u \big(x^{2m}\langle x \rangle_N^{2\theta} u\big)\, dx =&\int_{\mathbb{R}} \partial_x|D_x|^{a}(x^m u)(x^m u)\langle x\rangle_N^{2\theta}\, dx \\
				&+ \sum_{k=1}^m \binom{m}{k} \int_{\mathbb{R}} \widetilde{|D_x|^{a+1,k}}(x^{m-k}u)(x^m u)\langle x\rangle^{2\theta}_N \,dx\\
				=:& \mathcal{B}_{1, 0}+\sum_{k=1}^m \binom{m}{k}\mathcal{B}_{1, k}.
			\end{aligned}
		\end{equation}
		We adopt the convention that when $m=0$ the empty sum $\sum_{k=1}^{0}(\cdots)=0$. We remark that $\mathcal{B}_{1,k}$ is estimated under the following dispersion regimes: when $1<a\leq \frac{3}{2}$, we consider $m=0,1,2,3$, and hence only the terms with $k=1,2,3$. When $\frac{3}{2}<a<2$, we take $m=0,1,2,3,4$, and estimate the terms with $k=1,2,3,4$.
		\\ \\  
		To avoid repetition, we restrict ourselves to control $\mathcal{B}_{1,0}$, which captures the $a/2$ limitation of regularity. Despite this, the ideas developed in the proof of Theorem \ref{main} can be adapted to bound the remaining factors $\mathcal{B}_{1,k}$, $k\geq 1$.
		\\ \\
		\underline{\bf Estimate for $\mathcal{B}_{1, 0}$}. With the aim of applying \eqref{eqextenden2}, we write $\partial_x|D_x|^{a}=-\mathcal{H}|D_x|^{a}$ to get
		\begin{equation}\label{eqnonliner1}
			\begin{aligned}
				\langle x\rangle_N^{\theta}&\partial_x  |D_x|^{a}(x^m u) \\
				=&-[\langle x\rangle_N^{\theta};\mathcal{H}]|D_x|^{a+1}(x^m u)+\mathcal{B}_{1,0,1}\\
				&+\mathcal{H}\Big((a+1)\partial_x\big(\langle x\rangle_N^{\theta}\big)\mathcal{H}|D_x|^{a}(x^m u)-\frac{(a+1)a}{2}\partial_x^2\big(\langle x\rangle_N^{\theta}\big)|D_x|^{a-1}(x^m u)\Big)\\
				&-\mathcal{H}|D_x|^{a+1}(\langle x\rangle_N^{\theta}x^m u),
			\end{aligned}   
		\end{equation}
		where
		\begin{equation*}
			\mathcal{B}_{1,0,1}:=\mathcal{H}\Big([|D_x|^{a+1};\langle x\rangle_N^{\theta}](x^m u)-(a+1)\partial_x\big(\langle x\rangle_N^{\theta}\big)\mathcal{H}|D_x|^{a}(x^m u)+\frac{(a+1)a}{2}\partial_x^2\big(\langle x\rangle_N^{\theta}\big)|D_x|^{a-1}(x^m u)\Big),    
		\end{equation*}
		and by Lemma \ref{extendenlemma}, it satisfies   
		\begin{equation*}
			\begin{aligned}
				\|\mathcal{B}_{1,0,1}\|_{2}\leq c \||D_x|^{a+1}\big(\langle x\rangle_N^{\theta}\big)\|_{\infty}\|x^m u\|_{2}.        
			\end{aligned}
		\end{equation*}
		Note that by previous arguments, $\||D_x|^{a+1}\big(\langle x\rangle_N^{\theta}\big)\|_{\infty}\leq c <\infty$, where $c$ is independent of $N$.
		
		Returning to the integral defining $\mathcal{B}_{1,0}$ and using that $\mathcal{H}|D_x|^{a+1}$ is skew-symmetric, we conclude that the last term in \eqref{eqnonliner1} does not contribute to the estimate. 
		
		To control the first term on the right-hand side of \eqref{eqnonliner1}, we use \eqref{fracComwell2} to get 
		\begin{equation*}
			\|[\langle x\rangle_N^{\theta};\mathcal{H}]|D_x|^{a+1}(x^m u)\|_{2}\leq c \||D_x|^{a+1}\big(\langle x\rangle_N^{\theta}\big)\|_{\infty}\|x^m u\|_{2}.   
		\end{equation*}
		We now restrict our attention to the third-to-last term on the right-hand side of \eqref{eqnonliner1}. Using commutators, we rewrite
		\begin{equation}\label{eqnonliner2}
			\begin{aligned}
				\partial_x\big(\langle x\rangle_{N}^{\theta}\big)\mathcal{H}|D_x|^{a}(x^m u)=&[\partial_x\big(\langle x\rangle_{N}^{\theta}\big);|D_x|^{a}]\mathcal{H}(x^m u)+|D_x|^{a}\big(\partial_x\big(\langle x\rangle_{N}^{\theta}\big)\mathcal{H}(x^m u)\big)\\
				=& [\partial_x\big(\langle x\rangle_{N}^{\theta}\big);|D_x|^{a}]\mathcal{H}(x^m u)+a\partial_x^2\big(\langle x\rangle_N^{\theta}\big)\mathcal{H}^2|D_x|^{a-1}(x^m u)\\
				&-a\partial_x^2\big(\langle x\rangle_N^{\theta}\big)\mathcal{H}^2|D_x|^{a-1}(x^m u)+|D_x|^{a}\big(\partial_x\big(\langle x\rangle_{N}^{\theta}\big)\mathcal{H}(x^m u)\big).
			\end{aligned}    
		\end{equation}
		Recall from \eqref{eqnonliner1} that the Hilbert transform $\mathcal{H}$ acts on the identity above. We now proceed to estimate the terms arising in the above decomposition.
		Using that $\mathcal{H}$ is a bounded operator in $L^2$, we apply Lemma \ref{extendenlemma} inequality \eqref{eqextenden1} to infer
		\begin{equation*}
			\begin{aligned}
				\|\mathcal{H}\big([\partial_x\big(\langle x\rangle_{N}^{\theta}\big),|D_x|^{a}]\mathcal{H}(x^m u)+a\partial_x^2\big(\langle x\rangle_N^{\theta}\big)\mathcal{H}^2|D_x|^{a-1}(x^m u) \big)\|_{2}\leq c  \||D_x|^{a}\partial_x\big(\langle x\rangle_N^{\theta}\big)\|_{\infty}\|x^m u\|_{2}.   
			\end{aligned}    
		\end{equation*}
		Observe that arguing as above, we have $ \||D_x|^{a}\partial_x\big(\langle x\rangle_N^{\theta}\big)\|_{\infty}= \|\mathcal{H}|D_x|^{a+1}\big(\langle x\rangle_N^{\theta}\big)\|_{\infty}\leq c$ with $c>0$ independent of $N$.  Next, using that $\mathcal{H}^2=-I$, we return to the integral defining $\mathcal{B}_{1,0}$ to deduce
		\begin{equation*}
			\begin{aligned}
				\int_{\mathbb{R}} \mathcal{H}\Big(\partial_x^2\big(\langle x\rangle_N^{\theta}\big)\mathcal{H}^2|D_x|^{a-1}(x^m u)\Big) \big((x^m u)\langle x\rangle_N^{\theta}\big)\, dx=&-\int_{\mathbb{R}} \big([\mathcal{H};\partial_x^2\big(\langle x\rangle_N^{\theta}\big)]|D_x|^{a-1}(x^m u)\big)\big((x^m u)\langle x\rangle_N^{\theta}\big)\, dx \\
				&-\int_{\mathbb{R}} \mathcal{H}(x^m u)|D_x|^{a-1}\big((x^m u)\langle x\rangle_N^{\theta}\partial_x^2\big(\langle x\rangle_N^{\theta}\big)\big)\, dx.
			\end{aligned}   
		\end{equation*}
		By  Cauchy-Schwarz inequality and \eqref{fracComwell2}, we get 
		\begin{equation*}
			\begin{aligned}
				\Big|\int_{\mathbb{R}} \big([\mathcal{H};\partial_x^2 &\big(\langle x\rangle_N^{\theta}\big)]|D_x|^{a-1}(x^m u)\big)\big((x^m u)\langle x\rangle_N^{\theta}\big)\, dx\Big|\\
				\leq c & \|[\mathcal{H};\partial_x^2\big(\langle x\rangle_N^{\theta}\big)]|D_x|^{a-1}(x^m u)\|_{2} \|x^m \langle x\rangle_N^{\theta} u\|_{2}\\
				\leq c & \||D_x|^{a-1}\partial_x^2\big(\langle x\rangle_N^{\theta}\big)\|_{\infty}\|x^m u\|_{2}\|x^m \langle x\rangle_N^{\theta} u\|_{2}.
			\end{aligned}   
		\end{equation*}
		If $m\geq 1$, we notice that $|x^m\langle x\rangle_N^{\theta}\partial_x^2\big(\langle x\rangle_N^{\theta}\big)|\leq c \langle x\rangle_N^{m-1+\theta}$ with constant independent of $N$. Thus, by fractional  Leibniz 
		\begin{equation*}
			\begin{aligned}
				\||D_x|^{a-1}\big((x^m u)\langle x\rangle_N^{\theta}\partial_x^2\big(\langle x\rangle_N^{\theta}\big)\big)\|_{2}= &\||D_x|^{a-1}\Big(\frac{x^m \langle x \rangle_N^{\theta}\partial_x^2\big(\langle x\rangle_N^{\theta}\big)}{\langle x\rangle_N^{m-1+\theta}} (\langle x\rangle^{m-1+\theta}_N u )\Big)\|_{2}\\
				\leq & c \Big\||D_x|^{a-1}\Big(\frac{x^m \langle x \rangle_N^{\theta}\partial_x^2\big(\langle x\rangle_N^{\theta}\big)}{\langle x\rangle_N^{m-1+\theta}}\Big)\Big\|_{\infty} \|\langle x\rangle^{m-1+\theta}_N u \|_{2}\\
				&+c \Big\|\frac{x^m \langle x \rangle_N^{\theta}\partial_x^2\big(\langle x\rangle_N^{\theta}\big)}{\langle x\rangle^{m-1+\theta}}\Big\|_{\infty} \||D_x|^{a-1}\big(\langle x\rangle^{m-1+\theta}_N u \big)\|_{2} \\
				\leq & c \|J^{a-1}\big(\langle x\rangle_N^{m-1+\theta} u\big)\|_{2}.
			\end{aligned}    
		\end{equation*}
		From \eqref{inter1x}, which is valid for weights $\langle x \rangle_N$ with a constant independent of $N$ (see \cite{FONLP}) it follows that  
		\begin{equation*}
			\begin{aligned}
				\Big|\int_{\mathbb{R}} \mathcal{H}(x^m u) &|D_x|^{a-1}\big((x^m u)\langle x\rangle_N^{\theta}\partial_x^2\big(\langle x\rangle_N^{\theta}\big)\big)\, dx\Big|\\
				\leq c & \|x^m u\|_{2}\||D_x|^{a-1}\big((x^m u)\langle x\rangle_N^{\theta}\partial_x^2\big(\langle x\rangle_N^{\theta}\big)\big)\|_{2}\\
				\leq c  & \|x^m u\|_{2}\big(\|J^{(a-1)(m+\theta)}u\|_{2}+\|\langle x\rangle_N^{m+\theta} u\|_{2}\big)\\
				\leq c & \|x^m u\|_{2}\big(\|J^{(a-1)(m+\theta)}u\|_{2}+\|\langle x\rangle_N^{\theta} u\|_{2}+\|x^m\langle x\rangle_N^{\theta}  u\|_{2}\big).
			\end{aligned}    
		\end{equation*}
		If $m=0$, notice that $|\langle x\rangle_N^{\theta}\partial_x^2\big(\langle x\rangle_N^{\theta}\big)|\leq c$ with $c>0$ independent of $N$, then arguing as above, one has
		\begin{equation*}
			\begin{aligned}
				\||D_x|^{a-1}\big(u \langle x\rangle_N^{\theta}\partial_x^2\big(\langle x\rangle_N^{\theta}\big)\big)\|_{2}\leq c \|J^{a-1} u \|_{2},
			\end{aligned}    
		\end{equation*}
		so that
		\begin{equation*}
			\begin{aligned}
				\Big|\int_{\mathbb{R}} \mathcal{H} u  |D_x|^{a-1}\big( u\langle x\rangle_N^{\theta}\partial_x^2\big(\langle x\rangle_N^{\theta}\big)\big)\, dx\leq c \|u\|_{2}\|J^{(a-1)}u\|_{2},
			\end{aligned}    
		\end{equation*}
		which is a valid estimate provided that $\|J^{(a-1)}u\|_{2}\leq \|J^{a/2}u\|_{2}$. We continue with the last term in \eqref{eqnonliner2}. We write
		\begin{equation}\label{decompAux}
			\begin{aligned}
				\int_{\mathbb{R}} & \mathcal{H}|D_x|^{a}\big(\partial_x\big(\langle x\rangle_{N}^{\theta}\big) \mathcal{H}(x^m u)\big) (x^m u)\langle x \rangle_N^{\theta}\, dx  \\
				=&\int_{\mathbb{R}} \mathcal{H}|D_x|^{a}\big([\partial_x\big(\langle x\rangle_{N}^{\theta}\big); \mathcal{H}](x^m u)\big) (x^m u)\langle x \rangle_N^{\theta}\, dx -\int_{\mathbb{R}} |D_x|^{a}\big(\partial_x\big(\langle x\rangle_{N}^{\theta}\big) x^m u\big) (x^m u)\langle x \rangle_N^{\theta}\, dx.
			\end{aligned}    
		\end{equation}
		If $m\geq 1$, given that $|\partial_x\big(\langle x\rangle_{N}^{\theta}\big)x^m|\leq c \langle x\rangle^{m-1+\theta}_N$, we arrive at
		\begin{equation*}
			\begin{aligned}
				\Big|\int_{\mathbb{R}} &\mathcal{H}|D_x|^{a}\big(\partial_x\big(\langle x\rangle_{N}^{\theta}\big) \mathcal{H}(x^m u)\big) (x^m u)\langle x \rangle_N^{\theta}\, dx \Big|\\
				\leq &\Big|\int_{\mathbb{R}} \mathcal{H}|D_x|^{a}\big([\partial_x\big(\langle x\rangle_{N}^{\theta}\big);\mathcal{H}](x^m u)\big) (x^m u)\langle x \rangle_N^{\theta}\, dx \Big|\\
				&+\Big|\int_{\mathbb{R}} |D_x|^{a}\big(\partial_x\big(\langle x\rangle_{N}^{\theta}\big)(x^m u)\big) (x^m u)\langle x \rangle_N^{\theta}\, dx \Big|\\
				\leq & c\|J^{-2}|D_x|^{a} J^2 [\partial_x\big(\langle x\rangle_{N}^{\theta}\big);\mathcal{H}](x^m u)\|_{2}\|x^m \langle x \rangle^{\theta}_N u\|_{2}+c\|J^{a}\big(\langle x \rangle_N^{m-1+\theta}u\big)\|_{2}\|x^m \langle x \rangle^{\theta}_N u\|_{2}\\
				\leq & c \big( \|J^{a(m+\theta)}u\|_{2}+\|\langle x \rangle_N^{\theta}u\|_{2}+\|x^m \langle x \rangle^{\theta}_N u\|_{2}\big)\|x^m \langle x \rangle^{\theta}_N u\|_{2}.
			\end{aligned}    
		\end{equation*}
		When $m=0$, to preserve regularity at the level $a/2$, we further decompose \eqref{decompAux} as follows 
		\begin{equation*}
			\begin{aligned}
				\int_{\mathbb{R}} |D_x|^{a}\big(\partial_x\big(\langle x\rangle_{N}^{\theta}\big) u\big)  u\langle x \rangle_N^{\theta}\, dx=&\int_{\mathbb{R}} [|D_x|^{a};\partial_x\big(\langle x\rangle_{N}^{\theta}\big)] u \, (u\langle x \rangle_N^{\theta})\, dx\\
				&+\int_{\mathbb{R}} \partial_x\big(\langle x\rangle_{N}^{\theta}) |D_x|^{a} u \, (u\langle x \rangle_N^{\theta})\, dx,
			\end{aligned}    
		\end{equation*}
		and for the last term above, we use the decomposition
		\begin{equation*}
			\begin{aligned}
				\int_{\mathbb{R}} \partial_x\big(\langle x\rangle_{N}^{\theta}) |D_x|^{a} u \, (u\langle x \rangle_N^{\theta})\, dx=&\int_{\mathbb{R}} \frac{\langle x \rangle_N^{\theta}\partial_x\big(\langle x\rangle_{N}^{\theta})}{\langle x \rangle_N^{2\theta-1}} [\langle x \rangle_N^{\theta-1/2};|D_x|^{a}] u \, (\langle x \rangle_N^{\theta-1/2} u)\, dx\\
				&+ \int_{\mathbb{R}} \frac{\langle x \rangle_N^{\theta}\partial_x\big(\langle x\rangle_{N}^{\theta})}{\langle x \rangle_N^{2\theta-1}} |D_x|^{a}\big(\langle x \rangle_N^{\theta-1/2} u \big)\, (\langle x \rangle_N^{\theta-1/2} u)\, dx.  
			\end{aligned}    
		\end{equation*}
		Since $\frac{\langle x \rangle_N^{\theta}\partial_x\big(\langle x\rangle_{N}^{\theta})}{\langle x \rangle_N^{2\theta-1}}$ is a bounded with bounded derivative (both independent of $N$), using Theorem \ref{kprefi} and interpolation, we obtain 
		\begin{equation*}
			\begin{aligned}
				\Big|\int_{\mathbb{R}} &\mathcal{H}|D_x|^{a}\big(\partial_x\big(\langle x\rangle_{N}^{\theta}\big) \mathcal{H}u\big) \, u\langle x \rangle_N^{\theta}\, dx \Big|\\
				\leq & c\|J^{-2}|D_x|^{a}J^2 [\partial_x\big(\langle x\rangle_{N}^{\theta}\big);\mathcal{H}] u\|_{2}\| \langle x \rangle^{\theta}_N u\|_{2}+\|u\|_2\|\langle x \rangle_N^{\theta-1/2}u\|_2+\|J^{a/2}\big( \langle x \rangle^{\theta-1/2}_N u\big)\|_{2}^2\\
				\leq & c \big( \|J^{\max\{a/2,\,a\theta\}}u\|_{2}+\|\langle x \rangle_N^{\theta}u\|_{2}\big)\big(\|\langle x \rangle^{\theta}_N u\|_{2}+\|J^{\max\{a/2,\,a\theta\}}u\|_{2}\big).
			\end{aligned}    
		\end{equation*}
		It is worth emphasizing that if $0<\theta<\frac{1}{2}$, the estimate $\|J^{a/2}\big( \langle x \rangle^{\theta-1/2}_N u\big)\|_{2}^2\leq \|J^{a/2} u\|_{2}^2$ forces us to have $J^{\max\{a/2,\,a\theta\}}u \in L^2(\mathbb{R})$. This completes the bound of \eqref{eqnonliner2} and, in turn, the study of the third-to-last term of \eqref{eqnonliner1}. 
		
		Finally, we focus on the penultimate term in \eqref{eqnonliner1}. Using that $\mathcal{H}$ is a skew symmetric operator, we write
		\begin{equation*}
			\begin{aligned}
				\int_{\mathbb{R}} \mathcal{H}\Big(\partial_x^2\big(\langle x\rangle_N^{\theta}\big)|D_x|^{a-1}(x^m u)\Big) \big((x^m u)\langle x\rangle_N^{\theta}\big)\, dx=&\int_{\mathbb{R}} \big([\mathcal{H},\partial_x^2\big(\langle x\rangle_N^{\theta}\big)] |D_x|^{a-1}(x^m u)\big)\big((x^m u)\langle x\rangle_N^{\theta}\big)\, dx \\
				&-\int_{\mathbb{R}} (x^m u)\mathcal{H}|D_x|^{a-1}\big((x^m u)\langle x\rangle_N^{\theta}\partial_x^2\big(\langle x\rangle_N^{\theta}\big)\big)\, dx.
			\end{aligned}   
		\end{equation*}
		By the arguments employed above, we conclude that if $m\geq 1$
		\begin{equation*}
			\begin{aligned}
				\Big|\int_{\mathbb{R}} \mathcal{H}&\Big(\partial_x^2\big(\langle x\rangle_N^{\theta}\big)|D_x|^{a-1}(x^m u)\Big) \big((x^m u)\langle x\rangle_N^{\theta}\big)\, dx\Big|\\
				\leq c& \big(\|J^{(a-1)(m+\theta)}u\|_{2}+\|\langle x\rangle_N^{\theta} u\|_{2}+\|x^m \langle x\rangle_N^{\theta}  u\|_{2}\big)\|x^m \langle x\rangle_N^{\theta} u \|_{2},
			\end{aligned}   
		\end{equation*}
		and if $m=0$,
		\begin{equation*}
			\begin{aligned}
				\Big|\int_{\mathbb{R}} \mathcal{H}\Big(\partial_x^2\big(\langle x\rangle_N^{\theta}\big)|D_x|^{a-1} u\Big) \big( \langle x\rangle_N^{\theta} u\big)\, dx\Big|\leq c\|J^{(a-1)}u\|_{2}\|\langle x\rangle_N^{\theta} u\|_{2}.
			\end{aligned}   
		\end{equation*}
		Collecting all the previous estimates, we deduce
		\begin{equation*}
			|\mathcal{B}_{1, 0}|\leq c \big(\|J^{\max\{a/2,\, a(m+\theta)\}}u\|_{2}+\|\langle x\rangle_N^{\theta} u\|_{2}+\|x^m \langle x\rangle_N^{\theta} u \|_{2}\big)\big(\|J^{\max\{a/2,\, a(m+\theta)\}}u\|_{2}+\|x^m \langle x\rangle_N^{\theta} u \|_{2}\big).  
		\end{equation*}  
		This completes the proof
	\end{proof}
	\subsection{Proof of Theorem \ref{AA2} parts (ii) and (iv)}
	Next, we turn to the unique continuation principles. In contrast with parts (i) and (iii), the proof of parts (ii) and (iv) relies on the integral formulation of \eqref{dgbo}. Roughly speaking, we show that if the solution has the critical decay at two different times, then a term of the form $\widehat{u_0}(0)\,\xi|\xi|^{a-2}$ in the frequency variable, which comes from the lack of smoothness at the origin of the symbol $i\xi|\xi|^{a}$, must vanish.

	\begin{proof}[Proof of Theorem \ref{AA2} parts (ii) and (iv)]
		
		For the sake of brevity, we will only establish (ii), as the proof of (iv) is analogous and follows from applying $J^{-12}$ together with the identities in \cite{FONLP}. 
		
		To establish (ii), we begin with a further reduction. We consider only the case $1<a<\frac{3}{2}$, i.e., the case where $\frac{5}{2}<\frac{3}{2}+a<3$. We have that the result for dispersion $\frac{3}{2}\leq a<2$ follows similarly, with an additional local derivative estimate, so we omit further details.
		
		Let $u\in C([-T,T] : H^{a(3/2+a)^{-}}(\mathbb{R}))$ be a solution of the IVP \eqref{dgbo}. We assume that there exist two different times $t_1, t_2\in [-T,T]$ such that $|x|^{3/2+a}u(x,t_j)\in L^2(\mathbb{R}),\,j=1,2$. By the persistence results in (i) Theorem \ref{AA2}, we also have
		\begin{equation*}
			|x|^{(3/2+a)^{-}}u\in C([-T,T]: L^2(\mathbb{R})).    
		\end{equation*}
		Without loss of generality, we will assume that $t_1=0$ and $t_2>0$, that is, $u_0, u(t_2)\in L^2(|x|^{3+2\alpha}\, dx)$.
		
		We first require the following preliminary results.
		\begin{thm}[\cite{Stein1961}]
			\label{TheoSteDer}
			Let $0<b<1$ and $\frac{2}{1+2b}<p<\infty$, and set
			\begin{equation*}
				\mathcal{D}^b f(x)
				=\left(
				\int_{\mathbb{R}}
				\frac{|f(x)-f(y)|^2}{|x-y|^{1+2b}}\,dy
				\right)^{1/2}.
			\end{equation*}
			Then $f\in L_b^p(\mathbb{R})=J^{-b}L^p(\mathbb{R})$ if and only if
			$f,\mathcal{D}^b f\in L^p(\mathbb{R})$, and the norm equivalences hold
			\begin{equation*}
				\|J^b f\|_{p}
				\sim \big(\|f\|_{p}+\|\mathcal{D}^b f\|_{p} \big)
				\sim \big(\|f\|_{p}+\||D_x|^b f\|_{p}\big).
			\end{equation*}
			In particular,
			$\|\mathcal{D}^b f\|_{2}=c_b\||D_x|^b f\|_{2}$. Moreover,
			\begin{equation}\label{prelimneq}
				\begin{aligned}
					\|\mathcal{D}^b(fg)\|_{2} \leq & c_{b}\big( \|f\mathcal{D}^b g\|_{2}+\|g\mathcal{D}^b f\|_{2}\big),\\
					\|\mathcal{D}^b h\|_{\infty}\leq & c_b\big(\|h\|_{\infty}+\|h'\|_{\infty}\big).
				\end{aligned}    
			\end{equation}
		\end{thm}
		\begin{prop}\label{propUCPtwo1}
			It follows
			\begin{equation*}
				[x; e^{t\partial_x |D_x|^{a}}]f=-(1+a)t e^{t\partial_x |D_x|^{a}} |D_x|^{a}f.    
			\end{equation*}
			As a consequence,
			\begin{equation*}
				\begin{aligned}
					x^2 e^{t\partial_x|D_x|^{a}}f =&-(1+a) t e^{t\partial_x |D_x|^{a}}(x|D_x|^{a}f)
					+(1+a)^2 t^2 e^{t\partial_x |D_x|^{a}}|D_x|^{2a}f\\    
					&-(1+a) t e^{t\partial_x |D_x|^{a}}|D_x|^{a}(xf)+e^{t\partial_x |D_x|^{a}}(x^2f).
				\end{aligned} 
			\end{equation*}
			Where $e^{t\partial_x |D_x|^{a}}$ denotes the linear solution operator associated to $\partial_t v-\partial_x |D_x|^a v=0$.
		\end{prop}
		\begin{prop}\label{propUCPtwo2}
			Let $0<b<1$, then
			\begin{equation*}
				\begin{aligned}
					\|\langle x\rangle^{b}e^{t\partial_x |D_x|^{a}}f\|_{2}\leq c_{a,b} \langle t \rangle^b\big(\||D_x|^{a b} f\|_{2}  +\|\langle x\rangle^{b}f\|_{2} \big).  
				\end{aligned}   
			\end{equation*}    
		\end{prop}
		\begin{proof}
			Passing to the frequency domain,  the proof follows from \cite[Proposition 2.7]{FONLP}.
		\end{proof}
		\begin{prop}\label{propUCPtwo3}
			Let $0<b\leq 1$, $\beta\geq 1$, $m\in \mathbb{Z}^{+}$ with $\beta-m\leq 0$. Then 
			\begin{equation*}
				\|\langle x \rangle^{b}|D_x|^{\beta}J^{-m} f\|_{2}\leq c_{b,\beta,m}\|\langle x\rangle^b f\|_{2}.   
			\end{equation*}
			The same result holds for $\beta =0$. Moreover, the same conclusions are valid by replacing $|D_x|^{\beta}$ by $\partial_x |D_x|^{\beta-1}$.  
		\end{prop}
		With these results at hand, we continue with the proof of Theorem \ref{AA2} (iii). Since $3/2+a>2$, one has that $x^2 u\in C([-T,T]:L^2(\mathbb{R}))$. Hence, we use the integral formulation of $u$ times $x^2$, that is,
		\begin{equation}
			\begin{aligned}
				x^2 u(x,t)=x^2 e^{t\partial_x |D_x|^{a}}u_0-\int_0^t x^2 e^{(t-\tau)\partial_x |D_x|^{a}}(u\partial_x u)(\tau)\, d\tau.  
			\end{aligned}    
		\end{equation}
		The strategy is the following. Applying $\langle x\rangle^{a-1/2}J^{-5}$ to the integral equation above, we first show that the Duhamel term has the required decay. Therefore, the decay of $u(t_2)$ is equivalent to the decay of the linear part $x^{2}e^{t_2\partial_x|D_x|^{a}}u_0$, which can be computed explicitly by means of Proposition \ref{propUCPtwo1}.
		
		We first claim that
		\begin{equation}\label{integraleqpart}
			\begin{aligned}
				\Big\|\langle x \rangle^{a-1/2}J^{-5}&\int_0^t x^2 e^{(t-\tau)\partial_x |D_x|^{a}}(u\partial_x u)(\tau)\, d\tau \Big\|_{L^{\infty}_TL^2_x}\leq c_{a} \langle T\rangle^{a+5/2} \sup_{\tau\in [-T,T]}\triple u(\tau)\triple,    
			\end{aligned}
		\end{equation}
		where for simplicity, we define the norm $\triple\cdot\triple$ by
		\begin{equation*}
			\triple f \triple=\|J^{a(3/2+a)^{-}} f\|_{2}+\|\,|x|^{(3/2+a)^{-}}f\|_{2}.
		\end{equation*}
		We apply Proposition \ref{propUCPtwo1} to get
		\begin{equation*}
			\begin{aligned}
				\|\langle x\rangle^{a-1/2}J^{-5}\big(x^2 e^{(t-\tau)\partial_x |D_x|^{a}}(u\partial_x u)\big)\|_{2}\leq &c \langle t-\tau\rangle^2\big(\|\langle x\rangle^{a-1/2}J^{-5}e^{(t-\tau)\partial_x |D_x|^{a}}x|D_x|^{a}(u\partial_x u)\|_{2} \\
				&+\|\langle x\rangle^{a-1/2}J^{-5}e^{(t-\tau)\partial_x |D_x|^{a}} |D_x|^{2a}(u\partial_x u)\|_{2}\\
				&+\|\langle x\rangle^{a-1/2}J^{-5}e^{(t-\tau)\partial_x |D_x|^{a}}|D_x|^{a}(x u\partial_x u)\|_{2}\\
				&+\|\langle x\rangle^{a-1/2}J^{-5}e^{(t-\tau)\partial_x |D_x|^{a}}(x^2 u\partial_x u)\|_{2}\big)\\
				=:&\langle t-\tau\rangle^2\big(\mathcal{E}_1+\mathcal{E}_2+\mathcal{E}_3+\mathcal{E}_4\big).
			\end{aligned}    
		\end{equation*}
		Let us estimate the factors $\mathcal{E}_j$, $j=1,\dots,4$. To control $\mathcal{E}_1$, we write
		\begin{equation*}
			\begin{aligned}
				x|D_x|^{a}(u\partial_x u)=\frac{1}{2}[x;\partial_x|D_x|^{a}](u^2)+\frac{1}{2}\partial_x|D_x|^{a}(xu^2) =-\frac{1+a}{2}|D_x|^{a}(u^2)+\frac{1}{2}\partial_x|D_x|^{a}(xu^2).
			\end{aligned}    
		\end{equation*}
		Applying Propositions \ref{propUCPtwo3} and \ref{propUCPtwo2} we have
		\begin{equation*}
			\begin{aligned}
				\mathcal{E}_1\leq& c \|\langle x\rangle^{a-1/2}J^{-3}|D_x|^{a}e^{(t-\tau)\partial_x |D_x|^{a}}J^{-2}(u^2)\|_{2}+c\|\langle x\rangle^{a-1/2}J^{-3}\partial_x|D_x|^{a}e^{(t-\tau)\partial_x |D_x|^{a}}J^{-2}(x u^2)\|_{2}\\
				\leq & c \langle t-\tau \rangle^{a-1/2}\big(\|\langle x\rangle^{a-1/2}J^{-2}(u^2)\|_{2}+\||D_x|^{a(a-1/2)}J^{-2}(u^2)\|_{2}\\
				&+\|\langle x\rangle^{a-1/2}J^{-2}(xu^2)\|_{2}+\||D_x|^{a(a-1/2)}J^{-2}(xu^2)\|_{2}\big)\\
				\leq & c  \langle t-\tau \rangle^{a-1/2}\|\langle x\rangle^{a+1/2}(u^2)\|_{2}\\
				\leq & c  \langle t-\tau \rangle^{a-1/2}\|u\|_{\infty}\|\langle x\rangle^{a+1/2} u\|_{2}\leq c \langle t-\tau \rangle^{a-1/2}\triple u\triple^2.
			\end{aligned}    
		\end{equation*}
		Once again, by Propositions \ref{propUCPtwo3} and \ref{propUCPtwo2}, and fractional Leibniz rule, we get
		\begin{equation*}
			\begin{aligned}
				\mathcal{E}_2\leq& c \|\langle x\rangle^{a-1/2}J^{-5}\partial_x |D_x|^{2a}e^{(t-\tau)\partial_x |D_x|^{a}} (u^2)\|_{2}\\
				\leq & c\|\langle x\rangle^{a-1/2}e^{(t-\tau)\partial_x |D_x|^{a}} (u^2)\|_{2}\\
				\leq & c \langle t-\tau \rangle^{a-1/2}\big(\|\langle x\rangle^{a-1/2}(u^2)\|_{2}+\||D_x|^{a(a-1/2)}(u^2)\|_{2}\big)\\
				\leq & c  \langle t-\tau \rangle^{a-1/2}\|u\|_{\infty}\big(\|\langle x\rangle^{a+1/2} u\|_{2}+\||D_x|^{a(a-1/2)}u\|_{2}\big)\\
				\leq & c \langle t-\tau \rangle^{a-1/2}\triple u\triple^2.
			\end{aligned}    
		\end{equation*}
		Writing $xu \partial_x u=\frac{1}{2}\partial_x (xu^2)-\frac{1}{2}u^2$, the estimate for $\mathcal{E}_3$ is similar to that of $\mathcal{E}_1$. To estimate $\mathcal{E}_4$, we write $x^2u \partial_x u=-xu^2+\frac{1}{2}\partial_x (x^2u^2)$ 
		\begin{equation*}
			\begin{aligned}
				\mathcal{E}_4\leq& c \|\langle x\rangle^{a-1/2}J^{-3} e^{(t-\tau)\partial_x |D_x|^{a}}J^{-2}(x u^2)\|_{2}+c\|\langle x\rangle^{a-1/2}J^{-3}\partial_x e^{(t-\tau)\partial_x |D_x|^{a}}J^{-2}(x^2 u^2)\|_{2}\\
				\leq & c \langle t-\tau \rangle^{a-1/2}\big(\|\langle x\rangle^{a-1/2}J^{-2}(x u^2)\|_{2}+\||D_x|^{a(a-1/2)}J^{-2}(x u^2)\|_{2}\\
				&+\|\langle x\rangle^{a-1/2}J^{-2}(x^2u^2)\|_{2}+\||D_x|^{a(a-1/2)}J^{-2}(x^2 u^2)\|_{2}\big)\\
				\leq & c  \langle t-\tau \rangle^{a-1/2}\|\langle x\rangle^{a+3/2}(u^2)\|_{2}.
			\end{aligned}    
		\end{equation*}
		Now, using Sobolev embedding $H^{1/4}(\mathbb{R})\hookrightarrow L^4(\mathbb{R})$, and complex interpolation \eqref{inter1x},  we have
		\begin{equation*}
			\begin{aligned}
				\|\langle x\rangle^{a+3/2}(u^2)\|_{2}\leq c \|\langle x\rangle^{a/2+3/4}u\|_{4}^2  \leq c & \|J^{1/4}(\langle x\rangle^{a/2+3/4}u)\|_{2}^2  \\
				\leq c & \|J^{a(3/2+a)^{-}}u\|_{2}^2+c\|\langle x \rangle^{(3/2+a)^{-}}u\|_{2}^2.
			\end{aligned}    
		\end{equation*}
		This completes the deduction of \eqref{integraleqpart}.
		
		Note that since $u_0,u(t_2)\in L^2(|x|^{3+2a}\, dx)$, Proposition \ref{propUCPtwo3} establishes that $\langle x\rangle^{a-1/2}J^{-5}(x^2 u(t_2))\in L^2(\mathbb{R})$. Consequently, by \eqref{integraleqpart}, we have that 
		\begin{equation*}
			\begin{aligned}
				\langle x\rangle^{a-1/2}\big(J^{-5} x^2 e^{t_2\partial_x |D_x|^{a}}u_0\big)\in L^2(\mathbb{R}).   
			\end{aligned}    
		\end{equation*}
		It remains to analyze the linear part of the integral equation.
		
		By Proposition \ref{propUCPtwo1}, we write
		\begin{equation*}
			\begin{aligned}
				\langle x\rangle^{a-1/2}&\big(J^{-5} x^2 e^{t_2\partial_x|D_x|^{a}}u_0\big)\\
				=&-(1+a) t_2 \langle x\rangle^{a-1/2}\big(J^{-5}e^{t_2\partial_x |D_x|^{a}}(x|D_x|^{a} u_0)\big)\\
				&+(1+a)^2 t_2^2 \langle x\rangle^{a-1/2}\big(J^{-5}e^{t_2\partial_x |D_x|^{a}}|D_x|^{2a}u_0\big)\\
				&-(1+a) t_2 \langle x\rangle^{a-1/2}\big(J^{-5} e^{t_2\partial_x |D_x|^{a}}|D_x|^{a}(xu_0)\big)+\langle x\rangle^{a-1/2}\big(J^{-5}e^{t_2\partial_x |D_x|^{a}}(x^2u_0)\big)\\
				=:&\widetilde{\mathcal{C}_1}+\widetilde{\mathcal{C}_2}+\widetilde{\mathcal{C}_3}+\widetilde{\mathcal{C}_4}.
			\end{aligned} 
		\end{equation*}
		From Propositions \ref{propUCPtwo3} and \ref{propUCPtwo2}, it is seen that
		\begin{equation*}
			\begin{aligned}
				\|\widetilde{\mathcal{C}_2}\|_{2}+\|\widetilde{\mathcal{C}_3}\|_{2}\leq & c\langle t_2 \rangle^{2}\big( \|\langle x\rangle^{a-1/2}e^{t_2\partial_x |D_x|^{a}}J^{-2}u_0\|_{2}+\|\langle x\rangle^{a-1/2}e^{t_2\partial_x |D_x|^{a}}J^{-2}(xu_0)\|_{2}\big) \\
				\leq & c \langle t_2 \rangle^{3/2+a}\big( \|\langle x\rangle^{a-1/2}J^{-2}u_0\|_{2}+\||D_x|^{a(a-1/2)}J^{-2}u_0\|_{2}+\|\langle x\rangle^{a-1/2}J^{-2}(xu_0)\|_{2}\\
				&+\||D_x|^{a(a-1/2)}J^{-2}(xu_0)\|_{2}\big)\\
				\leq  & c \langle t_2 \rangle^{3/2+a} \|\langle x \rangle^{a+1/2} u_0\|_{2}.
			\end{aligned}   
		\end{equation*}
		Similarly, we get
		\begin{equation*}
			\begin{aligned}
				\|\widetilde{\mathcal{C}_4}\|_{2}\leq c & \|\langle x\rangle^{a-1/2}e^{t_2\partial_x |D_x|^{a}}J^{-2}(x^2u_0)\|_{2} \\
				\leq c & \langle t_2 \rangle^{a-1/2}\big(\|\langle x\rangle^{a-1/2}J^{-2}(x^2u_0)\|_{2}+\||D_x|^{a(a-1/2)}J^{-2}(x^2u_0)\|_{2}\big)\\
				\leq c &\langle t_2 \rangle^{a-1/2}\|\langle x \rangle^{3/2+a} u_0\|_{2}.
			\end{aligned}   
		\end{equation*}
		We conclude that 
		\begin{equation*}
			\langle x\rangle^{a-1/2}J^{-5}(x^2 u(t_2))\in L^2(\mathbb{R})\, \text{ if and only if }\, \widetilde{\mathcal{C}_1}\in L^2(\mathbb{R}).   
		\end{equation*}
		Now, using that $[x;|D_x|^{a}]=a\partial_x |D_x|^{a-2}$, we write
		\begin{equation*}
			\begin{aligned}
				\widetilde{\mathcal{C}_1}=&-(1+a) t_2 \langle x\rangle^{a-1/2}\big(J^{-5}(e^{t_2\partial_x |D_x|^{a}}-1)(x|D_x|^{a} u_0)+a J^{-5}(\partial_x|D_x|^{a-2} u_0)+J^{-5}|D_x|^{a}(x u_0)\big)\\
				=&-(1+a) t_2 \langle x\rangle^{a-1/2}\big( a J^{-5}(e^{t_2\partial_x |D_x|^{a}}-1)\partial_x|D_x|^{a-2} u_0\\
				&+J^{-5}(e^{t_2\partial_x |D_x|^{a}}-1)|D_x|^{a}( xu_0)+a J^{-5}(\partial_x|D_x|^{a-2} u_0)+J^{-5}|D_x|^{a}(x u_0)\big)\\
				=:&-(1+a) t_2\big(\widetilde{\mathcal{C}_{11}}+\widetilde{\mathcal{C}_{12}}+\widetilde{\mathcal{C}_{13}}+\widetilde{\mathcal{C}_{14}}\big). 
			\end{aligned}    
		\end{equation*}
		By mean value inequality, $\langle \xi\rangle^{-5}(e^{it_2\xi|\xi|^{a}}-1)\xi|\xi|^{a-2}\in W^{1,\infty}(\mathbb{R})$, then we apply Theorem \ref{TheoSteDer} and \eqref{prelimneq} to deduce
		\begin{equation*}
			\begin{aligned}
				\|\widetilde{\mathcal{C}_{11}}\|_{2}\leq  & c\|J^{a-1/2}\big(\langle \xi\rangle^{-5}(e^{it_2\xi|\xi|^{a}}-1)\xi|\xi|^{a-2}\widehat{u_0}\big)\|_{2} \\
				\leq & c\|\langle \xi\rangle^{-5}(e^{it_2\xi|\xi|^{a}}-1)\xi|\xi|^{a-2}\widehat{u_0}\|_{2}+c\|\mathcal{D}^{a-1/2}\big(\langle \xi\rangle^{-5}(e^{it_2\xi|\xi|^{a}}-1)\xi|\xi|^{a-2}\big)\widehat{u_0}\|_{2}\\
				&+c\|\langle \xi\rangle^{-5}(e^{it_2\xi|\xi|^{a}}-1)\xi|\xi|^{a-2}\mathcal{D}^{a-1/2}\widehat{u_0}\|_{2}\\
				\leq & c\langle t_2 \rangle \big(\|\widehat{u_0}\|_{2}+\|\mathcal{D}^{a-1/2}\widehat{u_0}\|_{2}\big)\\
				\leq & c \langle t_2 \rangle\|\langle x \rangle^{a-1/2}u_0\|_{2}.
			\end{aligned}   
		\end{equation*}
		By similar arguments used above, we get
		\begin{equation*}
			\begin{aligned}
				\|\widetilde{\mathcal{C}_{12}}\|_{2}+\|\widetilde{\mathcal{C}_{14}}\|_{2}\leq &  c\langle t_2 \rangle\|\langle x \rangle^{a+1/2}u_0\|_{2}.
			\end{aligned}   
		\end{equation*}
		Consequently, it must follow that $\widetilde{\mathcal{C}_{13}}\in L^2(\mathbb{R})$, i.e., we have deduced that
		\begin{equation*}
			\langle x\rangle^{a-1/2}J^{-5}(x^2 u(t_2))\in L^2(\mathbb{R})\, \text{ if and only if }\, \langle x \rangle^{a-1/2}J^{-5}(\partial_x |D_x|^{a-2}u_0)\in L^2(\mathbb{R}).   
		\end{equation*}
		Therefore, the problem is reduced to the study of the behavior of $\widehat{u_0}$ near the origin.
		
		Using Plancherel's identity, we write $J^{-5}(\partial_x |D_x|^{a-2}u_0)$ as
		\begin{equation*}
			\begin{aligned}
				\mathcal{F}(\langle x \rangle^{a-1/2}J^{-5}(\partial_x |D_x|^{a-2}u_0))(\xi)=&iJ_{\xi}^{a-1/2}\big(\langle \xi \rangle^{-5}\xi|\xi|^{a-2}\chi\widehat{u_0}(\xi)+\langle \xi \rangle^{-5}\xi|\xi|^{a-2}(1-\chi)\widehat{u_0}(\xi)\big) \\
				=&iJ_{\xi}^{a-1/2}\big(\langle \xi \rangle^{-5}\xi|\xi|^{a-2}\chi \widehat{u_0}(0)+\langle \xi \rangle^{-5}\xi|\xi|^{a-2}\chi(\widehat{u_0}(\xi)-\widehat{u_0}(0))\\
				&+\langle \xi \rangle^{-5}\xi|\xi|^{a-2}(1-\chi)\widehat{u_0}(\xi)\big).
			\end{aligned}   
		\end{equation*}
		Since $\langle \xi \rangle^{-5}\xi|\xi|^{a-2}(1-\chi)\in W^{1,\infty}(\mathbb{R})$, using \eqref{prelimneq}, it is not hard to deduce that
		\begin{equation*}
			\|J_{\xi}^{a-1/2}\big(\langle \xi \rangle^{-5}\xi|\xi|^{a-2}(1-\chi)\widehat{u_0}\big)\|_{2}\leq c \|\langle x \rangle^{a-1/2}u_0\|_{2}.  
		\end{equation*}
		Next, we write
		\begin{equation*}
			\begin{aligned}
				\langle \xi \rangle^{-5}\xi|\xi|^{a-2}\chi(\widehat{u_0}(\xi)-\widehat{u_0}(0))=\langle \xi \rangle^{-5}\xi^2|\xi|^{a-2}\chi\int_0^1 \frac{d}{d\xi}\widehat{u}(\sigma\xi)\, d\sigma.    
			\end{aligned}   
		\end{equation*}
		By the properties of the fractional derivative $\mathcal{D}^{a-1/2}$ (see \eqref{prelimneq}), we get
		\begin{equation*}
			\begin{aligned}
				\|J_{\xi}^{a-1/2}&\big(\langle \xi \rangle^{-5}\xi|\xi|^{a-2}\chi(\widehat{u_0}(\xi)-\widehat{u_0}(0))\big)\|_{2}\\
				\leq & c \|\langle \xi \rangle^{-5}\xi|\xi|^{a-2}\chi(\widehat{u_0}(\xi)-\widehat{u_0}(0))\|_{2}+c\|\mathcal{D}^{a-1/2}\big(\langle \xi \rangle^{-5}\xi^2|\xi|^{a-2} \chi\big) \int_0^1 \frac{d}{d\xi}\widehat{u_0}(\sigma\xi)\, d\sigma  \|_{2}\\
				&+c\|\langle \xi \rangle^{-5}\xi^2|\xi|^{a-2} \chi\|_{\infty}\||D_{\xi}|^{a-1/2}\big( \int_0^1 \frac{d}{d\xi}\widehat{u_0}(\sigma\xi)\, d\sigma \big) \|_{2}\\
				\leq & c\|\widehat{u_0}\|_{\infty}+c\int_0^1\|\frac{d}{d\xi}\widehat{u_0}(\sigma \xi)\|_{2}\, d\sigma+c\int_0^1\sigma^{a-1/2}\||D_{\xi}|^{a-1/2}\frac{d}{d\xi}\widehat{u_0}(\sigma \xi)\|_{2}\, d\sigma\\
				\leq & c\|\langle x \rangle^{1/2^{+}}u_0\|_{2}+c\int_0^1\sigma^{-1/2}\|xu_0\|_{2}\, d\sigma+c\int_0^1\sigma^{a-1}\||x|^{a+1/2}u_0\|_{2}\, d\sigma\\
				\leq & c\|\langle x \rangle^{a+1/2}u_0\|_{2}.
			\end{aligned}
		\end{equation*}
		We conclude that
		\begin{equation*}
			\langle x\rangle^{a-1/2}J^{-5}(x^2 u(t_2))\in L^2(\mathbb{R})\, \text{ if and only if }\, J_{\xi}^{a-1/2}\big(\langle \xi\rangle^{-5}\xi|\xi|^{a-2}\chi \widehat{u_0}(0)\big)\in L^2(\mathbb{R}),
		\end{equation*}
		which by the fact that $\langle \xi\rangle^{-5}\xi|\xi|^{a-2}\chi \widehat{u_0}(0)\in L^2(\mathbb{R})$ and Theorem \ref{TheoSteDer} is equivalent to
		\begin{equation}\label{ConclusionUCP}
			\langle x\rangle^{a-1/2}J^{-5}(x^2 u(t_2))\in L^2(\mathbb{R})\, \text{ if and only if }\, \widehat{u_0}(0)\mathcal{D}^{a-1/2}\big(\langle \xi\rangle^{-5}(\xi|\xi|^{a-2}\chi )\in L^2(\mathbb{R}).
		\end{equation}
		One can prove that $\mathcal{D}^{a-1/2}\big(\langle \xi\rangle^{-5}\xi|\xi|^{a-2}\chi \big)\notin L^2(\mathbb{R})$. As a consequence, from \eqref{ConclusionUCP}, it must be the case that $\widehat{u_0}(0)=0$. The proof is complete. 
	\end{proof}	
		\section*{Use of AI-assisted tools}
	
	The figures in this paper were produced with the assistance of Claude Opus 5.5 (model identifier \texttt{claude-opus-5-5}), a large language model developed by Anthropic. The model was used to generate the TikZ and \texttt{pgfplots} code for the figures, based on the mathematical content provided by the authors. All figures were checked by the authors for mathematical accuracy and agreement with the results of the paper. The authors take full responsibility for the content of this work.

\end{document}